\documentclass[12pt,reqno,letterpaper]{amsart}
\usepackage[T1]{fontenc}
\usepackage[utf8]{inputenc}
\usepackage{lmodern}
\usepackage[margin=2.5cm]{geometry}
\usepackage{amsmath,amssymb,amsthm,mathtools,bm}
\usepackage{array,booktabs,longtable}
\usepackage{graphicx,tikz}
\usepackage{enumitem}
\usepackage[colorlinks=true,linkcolor=blue,citecolor=blue,urlcolor=blue,unicode]{hyperref}
\usepackage[switch,mathlines]{lineno}
\makeatletter
\@ifundefined{subjclassname@2020}{%
  \@namedef{subjclassname@2020}{\textup{2020} Mathematics Subject Classification}%
}{}
\makeatother

\numberwithin{equation}{section}
\newtheorem{theorem}{Theorem}[section]
\newtheorem{lemma}[theorem]{Lemma}
\newtheorem{proposition}[theorem]{Proposition}
\newtheorem{corollary}[theorem]{Corollary}
\theoremstyle{definition}
\newtheorem{definition}[theorem]{Definition}
\theoremstyle{remark}

\newif\ifshowlineno
\newcommand{\srcline}{\marginpar{\raisebox{0pt}{\tiny\ttfamily\the\inputlineno}}}
\allowdisplaybreaks[2]
\title[Finiteness of planar five-body central configurations]{Finiteness of planar central configurations
of five bodies with positive masses}
\author{Kefeng Liu }
\address[Kefeng Liu]{Mathematical Sciences Research Center, Chongqing University of Technology, Chongqing, 400054, China}
\email{kefengliu@cqut.edu.cn}
\author{Yangshanshan Liu}
\address[Yangshanshan Liu]{Mathematical Sciences Research Center, Chongqing University of Technology, Chongqing, 400054, China}
\email{lyss6133@126.com}
\date{October 20, 2026}
\hypersetup{pdfauthor={Kefeng Liu and Yangshanshan Liu}}
\begin{document}
\begin{abstract}
We prove that, for every choice of positive masses, the planar Newtonian five-body problem has
finitely many central configurations up to similarity. After normalization and removal of the
translation and rotation freedoms, every physical configuration is isolated in the complex solution
variety. The proof combines the classification of singular sequences by Albouy and Kaloshin with
a lower bound for the potential, asymptotic estimates, and local analytic reductions.
\end{abstract}
\maketitle
\nolinenumbers
\tableofcontents

\ifshowlineno\linenumbers\fi
\section{Introduction}\label{sec:1}

The Newtonian $n$-body problem describes the motion of $n$ point masses under their mutual gravitational attraction. We consider the planar problem with $n\ge2$ and fixed masses $m_1,\ldots,m_n>0$. If $q_i(t)\in\mathbb{R}^2$ is the position of body $i$ at time $t$, then, with the gravitational constant set to one, the equations of motion are 
\[
\ddot q_i=\sum_{j\ne i}m_j\frac{q_j-q_i}{|q_j-q_i|^3},
\qquad i=1,\ldots,n.
\]
These equations are defined away from collisions, where $q_i\ne q_j$ for $i\ne j$. Throughout this paper, planar configurations include collinear ones.

For a collision-free configuration $q=(q_1,\ldots,q_n)$, denote the total mass and the center of mass by
\[
M=\sum_{i=1}^n m_i,\qquad C=\frac{1}{M}\sum_{i=1}^n m_iq_i.
\]
A configuration is called a \emph{central configuration} if there exists a common constant $\lambda>0$ such that
\begin{equation}\label{eq:1.1}
\sum\limits_{j\ne i}m_j\frac{q_i-q_j}{|q_i-q_j|^3}
=\lambda(q_i-C),\qquad i=1,\ldots,n.
\end{equation}
Thus the Newtonian acceleration of each body is $-\lambda(q_i-C)$, with the same proportionality factor for all bodies.

Central configurations generate homographic solutions, in which the configuration remains similar to itself throughout the motion. Releasing the bodies from rest gives a homothetic collapse to their center of mass. A suitable choice of initial velocities produces a rigid rotation with angular speed $\sqrt{\lambda}$, called a relative equilibrium. Equivalently, on the space of configurations with center of mass at the origin and a fixed positive moment of inertia, the central configurations are precisely the critical points of the restricted Newtonian potential \cite[Theorem~C]{Smale70}.

Translations, rotations, reflections, and positive scalings preserve central configurations. We therefore count labeled configurations up to these similarities. Translation and positive scaling give a representative with $C=0$ and $\lambda=1$, unique up to rotation and reflection. The finiteness problem asks whether, for every fixed $n$ and every choice of positive masses, the number of such similarity classes is finite. The question arose in the work of Chazy and Wintner \cite{Chazy,Wintner} and became Smale's sixth problem for the twenty-first century \cite{Smale98}. Finiteness does not require every central configuration to be nondegenerate; the stronger nondegeneracy assertion already fails for positive masses \cite[p.~491]{Palmore}. The present paper addresses the case $n=5$.

Euler and Lagrange determined the three-body central configurations, and Moulton proved uniqueness of the collinear configuration for each ordering of arbitrary positive masses \cite{Moulton,Smale70}. Hampton and Moeckel established finiteness for all positive four-body masses \cite{HM}. In higher dimensions, Moeckel proved generic finiteness for Dziobek configurations, whose affine span has dimension $n-2$ \cite{Moeckel01}; Hampton and Jensen refined the spatial five-body result by giving explicit exceptional mass conditions \cite{HJ}. For the planar six-body problem, Chang and Chen developed symbolic algorithms to enumerate possible singular diagrams and asymptotic orders \cite{ChangChenI}, and to derive necessary mass relations for many of these diagrams \cite{ChangChenII}.

A related finiteness problem concerns stationary configurations of planar point vortices, including relative equilibria and self-similar collapse configurations. Yu \cite{Yu23} adapted the singularity analysis of Albouy and Kaloshin to the four-vortex problem and proved finiteness of normalized stationary configurations for arbitrary nonzero real vorticities, with the nonzero complex motion parameter fixed. For five vortices, Yu and Zhu \cite{YZ25} established finiteness results for relative equilibria and collapse configurations outside a specified codimension-two subvariety of vorticity space. In particular, they proved that the number of similarity classes of stationary configurations is finite whenever all vorticities have the same sign.

Within the planar five-body problem, Hampton obtained finiteness results for symmetric kites subject to explicit mass exceptions \cite{HamptonKite}. Lee and Santoprete classified the isolated configurations with equal masses, without thereby excluding positive-dimensional families \cite{LS}. Moczurad and Zgliczy\'{n}ski subsequently gave a computer-assisted proof of the complete classification for five, six, and seven equal masses, and proved that every such configuration has a reflection symmetry \cite{MZ19}. These results concern restricted shapes or masses. For unrestricted planar configurations, Albouy and Kaloshin \cite{AK} (hereafter AK) proved the following theorem.

\begin{theorem}[{Albouy--Kaloshin \cite[Theorem~6]{AK}}]\label{thm:AK}
There exists a closed algebraic exceptional set $\mathcal{A}$ of codimension two such that, for every mass vector $(m_1,\ldots,m_5)\in(\mathbb{R}_{>0})^5\setminus\mathcal{A}$, there are only finitely many normalized central configurations of the planar five-body problem.
\end{theorem}

Here normalized central configurations are the complex solutions in \cite[Definition~2]{AK}, with multiplier one and the rotation gauge imposed, including all signed-distance choices. Restricting this statement to the physical branch gives \cite[Theorem~2]{AK}. Their method follows a hypothetical complex family to a singular limit and classifies its possible degeneration diagrams.

The result below establishes physical finiteness for every positive mass vector, including those in $\mathcal{A}$. It also establishes that each physical configuration is isolated in the complex solution variety after normalization and removal of the translation and rotation freedoms.

\begin{theorem}[Main theorem]\label{thm:1.1}
For every $m_1,\ldots,m_5>0$, the physical planar Newtonian
central configurations have finitely many similarity classes.
After normalizing the multiplier to one, fixing the center of
mass at the origin, and imposing the rotation gauge, every
physical configuration is an isolated point of the resulting
complex solution variety.
\end{theorem}

The proof is by contradiction on a positive-dimensional irreducible complex component through a physical point. The common constant $U=I$ satisfies the physical lower bound \eqref{eq:3.1}. The construction in the proof of \cite[Theorem~6]{AK} gives a singular sequence along which some 4-product tends to infinity. By \cite[Proposition~8]{AK}, its balanced limit has one of the first eight diagrams. These steps precede the generic-mass exclusions in AK.

The physical lower bound excludes diagram~3 through an upper bound for $U$ in the limit, and diagrams~1 and~4 through estimates of the limiting binary contribution to $I$ (Sections~\ref{sec:3} and~\ref{sec:4}). For diagram~7, the divergent 4-product and torque identities control the mixed term in $I$ and force $I\to0$, contradicting its positive constant value (Section~\ref{sec:8}).

Diagram~6 is excluded in Section~\ref{sec:9} by reducing the function $G=U+I/2$ on centered configurations to one parameter; the derivative of the reduced function has a nonzero leading term, contradicting the remaining critical equation. In diagram~2, Lagrange limits violate the physical lower bound, whereas non-Lagrange Euler limits are ruled out using the limiting constraint on the crossed terms and a partial Hessian argument, with a separate treatment of the degenerate case (Section~\ref{sec:11}). The analysis of diagram~5 combines a classification of scales with estimates of $U$ and $I$, leading force balances, and Schur-complement estimates. For the remaining collinear limit, polynomial elimination yields a contradiction with the physical lower bound (Section~\ref{sec:5}). For diagram~8, a two-parameter analytic reduction forces paired masses, after which comparison of an upper bound for $I=U$ with the physical lower bound excludes the residual case (Section~\ref{sec:10}).

Section~\ref{sec:12} completes the proof of Theorem~\ref{thm:1.1}. No positional symmetry or nondegeneracy of every central configuration is assumed. Roberts constructed a continuum of five-body relative equilibria with one negative mass \cite{Roberts}; AK describe a positive-mass signed modification on a nonphysical distance sheet \cite{AK}. The common value of $U$ and $I$ on the latter component lies below the physical lower bound. Thus the physical hypothesis is essential, and the theorem does not assert that every complex signed fiber is zero-dimensional.

An initial draft of the material in Sections 3–9 was completed .

\section{Basic setup}\label{sec:2}

Write $q_{ij}=q_i-q_j$. We fix the following body-index sets throughout the paper:
\begin{equation}
\begin{aligned}
S_2&=\{1,2\},& T_2&=\{3,4\},\\
S_3&=\{3,4,5\},& T_3&=\{1,2,5\},\\
S_4&=\{1,2,3,4\},& T_4&=\{1,3,4,5\},\\
S_5&=\{1,2,3,4,5\}.
\end{aligned}
\label{eq:2.index_sets}
\end{equation}
After any stated relabeling of bodies, these sets refer to the resulting labels. Arbitrary subsets and variable partitions will be specified separately.
For an index set $S\subseteq S_5$, let
\[
e_k(S)=\sum\limits_{i_1<\cdots<i_k\atop\ i_1,\dots,i_k\in S}m_{i_1}\cdots m_{i_k}
\]
denote the $k$-th elementary symmetric function of the masses indexed by $S$; in particular
\[
e_1(S)=\sum\limits_{i\in S}m_i,\qquad
e_2(S)=\sum\limits_{i<j\atop i,j\in S}m_im_j .
\]
The total mass is $M=e_1(S_5)$ and the pair sum of the whole system is $e_2(S_5)$.

For every nonempty $S\subseteq S_5$, we call the bodies indexed by $S$ a cluster and define its center of mass by
\[
C_S=\frac{1}{e_1(S)}\sum\limits_{i\in S}m_iq_i.
\]
Thus $C=C_{S_5}$ in \eqref{eq:1.1}. In the complex coordinates introduced below, $C_{S,z}$ and $C_{S,w}$ denote the $z$ and $w$ components of $C_S$.
For every nonempty $S\subseteq S_5$ and every $i\in S_5$, define the coordinates relative to the center of mass of $S$ by
\begin{equation}\label{eq:2.relative_coordinates}
q_{S,i}:=q_i-C_S,\qquad \zeta_{S,i}:=z_i-C_{S,z},\qquad \omega_{S,i}:=w_i-C_{S,w}.
\end{equation}
Thus $q_{S,i}=(\zeta_{S,i},\omega_{S,i})$ in $(z,w)$ components. The first subscript specifies the reference cluster, and $i$ may lie inside or outside it.
In particular, the positions of body 5 relative to the centers of mass of the two pairs are
\[
q_{S_2,5}=q_5-C_{S_2},\qquad q_{T_2,5}=q_5-C_{T_2}.
\]
Their components are $(\zeta_{S_2,5},\omega_{S_2,5})$ and $(\zeta_{T_2,5},\omega_{T_2,5})$, respectively. The differences $C_{S_2}-q_5$ and $C_{T_2}-q_5$ have the opposite sign.

Hats denote quantities obtained by scaling or normalization; each local chart specifies the scaling. A superscript $0$ denotes a limiting value or evaluation at a specified limiting critical point. For example, $\widehat{\zeta}_{S_3,i}$ and $\widehat{\zeta}_{S_3,i}^{0}$ denote a rescaled coordinate and its limit. Powers of a limiting value are written with parentheses, as in $(r_{ij}^{0})^2$. A tilde denotes a function obtained by solving some critical equations and substituting the solution.

In expansions of $G$, the symbols $G_{(\mathfrak{a},i)}$ and $G_{(\mathfrak{c},i)}$ denote the coefficients of $\mathfrak{a}^i$ and $\mathfrak{c}^i$, respectively. For the ordered parameter pair $(\mathfrak{b}_1,\mathfrak{b}_2)$, $G_{\mathfrak{b}}^{i,j}$ denotes the coefficient of $\mathfrak{b}_1^i\mathfrak{b}_2^j$ in the analytic part specified below. Tildes distinguish coefficients of the reduced function. The parameters and any evaluation at a limiting critical point are specified with each expansion.

\subsection{The complex central configuration equations}\label{sec:2.1}

\emph{Physical configurations and normalization.}
A collision-free central configuration is called \emph{physical} if its Cartesian positions lie in $\mathbb{R}^2$ and its signed distances are the positive Euclidean distances $r_{ij}=|q_i-q_j|>0$ for every distinct pair. This condition does not specify a scaling or rotation normalization. Write
\[
U=\sum\limits_{i<j}\frac{m_im_j}{r_{ij}},\qquad I=\sum\limits_i m_iq_i\cdot q_i,
\]
after placing the center of mass at the origin. For a physical configuration, pairwise summation of \eqref{eq:1.1} gives $\lambda I=U>0$. Dilating positions by $\lambda^{1/3}$ makes the multiplier one, uniquely within its positive scaling class. Common positive scaling of all masses preserves central-configuration shapes. We retain arbitrary positive masses. In this paper, \emph{multiplier-one}, or \emph{normalized} when no rotation condition is stated, refers to the centered equations with $\lambda=1$. The rotation normalization is specified separately.

Throughout the paper, $\mathrm{i}$ denotes the imaginary unit, with $\mathrm{i}^2=-1$. In complex coordinates $z_i=x_i+\mathrm{i}y_i$, $w_i=x_i-\mathrm{i}y_i$, regarded henceforth as independent, and writing $\delta_{ij}=r_{ij}^{-3}$ for the inverse cubes, the squared distance and force equations are
\begin{equation}\label{eq:2.1}
r_{ij}^2=(z_i-z_j)(w_i-w_j),
\end{equation}
\begin{equation}\label{eq:2.2}
z_i=\sum\limits_{j\ne i}m_j(z_i-z_j)\delta_{ij},\qquad
w_i=\sum\limits_{j\ne i}m_j(w_i-w_j)\delta_{ij}.
\end{equation}
Every $r_{ij}$ is nonzero, and introducing $\delta_{ij}$ as a variable makes this an affine polynomial system without collision points. Its center of mass is at the origin, either as a consequence of the force equations or as an explicitly imposed equation.

In these coordinates, Cartesian reality is equivalent to $w_i=\overline{z_i}$ for all $i\in S_5$. A \emph{physical point} of the complex solution variety satisfies these relations and $r_{ij}>0$ for every distinct pair. A component containing a physical point need not consist of physical points.

\emph{Signed-distance branches.}
Following the complexification in \cite[Section~2]{AK}, retain both choices of each signed distance. On a sufficiently small simply connected neighborhood where every $(z_i-z_j)(w_i-w_j)$ is nonzero, a \emph{signed-distance branch}, also called a \emph{sheet}, consists of a holomorphic choice of $r_{ij}$ satisfying \eqref{eq:2.1} for each unordered pair. The choice is fixed by its value at one point. At a physical point the physical choice has every $r_{ij}>0$.

For distinct $i,j\in S_5$, put
\[
z_{ij}=z_i-z_j,\qquad w_{ij}=w_i-w_j,\qquad
Z_{ij}=z_{ij}\delta_{ij},\qquad W_{ij}=w_{ij}\delta_{ij}.
\]
Following the terminology of \cite[p.~542]{AK}, we call $z_{ij}$ the \emph{$z$-separation} and $w_{ij}$ the \emph{$w$-separation} of the ordered pair $(i,j)$. These are complex coordinate differences, not their moduli.

The signed distance is attached to the unordered pair, so $r_{ji}=r_{ij}$ and $\delta_{ji}=\delta_{ij}$. Hence $Z_{ji}=-Z_{ij}$ and $W_{ji}=-W_{ij}$. Our differences and pair components have the opposite sign to those with the same ordered indices in \cite{AK}. In this notation, \eqref{eq:2.2} becomes
\[
z_i=\sum\limits_{j\in S_5\setminus\{i\}}m_jZ_{ij},\qquad
w_i=\sum\limits_{j\in S_5\setminus\{i\}}m_jW_{ij},\qquad i\in S_5.
\]
For vectors $a=(a_z,a_w)$ and $b=(b_z,b_w)$, define the bilinear Euclidean product and the two-dimensional skew product by
\[
a\cdot b=\tfrac12(a_zb_w+a_wb_z),\qquad [a,b]=a_zb_w-a_wb_z.
\]
Neither product uses complex conjugation. In particular, $a^2=a\cdot a=a_za_w$ and $I=\sum\limits_i m_iz_iw_i$. The words ``distance'' and ``potential'' on this complex cover refer to the chosen signed variables, not to their moduli.

For a nonempty cluster $S\subseteq S_5$, define its moment of inertia $I_S$, potential $U_S$, and the contribution $f_i^S$ of the bodies in $S\setminus\{i\}$ to the right-hand side of the central configuration equation for body $i$ by
\[
\begin{aligned}
I_S&:=\sum\limits_{i\in S}m_iq_{S,i}\cdot q_{S,i}
=\frac{1}{e_1(S)}\sum\limits_{\substack{i<j\\i,j\in S}}m_im_jr_{ij}^2,\\
U_S&:=\sum\limits_{\substack{i<j\\i,j\in S}}\frac{m_im_j}{r_{ij}},\\
f_i^S&:=\sum\limits_{j\in S\setminus\{i\}}m_jq_{ij}\delta_{ij},\qquad i\in S.
\end{aligned}
\]
The notation $I_S$ always refers to the moment of inertia about the cluster's own center of mass.
The components of $f_i^S$ are denoted $f_{i,z}^S,f_{i,w}^S$.
The sum defining $f_i^S$ is exact; no assumption on the relative orders of its terms is part of this definition.
For the centered whole system, $U=U_{S_5}$ and $I=I_{S_5}$.

For every nonempty $S\subseteq S_5$, multiply the equations \eqref{eq:2.2} for $i\in S$ by $m_i$ and sum. Antisymmetry cancels every pair internal to $S$, giving
\[
\begin{aligned}
e_1(S)C_{S,z}&=\sum\limits_{\substack{i\in S\\j\in S_5\setminus S}}m_im_jZ_{ij},\\
e_1(S)C_{S,w}&=\sum\limits_{\substack{i\in S\\j\in S_5\setminus S}}m_im_jW_{ij}.
\end{aligned}
\]
For each $i\in S$, separating the internal and external terms in its individual equations also gives
\[
\begin{aligned}
f_{i,z}^S&=\sum\limits_{j\in S\setminus\{i\}}m_jZ_{ij}
=z_i-\sum\limits_{j\in S_5\setminus S}m_jZ_{ij},\\
f_{i,w}^S&=\sum\limits_{j\in S\setminus\{i\}}m_jW_{ij}
=w_i-\sum\limits_{j\in S_5\setminus S}m_jW_{ij}.
\end{aligned}
\]

\begin{lemma}[Constancy of $U$ and $I$]\label{lem:2.1}
On every irreducible component of the normalized complex solution variety, $U=I$ is constant. If the component contains a physical point, the common constant is positive real.
\end{lemma}
\begin{proof}
Multiply the $z$ equation of \eqref{eq:2.2} by $m_iw_i$ and sum over $i$. Each unordered pair contributes $m_im_jw_{ij}Z_{ij}=m_im_j/r_{ij}$, so $I=U$. Along the smooth locus of the distance cover,
\[
dU=-\frac12\sum\limits_i m_i(w_i\,dz_i+z_i\,dw_i)=-\tfrac12dI
\]
at a central configuration. The first equality uses both force equations and follows by differentiating \eqref{eq:2.1}; its distance derivative is legitimate because $r_{ij}\ne0$. Since $I=U$ on the solution variety, $dU=0$ on the smooth locus of each irreducible component. Over characteristic zero this makes the regular function $U$ constant on the component. Equivalently, use the connected dense smooth locus and continuation. Positivity at a physical point is immediate. This is the invariant in \cite[Lemma 2]{AK}.
\end{proof}

We first impose the rotation gauge $z_{12}=w_{12}$, possible whenever $r_{12}\ne0$. With this gauge, our centered multiplier-one solutions are the \emph{normalized central configurations} of \cite[Definition~2]{AK}; those that are physical are the \emph{positive normalized central configurations} of \cite[Definitions~1--2]{AK}. Later we use the balanced normalization of \cite[Section~3, pp.~547--548]{AK}, specifying that change when it is made. A fixed collection of squared distances determines only finitely many configurations in this gauge. Indeed, $y_1=y_2$ and $(x_1-x_2)^2=r_{12}^2\ne0$ first determine their difference up to sign. Differences of squared distances then determine each $x_j-x_1$; each $(y_j-y_1)^2$ is determined, leaving finitely many choices. The center-of-mass equations determine the remaining translation. Each squared-distance tuple also has only finitely many signed lifts.

\emph{Singular sequences.}
Following \cite[Definition~3 and the preceding construction]{AK}, consider a sequence of normalized central configurations in the rotation gauge just imposed. Form the two augmented coordinate lists
\[
\bigl((z_i)_{i\in S_5},(Z_{ij})_{i<j,\ i,j\in S_5}\bigr),\qquad
\bigl((w_i)_{i\in S_5},(W_{ij})_{i<j,\ i,j\in S_5}\bigr).
\]
Extract a subsequence so that, in each list, a component attaining the maximum modulus has a fixed index and the list divided by this maximum modulus converges. In the unbounded case, the resulting subsequence is called a \emph{singular sequence}. After further subsequence selection, at least one of the two maxima tends to infinity. The balanced normalization used below makes both maxima equal to $\epsilon^{-2}$, with $\epsilon\to0$.

Following \cite[Section~9, p.~580]{AK}, a \emph{4-product} is a quantity
\begin{equation}\label{eq:2.3}
p_{ij}=r_{ij}^2r_{kl}^2r_{lm}^2r_{mk}^2,\qquad
\{i,j,k,l,m\}=S_5,
\end{equation}
where the five indices are distinct.

\begin{lemma}[A singular sequence with a 4-product tending to infinity]\label{lem:2.2}
If a fixed positive five-body mass vector has infinitely many physical similarity classes, its normalized, rotation-gauged complex solution variety has an irreducible component containing physical points and a singular sequence on that component along which $U=I$ is a positive real constant and some 4-product \eqref{eq:2.3} tends to infinity. After the balanced normalization of \cite{AK} and further subsequence selection, the sequence approaches one of diagrams 1--8 in Figure~\ref{fig:ak-five-body} (\cite[Figure~11]{AK}).
\end{lemma}
\begin{proof}
There are finitely many irreducible components. The finite-fiber observation above implies that one component containing infinitely many physical gauge points has a nonconstant squared-distance function. A nonconstant regular function on an irreducible complex affine variety has constructible dense image in $\mathbb{C}$ \cite[Proposition~2.31, p.~37]{Mumford}, hence an image containing the complement of finitely many points. We may push that distance square to zero while staying on the component.

The opening argument of \cite[Theorem 6]{AK}, using their Proposition 9, then produces a 4-product \eqref{eq:2.3} tending to zero. It is nonzero on the collision-free variety, so it is nonconstant on the same component. Push it instead to infinity. This gives a singular sequence along which that 4-product tends to infinity. Their Proposition 8 excludes diagrams 9--16 for this sequence, leaving diagrams 1--8. These steps precede the generic-mass exclusions in their theorem; they do not assume that the masses avoid its exceptional set. Lemma~\ref{lem:2.1} preserves the physical invariant.
\end{proof}

\subsection{AK diagrams and asymptotic conventions}\label{sec:2.2}

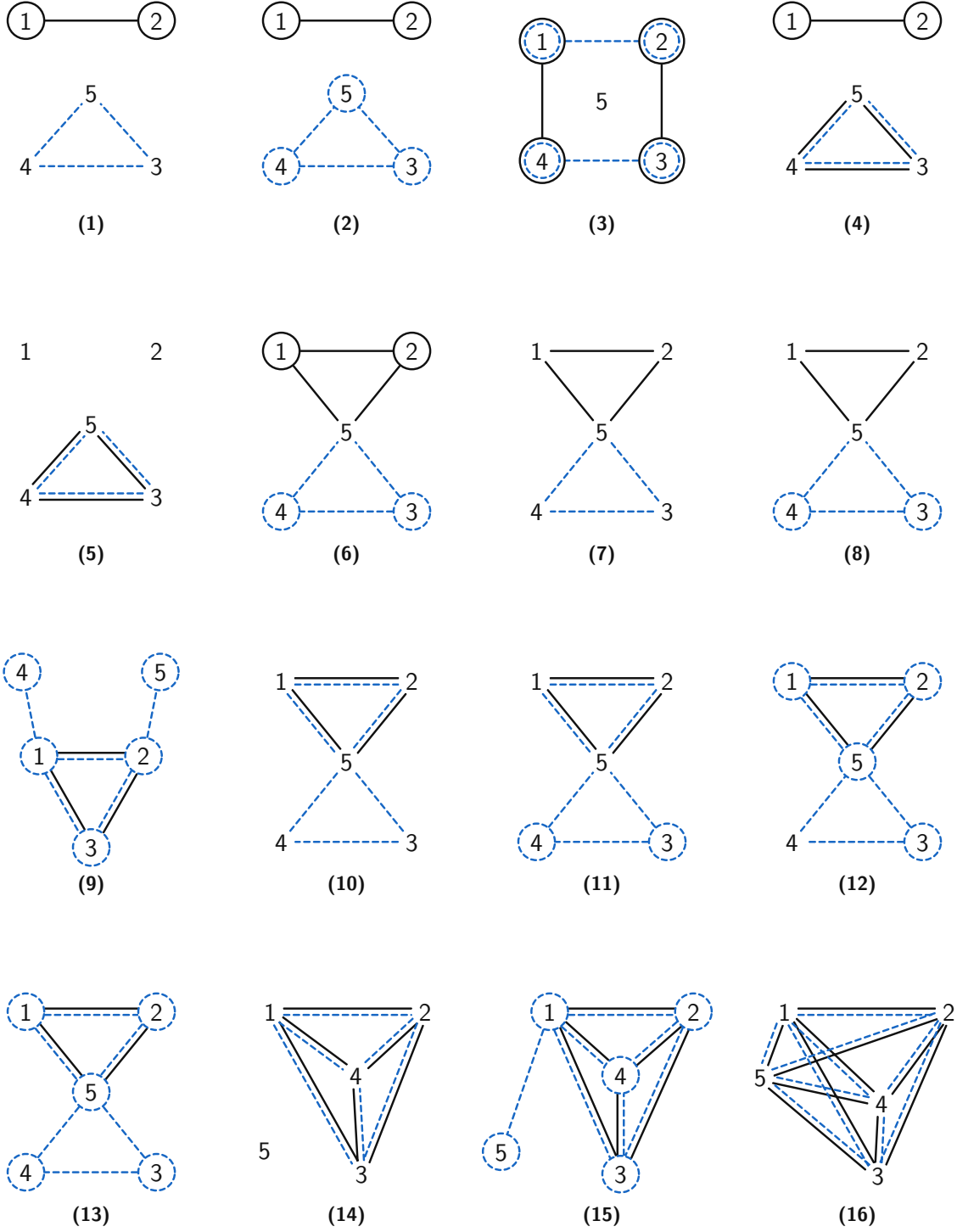
\begin{figure}[p]
\centering
\begingroup
\definecolor{akZ}{HTML}{111111}
\definecolor{akW}{HTML}{1767C4}
\resizebox{.96\linewidth}{!}{%
\begin{tikzpicture}[x=1bp,y=-1bp,line width=1.12bp,line cap=round,
  every node/.style={inner sep=0pt,outer sep=0pt,anchor=base,text=akZ,
    font=\sffamily\fontsize{13.5bp}{13.5bp}\selectfont},
  ak z/.style={draw=akZ},
  ak w/.style={draw=akW,dash pattern=on 3bp off 2.25bp},
  ak wc/.style={draw=akW,dash pattern=on 2.7bp off 2.1bp}]
\path[use as bounding box] (24,94) rectangle (571.276,758);
\draw[ak z] (67.7594,110.0000) -- (117.0594,110.0000);
\draw[ak w] (119.8594,188.0000) -- (64.9594,188.0000) (122.3667,182.3810) -- (97.4522,154.6190) (62.4522,182.3810) -- (87.3667,154.6190);
\path[fill=white] (57.4094,110.0000) circle[radius=9.8500bp] (127.4094,110.0000) circle[radius=9.8500bp] (127.4094,188.0000) circle[radius=7.0500bp] (57.4094,188.0000) circle[radius=7.0500bp] (92.4094,149.0000) circle[radius=7.0500bp];
\draw[ak z] (57.4094,110.0000) circle[radius=9.8000bp] (127.4094,110.0000) circle[radius=9.8000bp];
\node at (57.4094,114.7000) {1}; \node at (127.4094,114.7000) {2}; \node at (127.4094,192.7000) {3}; \node at (57.4094,192.7000) {4}; \node at (92.4094,153.7000) {5};
\node[font=\sffamily\bfseries\fontsize{10.5bp}{10.5bp}\selectfont] at (92.4094,222.0000) {(1)};
\draw[ak z] (204.5783,110.0000) -- (253.8783,110.0000);
\draw[ak w] (253.8783,188.0000) -- (204.5783,188.0000) (257.3155,180.2971) -- (236.1412,156.7029) (201.1412,180.2971) -- (222.3155,156.7029);
\path[fill=white] (194.2283,110.0000) circle[radius=9.8500bp] (264.2283,110.0000) circle[radius=9.8500bp] (264.2283,188.0000) circle[radius=9.8500bp] (194.2283,188.0000) circle[radius=9.8500bp] (229.2283,149.0000) circle[radius=9.8500bp];
\draw[ak z] (194.2283,110.0000) circle[radius=9.8000bp] (264.2283,110.0000) circle[radius=9.8000bp];
\draw[ak wc] (264.2283,188.0000) circle[radius=9.8000bp] (194.2283,188.0000) circle[radius=9.8000bp] (229.2283,149.0000) circle[radius=9.8000bp];
\node at (194.2283,114.7000) {1}; \node at (264.2283,114.7000) {2}; \node at (264.2283,192.7000) {3}; \node at (194.2283,192.7000) {4}; \node at (229.2283,153.7000) {5};
\node[font=\sffamily\bfseries\fontsize{10.5bp}{10.5bp}\selectfont] at (229.2283,222.0000) {(2)};
\draw[ak z] (334.0472,133.5500) -- (334.0472,172.4500) (398.0472,133.5500) -- (398.0472,172.4500);
\draw[ak w] (346.5972,121.0000) -- (385.4972,121.0000) (385.4972,185.0000) -- (346.5972,185.0000);
\path[fill=white] (334.0472,121.0000) circle[radius=12.0500bp] (398.0472,121.0000) circle[radius=12.0500bp] (398.0472,185.0000) circle[radius=12.0500bp] (334.0472,185.0000) circle[radius=12.0500bp] (366.0472,153.0000) circle[radius=7.0500bp];
\draw[ak z] (334.0472,121.0000) circle[radius=12.0000bp] (398.0472,121.0000) circle[radius=12.0000bp] (398.0472,185.0000) circle[radius=12.0000bp] (334.0472,185.0000) circle[radius=12.0000bp];
\draw[ak wc] (334.0472,121.0000) circle[radius=9.3000bp] (398.0472,121.0000) circle[radius=9.3000bp] (398.0472,185.0000) circle[radius=9.3000bp] (334.0472,185.0000) circle[radius=9.3000bp];
\node at (334.0472,125.7000) {1}; \node at (398.0472,125.7000) {2}; \node at (398.0472,189.7000) {3}; \node at (334.0472,189.7000) {4}; \node at (366.0472,157.7000) {5};
\node[font=\sffamily\bfseries\fontsize{10.5bp}{10.5bp}\selectfont] at (366.0472,222.0000) {(3)};
\draw[ak z] (478.2161,110.0000) -- (527.5161,110.0000) (530.5100,189.7000) -- (475.2223,189.7000) (531.6877,183.6607) -- (506.5142,155.6102) (471.5142,181.3898) -- (496.6877,153.3393);
\draw[ak w] (530.5100,186.3000) -- (475.2223,186.3000) (534.2181,181.3898) -- (509.0446,153.3393) (474.0446,183.6607) -- (499.2181,155.6102);
\path[fill=white] (467.8661,110.0000) circle[radius=9.8500bp] (537.8661,110.0000) circle[radius=9.8500bp] (537.8661,188.0000) circle[radius=7.0500bp] (467.8661,188.0000) circle[radius=7.0500bp] (502.8661,149.0000) circle[radius=7.0500bp];
\draw[ak z] (467.8661,110.0000) circle[radius=9.8000bp] (537.8661,110.0000) circle[radius=9.8000bp];
\node at (467.8661,114.7000) {1}; \node at (537.8661,114.7000) {2}; \node at (537.8661,192.7000) {3}; \node at (467.8661,192.7000) {4}; \node at (502.8661,153.7000) {5};
\node[font=\sffamily\bfseries\fontsize{10.5bp}{10.5bp}\selectfont] at (502.8661,222.0000) {(4)};
\draw[ak z] (120.0533,366.7000) -- (64.7656,366.7000) (121.2310,360.6607) -- (96.0575,332.6102) (61.0575,358.3898) -- (86.2310,330.3393);
\draw[ak w] (120.0533,363.3000) -- (64.7656,363.3000) (123.7614,358.3898) -- (98.5879,330.3393) (63.5879,360.6607) -- (88.7614,332.6102);
\path[fill=white] (57.4094,287.0000) circle[radius=7.0500bp] (127.4094,287.0000) circle[radius=7.0500bp] (127.4094,365.0000) circle[radius=7.0500bp] (57.4094,365.0000) circle[radius=7.0500bp] (92.4094,326.0000) circle[radius=7.0500bp];
\node at (57.4094,291.7000) {1}; \node at (127.4094,291.7000) {2}; \node at (127.4094,369.7000) {3}; \node at (57.4094,369.7000) {4}; \node at (92.4094,330.7000) {5};
\node[font=\sffamily\bfseries\fontsize{10.5bp}{10.5bp}\selectfont] at (92.4094,399.0000) {(5)};
\draw[ak z] (204.5783,287.0000) -- (253.8783,287.0000) (200.7620,295.0271) -- (224.4622,324.1445) (257.6947,295.0271) -- (233.9944,324.1445);
\draw[ak w] (253.8783,373.0000) -- (204.5783,373.0000) (257.6947,364.9729) -- (233.9944,335.8555) (200.7620,364.9729) -- (224.4622,335.8555);
\path[fill=white] (194.2283,287.0000) circle[radius=9.8500bp] (264.2283,287.0000) circle[radius=9.8500bp] (264.2283,373.0000) circle[radius=9.8500bp] (194.2283,373.0000) circle[radius=9.8500bp] (229.2283,330.0000) circle[radius=7.0500bp];
\draw[ak z] (194.2283,287.0000) circle[radius=9.8000bp] (264.2283,287.0000) circle[radius=9.8000bp];
\draw[ak wc] (264.2283,373.0000) circle[radius=9.8000bp] (194.2283,373.0000) circle[radius=9.8000bp];
\node at (194.2283,291.7000) {1}; \node at (264.2283,291.7000) {2}; \node at (264.2283,377.7000) {3}; \node at (194.2283,377.7000) {4}; \node at (229.2283,334.7000) {5};
\node[font=\sffamily\bfseries\fontsize{10.5bp}{10.5bp}\selectfont] at (229.2283,399.0000) {(6)};
\draw[ak z] (338.5972,287.0000) -- (393.4972,287.0000) (335.8133,292.8555) -- (361.2811,324.1445) (396.2811,292.8555) -- (370.8133,324.1445);
\draw[ak w] (393.4972,373.0000) -- (338.5972,373.0000) (396.2811,367.1445) -- (370.8133,335.8555) (335.8133,367.1445) -- (361.2811,335.8555);
\path[fill=white] (331.0472,287.0000) circle[radius=7.0500bp] (401.0472,287.0000) circle[radius=7.0500bp] (401.0472,373.0000) circle[radius=7.0500bp] (331.0472,373.0000) circle[radius=7.0500bp] (366.0472,330.0000) circle[radius=7.0500bp];
\node at (331.0472,291.7000) {1}; \node at (401.0472,291.7000) {2}; \node at (401.0472,377.7000) {3}; \node at (331.0472,377.7000) {4}; \node at (366.0472,334.7000) {5};
\node[font=\sffamily\bfseries\fontsize{10.5bp}{10.5bp}\selectfont] at (366.0472,399.0000) {(7)};
\draw[ak z] (475.4161,287.0000) -- (530.3161,287.0000) (472.6322,292.8555) -- (498.1000,324.1445) (533.1000,292.8555) -- (507.6322,324.1445);
\draw[ak w] (527.5161,373.0000) -- (478.2161,373.0000) (531.3325,364.9729) -- (507.6322,335.8555) (474.3998,364.9729) -- (498.1000,335.8555);
\path[fill=white] (467.8661,287.0000) circle[radius=7.0500bp] (537.8661,287.0000) circle[radius=7.0500bp] (537.8661,373.0000) circle[radius=9.8500bp] (467.8661,373.0000) circle[radius=9.8500bp] (502.8661,330.0000) circle[radius=7.0500bp];
\draw[ak wc] (537.8661,373.0000) circle[radius=9.8000bp] (467.8661,373.0000) circle[radius=9.8000bp];
\node at (467.8661,291.7000) {1}; \node at (537.8661,291.7000) {2}; \node at (537.8661,377.7000) {3}; \node at (467.8661,377.7000) {4}; \node at (502.8661,334.7000) {5};
\node[font=\sffamily\bfseries\fontsize{10.5bp}{10.5bp}\selectfont] at (502.8661,399.0000) {(8)};
\draw[ak z] (74.6189,502.3000) -- (110.2000,502.3000) (70.9508,512.0208) -- (88.8202,543.2923) (116.8202,513.7077) -- (98.9508,544.9792);
\draw[ak w] (74.6189,505.7000) -- (110.2000,505.7000) (67.9987,513.7077) -- (85.8681,544.9792) (62.3796,493.8510) -- (57.4393,469.1490) (113.8681,512.0208) -- (95.9987,543.2923) (122.4393,493.8510) -- (127.3796,469.1490);
\path[fill=white] (64.4094,504.0000) circle[radius=9.8500bp] (120.4094,504.0000) circle[radius=9.8500bp] (92.4094,553.0000) circle[radius=9.8500bp] (55.4094,459.0000) circle[radius=9.8500bp] (129.4094,459.0000) circle[radius=9.8500bp];
\draw[ak wc] (64.4094,504.0000) circle[radius=9.8000bp] (120.4094,504.0000) circle[radius=9.8000bp] (92.4094,553.0000) circle[radius=9.8000bp] (55.4094,459.0000) circle[radius=9.8000bp] (129.4094,459.0000) circle[radius=9.8000bp];
\node at (64.4094,508.7000) {1}; \node at (120.4094,508.7000) {2}; \node at (92.4094,557.7000) {3}; \node at (55.4094,463.7000) {4}; \node at (129.4094,463.7000) {5};
\node[font=\sffamily\bfseries\fontsize{10.5bp}{10.5bp}\selectfont] at (92.4094,576.0000) {(9)};
\draw[ak z] (201.5845,462.3000) -- (256.8722,462.3000) (200.1905,468.6320) -- (225.9031,500.2217) (260.9031,470.7783) -- (235.1905,502.3680);
\draw[ak w] (201.5845,465.7000) -- (256.8722,465.7000) (197.5536,470.7783) -- (223.2662,502.3680) (258.2662,468.6320) -- (232.5536,500.2217) (256.6783,550.0000) -- (201.7783,550.0000) (259.4622,544.1445) -- (233.9944,512.8555) (198.9944,544.1445) -- (224.4622,512.8555);
\path[fill=white] (194.2283,464.0000) circle[radius=7.0500bp] (264.2283,464.0000) circle[radius=7.0500bp] (264.2283,550.0000) circle[radius=7.0500bp] (194.2283,550.0000) circle[radius=7.0500bp] (229.2283,507.0000) circle[radius=7.0500bp];
\node at (194.2283,468.7000) {1}; \node at (264.2283,468.7000) {2}; \node at (264.2283,554.7000) {3}; \node at (194.2283,554.7000) {4}; \node at (229.2283,511.7000) {5};
\node[font=\sffamily\bfseries\fontsize{10.5bp}{10.5bp}\selectfont] at (229.2283,576.0000) {(10)};
\draw[ak z] (338.4034,462.3000) -- (393.6911,462.3000) (337.0094,468.6320) -- (362.7220,500.2217) (397.7220,470.7783) -- (372.0094,502.3680);
\draw[ak w] (338.4034,465.7000) -- (393.6911,465.7000) (334.3725,470.7783) -- (360.0851,502.3680) (395.0851,468.6320) -- (369.3725,500.2217) (390.6972,550.0000) -- (341.3972,550.0000) (394.5136,541.9729) -- (370.8133,512.8555) (337.5809,541.9729) -- (361.2811,512.8555);
\path[fill=white] (331.0472,464.0000) circle[radius=7.0500bp] (401.0472,464.0000) circle[radius=7.0500bp] (401.0472,550.0000) circle[radius=9.8500bp] (331.0472,550.0000) circle[radius=9.8500bp] (366.0472,507.0000) circle[radius=7.0500bp];
\draw[ak wc] (401.0472,550.0000) circle[radius=9.8000bp] (331.0472,550.0000) circle[radius=9.8000bp];
\node at (331.0472,468.7000) {1}; \node at (401.0472,468.7000) {2}; \node at (401.0472,554.7000) {3}; \node at (331.0472,554.7000) {4}; \node at (366.0472,511.7000) {5};
\node[font=\sffamily\bfseries\fontsize{10.5bp}{10.5bp}\selectfont] at (366.0472,576.0000) {(11)};
\draw[ak z] (478.0756,462.3000) -- (527.6567,462.3000) (475.6295,470.8449) -- (497.7397,498.0088) (532.7397,472.9912) -- (510.6295,500.1551);
\draw[ak w] (478.0756,465.7000) -- (527.6567,465.7000) (472.9926,472.9912) -- (495.1028,500.1551) (530.1028,470.8449) -- (507.9926,498.0088) (527.5161,550.0000) -- (475.4161,550.0000) (531.3325,541.9729) -- (509.3998,515.0271) (472.6322,544.1445) -- (496.3325,515.0271);
\path[fill=white] (467.8661,464.0000) circle[radius=9.8500bp] (537.8661,464.0000) circle[radius=9.8500bp] (537.8661,550.0000) circle[radius=9.8500bp] (467.8661,550.0000) circle[radius=7.0500bp] (502.8661,507.0000) circle[radius=9.8500bp];
\draw[ak wc] (467.8661,464.0000) circle[radius=9.8000bp] (537.8661,464.0000) circle[radius=9.8000bp] (537.8661,550.0000) circle[radius=9.8000bp] (502.8661,507.0000) circle[radius=9.8000bp];
\node at (467.8661,468.7000) {1}; \node at (537.8661,468.7000) {2}; \node at (537.8661,554.7000) {3}; \node at (467.8661,554.7000) {4}; \node at (502.8661,511.7000) {5};
\node[font=\sffamily\bfseries\fontsize{10.5bp}{10.5bp}\selectfont] at (502.8661,576.0000) {(12)};
\draw[ak z] (67.6189,639.3000) -- (117.2000,639.3000) (65.1728,647.8449) -- (87.2830,675.0088) (122.2830,649.9912) -- (100.1728,677.1551);
\draw[ak w] (67.6189,642.7000) -- (117.2000,642.7000) (62.5359,649.9912) -- (84.6461,677.1551) (119.6461,647.8449) -- (97.5359,675.0088) (117.0594,727.0000) -- (67.7594,727.0000) (120.8758,718.9729) -- (98.9431,692.0271) (63.9431,718.9729) -- (85.8758,692.0271);
\path[fill=white] (57.4094,641.0000) circle[radius=9.8500bp] (127.4094,641.0000) circle[radius=9.8500bp] (127.4094,727.0000) circle[radius=9.8500bp] (57.4094,727.0000) circle[radius=9.8500bp] (92.4094,684.0000) circle[radius=9.8500bp];
\draw[ak wc] (57.4094,641.0000) circle[radius=9.8000bp] (127.4094,641.0000) circle[radius=9.8000bp] (127.4094,727.0000) circle[radius=9.8000bp] (57.4094,727.0000) circle[radius=9.8000bp] (92.4094,684.0000) circle[radius=9.8000bp];
\node at (57.4094,645.7000) {1}; \node at (127.4094,645.7000) {2}; \node at (127.4094,731.7000) {3}; \node at (57.4094,731.7000) {4}; \node at (92.4094,688.7000) {5};
\node[font=\sffamily\bfseries\fontsize{10.5bp}{10.5bp}\selectfont] at (92.4094,753.0000) {(13)};
\draw[ak z] (195.5845,639.3000) -- (263.8722,639.3000) (193.3195,646.5752) -- (235.0997,720.7563) (195.1544,644.0053) -- (229.3232,669.2605) (270.1341,648.4703) -- (241.4893,721.7673) (266.9621,647.2291) -- (240.7951,671.2744) (235.1153,720.7517) -- (232.9468,682.4404);
\draw[ak w] (195.5845,642.7000) -- (263.8722,642.7000) (190.3570,648.2437) -- (232.1372,722.4248) (193.1335,646.7395) -- (227.3023,671.9947) (266.9674,647.2327) -- (238.3226,720.5297) (264.6616,644.7256) -- (238.4946,668.7709) (238.5099,720.5596) -- (236.3413,682.2483);
\path[fill=white] (188.2283,641.0000) circle[radius=7.0500bp] (271.2283,641.0000) circle[radius=7.0500bp] (237.2283,728.0000) circle[radius=7.0500bp] (234.2283,675.0000) circle[radius=7.0500bp] (185.2283,715.0000) circle[radius=7.0500bp];
\node at (188.2283,645.7000) {1}; \node at (271.2283,645.7000) {2}; \node at (237.2283,732.7000) {3}; \node at (234.2283,679.7000) {4}; \node at (185.2283,719.7000) {5};
\node[font=\sffamily\bfseries\fontsize{10.5bp}{10.5bp}\selectfont] at (229.2283,753.0000) {(14)};
\draw[ak z] (348.2567,639.3000) -- (404.8378,639.3000) (343.6916,649.6755) -- (373.5186,717.9636) (346.7893,646.5407) -- (369.5723,666.9255) (412.4223,651.0116) -- (381.7747,719.3792) (408.4688,648.9904) -- (384.8600,669.5724) (374.3472,717.7906) -- (374.3472,685.2094);
\draw[ak w] (348.2567,642.7000) -- (404.8378,642.7000) (340.5759,651.0364) -- (370.4029,719.3245) (344.5222,649.0745) -- (367.3052,669.4593) (334.6572,650.7791) -- (315.4373,706.2209) (409.3197,649.6208) -- (378.6722,717.9884) (406.2345,646.4276) -- (382.6257,667.0096) (377.7472,717.7906) -- (377.7472,685.2094);
\path[fill=white] (338.0472,641.0000) circle[radius=9.8500bp] (415.0472,641.0000) circle[radius=9.8500bp] (376.0472,728.0000) circle[radius=9.8500bp] (376.0472,675.0000) circle[radius=9.8500bp] (312.0472,716.0000) circle[radius=9.8500bp];
\draw[ak wc] (338.0472,641.0000) circle[radius=9.8000bp] (415.0472,641.0000) circle[radius=9.8000bp] (376.0472,728.0000) circle[radius=9.8000bp] (376.0472,675.0000) circle[radius=9.8000bp] (312.0472,716.0000) circle[radius=9.8000bp];
\node at (338.0472,645.7000) {1}; \node at (415.0472,645.7000) {2}; \node at (376.0472,732.7000) {3}; \node at (376.0472,679.7000) {4}; \node at (312.0472,720.7000) {5};
\node[font=\sffamily\bfseries\fontsize{10.5bp}{10.5bp}\selectfont] at (366.0472,753.0000) {(15)};
\draw[ak z] (471.2223,639.3000) -- (544.5100,639.3000) (468.9782,646.5560) -- (511.7102,721.7644) (470.3857,644.8076) -- (511.6783,683.7179) (462.8985,648.4877) -- (455.0211,669.6961) (550.5106,648.4273) -- (518.3431,722.9206) (548.8808,647.9347) -- (521.5915,685.0784) (545.4722,645.0149) -- (458.3734,675.1977) (512.5451,721.5665) -- (513.7916,697.2594) (507.1427,725.5653) -- (455.4008,682.0365) (508.3170,690.1130) -- (457.6994,679.2108);
\draw[ak w] (471.2223,642.7000) -- (544.5100,642.7000) (466.0221,648.2356) -- (508.7541,723.4440) (468.0540,647.2821) -- (509.3466,686.1924) (459.7112,647.3039) -- (451.8338,668.5123) (547.3892,647.0794) -- (515.2217,721.5727) (546.1408,645.9216) -- (518.8515,683.0653) (544.3589,641.8023) -- (457.2601,671.9851) (515.9407,721.7406) -- (517.1872,697.4335) (509.3314,722.9635) -- (457.5896,679.4347) (509.0329,686.7892) -- (458.4153,675.8870);
\path[fill=white] (463.8661,641.0000) circle[radius=7.0500bp] (551.8661,641.0000) circle[radius=7.0500bp] (513.8661,729.0000) circle[radius=7.0500bp] (515.8661,690.0000) circle[radius=7.0500bp] (450.8661,676.0000) circle[radius=7.0500bp];
\node at (463.8661,645.7000) {1}; \node at (551.8661,645.7000) {2}; \node at (513.8661,733.7000) {3}; \node at (515.8661,694.7000) {4}; \node at (450.8661,680.7000) {5};
\node[font=\sffamily\bfseries\fontsize{10.5bp}{10.5bp}\selectfont] at (502.8661,753.0000) {(16)};
\end{tikzpicture}%
}
\endgroup
\caption{The sixteen five-body diagrams, redrawn from \cite[Figure~11, p.~568]{AK}, with the original numbering and vertex labels. Black solid strokes and circles represent $z$-strokes and $z$-circles; blue dashed strokes and circles represent $w$-strokes and $w$-circles. A $zw$-edge is drawn with two parallel strokes.}
\label{fig:ak-five-body}
\end{figure}

We use the diagram numbering and vertex labels of \cite[Figure~11]{AK}. The sixteen five-body diagrams are reproduced in Figure~\ref{fig:ak-five-body}, using black solid strokes and circles for $z$ and blue dashed strokes and circles for $w$.

Lemma~\ref{lem:2.2} reduces the subsequent analysis to Diagrams 1--8. We first describe their structure. Following \cite[Definition~5 and Section~3]{AK}, an \emph{isolated component} has no $z$- or $w$-stroke joining it to its complement. We call a triangle \emph{collapsing} if its three mutual signed distances tend to zero as $\epsilon\to0$. By the AK distance estimates recalled below, a triangle of nonmaximal edges is collapsing, whereas a triangle whose three edges are maximal has its three signed distances bounded and bounded away from zero \cite[Section~9, p.~580]{AK}.

For the purposes of this paper, we group the first eight diagrams into four broad structural classes:
\begin{enumerate}[label=(\roman*)]
\item \emph{Quadrilateral with an isolated vertex (Diagram 3).}
Vertices $1,2,3,4$ form a quadrilateral with alternating $z$- and $w$-edges, and vertex $5$ is isolated.
\item \emph{Isolated collapsing triangle (Diagrams 1, 4, and 5).}
The triangle on $S_3$ is isolated and collapsing. Its edges are $w$-edges in Diagram 1 and $zw$-edges in Diagrams 4 and 5. In Diagrams 1 and 4, the remaining vertices $1,2$ are joined by a maximal $z$-edge; in Diagram 5, they are two isolated vertices.
\item \emph{Isolated triangle with maximal edges (Diagram 2).}
The triangle on $S_3$ consists of three maximal $w$-edges and forms an isolated component, accompanied by the isolated maximal $z$-edge $12$.
\item \emph{Butterfly (Diagrams 6, 7, and 8).}
The $z$-triangle on $T_3$ and the $w$-triangle on $S_3$ share vertex $5$. In Diagram 6, all three edges of each wing are maximal. In Diagram 7, both wings are collapsing. In Diagram 8, the $T_3$ wing is collapsing and all three edges of the $S_3$ wing are maximal.
\end{enumerate}

For the pair components defined in Section~\ref{sec:2.1}, the distance identity \eqref{eq:2.1} gives
\[
Z_{ij}=\frac{1}{r_{ij}w_{ij}},\qquad
W_{ij}=\frac{1}{r_{ij}z_{ij}},\qquad
w_{ij}Z_{ij}=z_{ij}W_{ij}=\frac{1}{r_{ij}}.
\]
Consequently, for $i,j,k,l\in S_5$ with $i\ne j$ and $k\ne l$,
\[
Z_{ij}^2z_{ij}w_{ij}^3=1,\qquad
W_{ij}^2z_{ij}^3w_{ij}=1,
\]
and
\[
\left(\frac{Z_{ij}}{Z_{kl}}\right)^2
=\frac{z_{kl}}{z_{ij}}\left(\frac{w_{kl}}{w_{ij}}\right)^3,\qquad
\left(\frac{W_{ij}}{W_{kl}}\right)^2
=\left(\frac{z_{kl}}{z_{ij}}\right)^3\frac{w_{kl}}{w_{ij}}.
\]
These are complex identities on the chosen signed-distance branches; taking moduli gives the corresponding scale comparisons without choosing any additional radical branches.

For a singular sequence, AK use the \emph{balanced normalization} to make the maxima of the two augmented coordinate lists equal to $\epsilon^{-2}$, with $\epsilon\to0$ \cite[Section~3, pp.~547--548]{AK}. This uses reciprocal complex rotations $(z,w)\mapsto(\mathfrak{a}z,\mathfrak{a}^{-1}w)$, which preserve every signed distance and $U,I$. Such a rotation may leave the original gauge component; no subsequent argument requires it to stay in that gauge.

We write $f\approx g$ for two-sided comparison of moduli, $f=o(g)$ for $f/g\to0$, and $f\sim g$ for $f/g\to1$. In the notation of \cite{AK} these three relations are $\approx$, $\prec$ and $\sim$. For the same relations we also use $f\preceq g$ for $f=O(g)$ (that is, $f/g$ is bounded), $f\succeq g$ for $g=O(f)$, and $f\succ g$ for $g=o(f)$. Thus
\[
g=O(f)\iff f\succeq g,\qquad g=o(f)\iff g\prec f\iff f\succ g .
\]
Finite sums preserve $O$ and $o$ upper bounds. Two-sided comparisons, however, do not in general pass from the summands to their sum, since complex terms may cancel. The constants implicit in $O(\cdot)$ and $\preceq$ depend only on the fixed masses and diagram data. All statements are made after the finite subsequence selections needed to fix diagram, signed branches, and leading ratios. Unless a different limiting variable is explicitly specified, all limits and asymptotic estimates in the diagram-degeneration arguments are taken along the selected subsequence as $\epsilon\to0$. For any complex-valued quantity considered below, convergence to infinity is understood in the Riemann sphere, equivalently as divergence of its modulus. In particular, after subsequence selection, each nonzero signed distance either tends to zero, has a finite nonzero limit, or tends to infinity. These conventions also follow by choosing algebraic curve branches and Puiseux parameters.

The following are the specific published scale facts used below. AK Estimates 1--2 \cite[Section 3, p.~548]{AK} give
\begin{equation}\label{eq:2.4}
r_{ij}\succeq\epsilon,\qquad
r_{ij}\approx\epsilon\iff ij\text{ is a }zw\text{ edge},
\end{equation}
and on a $zw$ edge both $z_{ij},w_{ij}\approx\epsilon$. In the original balanced AK coordinates, a $z$-stroke on $ij$ is equivalent to $Z_{ij}\approx\epsilon^{-2}$; absence of a $z$-stroke gives $Z_{ij}=o(\epsilon^{-2})$. The analogous statements hold for $w$-strokes and $W_{ij}$. Assuming that $ij$ carries a $z$-stroke,
\[
ij\text{ is a maximal }z\text{-edge}
\iff r_{ij}\approx1
\iff z_{ij}\approx\epsilon^{-2}
\iff w_{ij}\approx\epsilon^2.
\]
Exchanging $z$ and $w$ gives the corresponding maximal $w$-edge criterion. These statements use the original balanced normalization; a further reciprocal coordinate scaling need not preserve its two maxima. A vertex not circled in a coordinate has that coordinate $o(\epsilon^{-2})$, by the definition of the diagram.

In diagram 5 the only edges form the $zw$ triangle on $3,4,5$; all five vertices are uncircled and $1,2$ are isolated. Consequently all absolute coordinates are $o(\epsilon^{-2})$, the internal triangle differences are comparable to $\epsilon$ in both coordinates, and every cross distance from $S_2$ to $S_3$ is much larger than $\epsilon$. The last assertion follows from both directions of \eqref{eq:2.4}, not merely from its lower bound.

\begin{lemma}[A singular sequence on a component containing a physical point]\label{lem:2.3}
Let the five masses be positive. If a positive-dimensional irreducible component of the collision-free, normalized, rotation-gauged complex solution variety contains a physical point, it contains a singular sequence along which some 4-product tends to infinity. After balanced normalization and subsequence selection, this sequence approaches one of AK diagrams 1--8. Along this sequence, $U=I$ is a positive real constant satisfying the physical lower bound \eqref{eq:3.1}.
\end{lemma}
\begin{proof}
A fixed squared-distance tuple has only finitely many points in the gauge, as shown in Section 2. Thus some squared distance is nonconstant on the component. The proof of Lemma~\ref{lem:2.2} now applies: push this function to zero, use AK Proposition 9 to obtain a nonconstant 4-product, and push the 4-product to infinity. AK Proposition 8 excludes diagrams 9--16. Lemma~\ref{lem:2.1} and the single physical point give the asserted positive value of $U=I$ and the physical lower bound. No density of physical points in the component is required.
\end{proof}

\section{The exclusion of Diagram 3}\label{sec:3}

Throughout this section the masses are positive and the central-configuration multiplier is normalized to one. In the complex equations the variables $r_{ij}$ are signed square roots of $(z_i-z_j)(w_i-w_j)$; they need not be real or positive. We use
\[
\begin{gathered}
U=\sum\limits_{i<j}\frac{m_im_j}{r_{ij}},\qquad
I=\sum\limits_i m_iz_iw_i,\\
M=\sum\limits_i m_i,\qquad
\sum\limits_i m_iz_i=\sum\limits_i m_iw_i=0.
\end{gathered}
\]
The force convention is $q_i=\sum\limits_{j\ne i}m_j(q_i-q_j)r_{ij}^{-3}$. We use the asymptotic conventions of Section~\ref{sec:2.2}.

\subsection{The physical lower bound}\label{sec:3.1}

\begin{lemma}\label{lem:3.1}
At a physical normalized central configuration,
\begin{equation}\label{eq:3.1}
U=I\ge\frac{e_2(S_5)}{M^{1/3}}.
\end{equation}
On every irreducible component of the complex normalized solution variety containing a physical point, $U$ and $I$ are the same positive real constant, and this constant satisfies \eqref{eq:3.1}.
\end{lemma}
\begin{proof}
Pairing the normalized force equations with the mass-weighted positions gives $U=I$. The center condition also gives $\sum\limits_{i<j}m_im_jr_{ij}^2=MI$. For physical distances apply Jensen's inequality to the strictly convex function $s\mapsto s^{-1/2}$, with weights $m_im_j/e_2(S_5)$:
\[
\frac U{e_2(S_5)}\ge\left(\frac1{e_2(S_5)}\sum\limits_{i<j}m_im_jr_{ij}^2\right)^{-1/2}
=\left(\frac{MU}{e_2(S_5)}\right)^{-1/2}.
\]
Since $U>0$, this is equivalent to \eqref{eq:3.1}. Lemma~\ref{lem:2.1} makes $U=I$ constant on each irreducible normalized component. A physical point fixes its value and imposes the physical lower bound \eqref{eq:3.1}. Complex orthogonal coordinate changes preserve $U$ and $I$, even when they leave a particular rotation gauge; consequently the same constants persist under AK's balancing normalization and along its singular sequences.
\end{proof}

We refer to \eqref{eq:3.1} as the \emph{physical lower bound}. Since the multiplier is normalized to one, $U=I$, so this bound applies equally to the potential and the moment of inertia on every irreducible component containing a physical point.

\begin{lemma}[Potential estimate]\label{lem:3.2}
Partition $S_5$ into at least two nonempty sets $\mathcal I_1,\dots,\mathcal I_r$, with $r\geq2$, total masses $e_1(\mathcal I_j)$ and internal pair sums $e_2(\mathcal I_j)$. Suppose along a complex degeneration that, as $\epsilon\to0$, all distances between distinct sets diverge and every within-set distance satisfies $r_{kl}^3\to e_1(\mathcal I_j)$ for $k,l\in \mathcal I_j$, $k\ne l$, and $j=1,\dots,r$. Such a degeneration cannot have the constant potential of an irreducible component containing a physical point.
\end{lemma}
\begin{proof}
All choices of cube-root phases satisfy $|U|\le\sum\limits_j e_2(\mathcal I_j)/e_1(\mathcal I_j)^{1/3}$ in the limit. Denote the pair sum between distinct groups by 
$$e_2^{\mathrm{cross}}=\sum\limits_{1\leq p<q\leq r}e_1(\mathcal I_p)e_1(\mathcal I_q)>0,$$ 
so $e_2(S_5)=\sum\limits_{j=1}^{r} e_2(\mathcal I_j)+e_2^{\mathrm{cross}}$ and $M=\sum\limits_{j=1}^{r} e_1(\mathcal I_j)$. The tangent inequality $s^{1/3}\le1+(s-1)/3$, applied at $s=M/e_1(\mathcal I_j)$, gives
\[
\begin{aligned}
M^{1/3}\sum\limits_{j=1}^{r}\frac{e_2(\mathcal I_j)}{e_1(\mathcal I_j)^{1/3}}-\sum\limits_{j=1}^{r} e_2(\mathcal I_j)
&\le\sum\limits_{j=1}^{r}\frac{e_2(\mathcal I_j)(M-e_1(\mathcal I_j))}{3e_1(\mathcal I_j)}\\
&\le\sum\limits_{j=1}^{r}\frac{e_1(\mathcal I_j)(M-e_1(\mathcal I_j))}6=\frac {e_2^{\mathrm{cross}}}{3}<e_2^{\mathrm{cross}},
\end{aligned}
\]
where $e_2(\mathcal I_j)\le e_1(\mathcal I_j)^2/2$ was used. Hence $\sum\limits_{j=1}^{r}e_2(\mathcal I_j)/e_1(\mathcal I_j)^{1/3}<e_2(S_5)/M^{1/3}$, contrary to Lemma~\ref{lem:3.1}.
\end{proof}

\subsection{The potential estimate}\label{sec:3.2}

\begin{proposition}\label{prop:3.3}
AK diagram 3 cannot occur on a complex normalized component containing a physical central configuration, for any five positive masses.
\end{proposition}
\begin{proof}
Write the four quadrilateral masses in cyclic order as $m_1,m_2,m_3,m_4$ and the isolated mass as $m_5$. AK Section 8.3 imposes $m_1m_3=m_2m_4$ and refers to the clustering in Section 5.2: $z:12\ldots34$ and $w:14\ldots23$ \cite[pp.~558,575]{AK}. The four side distances have finite nonzero limits; the two diagonals and all four distances to the isolated vertex diverge. For example, the leading $z$ equations for the pair $14$ are
\[
z_1=m_4Z_{14}+o(z_{14}),\qquad
z_4=-m_1Z_{14}+o(z_{14}).
\]
Subtracting gives $z_{14}=(m_1+m_4)Z_{14}+o(z_{14})$, hence $r_{14}^3=z_{14}/Z_{14}\to m_1+m_4$. The same argument uses $Z_{23}$ for the other leading $z$ pair, and $W_{12},W_{34}$ for the two leading $w$ pairs. Thus, as $\epsilon\to0$,
\[
\begin{aligned}
r_{12}^3=\frac{w_{12}}{W_{12}}&\longrightarrow m_1+m_2,&
r_{23}^3=\frac{z_{23}}{Z_{23}}&\longrightarrow m_2+m_3,\\
r_{34}^3=\frac{w_{34}}{W_{34}}&\longrightarrow m_3+m_4,&
r_{14}^3=\frac{z_{14}}{Z_{14}}&\longrightarrow m_1+m_4.
\end{aligned}
\]
At this limit the four sides dominate the potential, and
\[
\begin{aligned}
U&\le\sum\limits_{i<j}\left|\frac{m_im_j}{r_{ij}}\right|
=\sum\limits_{i<j}\frac{m_im_j}{|r_{ij}|}\ \longrightarrow\ U_4,\\
U_4&:=\frac{m_1m_2}{(m_1+m_2)^{1/3}}+\frac{m_2m_3}{(m_2+m_3)^{1/3}}
+\frac{m_3m_4}{(m_3+m_4)^{1/3}}+\frac{m_4m_1}{(m_4+m_1)^{1/3}},
\end{aligned}
\]
because $r_{ij}^3\to m_i+m_j$ on the four sides, while the other six distances diverge.
Set
\[
\theta_i=\frac{m_i}{e_1(S_4)},\qquad i\in S_4.
\]
Since $e_1(S_4)=M-m_5$, positivity and the mass relation give
\[
\theta_i>0,\qquad \sum\limits_{i\in S_4}\theta_i=1,\qquad
\theta_1\theta_3=\theta_2\theta_4.
\]
Using these identities, we obtain
\[
\frac{e_2(S_4)}{e_1(S_4)^2}
=\sum\limits_{\substack{i<j\\i,j\in S_4}}\theta_i\theta_j
\]
and
\[
\begin{aligned}
\frac{U_4}{e_1(S_4)^{5/3}}
={}&(\theta_1+\theta_4)(\theta_2+\theta_3)
\bigl[(\theta_1+\theta_2)^{5/3}+(\theta_3+\theta_4)^{5/3}\bigr]\\
&+(\theta_1+\theta_2)(\theta_3+\theta_4)
\bigl[(\theta_1+\theta_4)^{5/3}+(\theta_2+\theta_3)^{5/3}\bigr].
\end{aligned}
\]
By concavity of $x^{2/3}$ for $x>0$ at $x=1$, we have $x^{2/3}\le(1+2x)/3$, hence $x^{5/3}+(1-x)^{5/3}\le(4x^2-4x+3)/3$ for $0<x<1$. Applying this to $x=\theta_1+\theta_2$ and $x=\theta_1+\theta_4$, and using $\theta_1\theta_3=\theta_2\theta_4$, gives
\[
\frac{U_4}{e_1(S_4)^{5/3}}
\le\frac{e_2(S_4)}{e_1(S_4)^2}-\frac23\theta_1\theta_3
<\frac{e_2(S_4)}{e_1(S_4)^2},
\]
where the strict inequality follows from $\theta_1\theta_3>0$.
Finally $(e_2(S_4)+e_1(S_4)m_5)/M^{1/3}$ is increasing for $m_5\ge0$ with $m_1,\dots,m_4$ fixed: its derivative has the positive numerator $e_1(S_4)^2+\tfrac23 e_1(S_4)m_5-e_2(S_4)/3$. It follows that
\[
U\le U_4<\frac{e_2(S_4)}{e_1(S_4)^{1/3}}<\frac{e_2(S_4)+e_1(S_4)m_5}{M^{1/3}}=\frac{e_2(S_5)}{M^{1/3}},
\]
contrary to \eqref{eq:3.1}.
\end{proof}

\section{The exclusion of Diagrams 1 and 4}\label{sec:4}

In both diagrams, $S_3$ forms an isolated collapsing triangle in the sense of Section~\ref{sec:2.2}; its three internal signed distances tend to zero along the singular sequence. The following comparison for the moment of inertia applies throughout the positive five-body mass cone.

\begin{proposition}\label{prop:4.1}
For any fixed positive masses $m_1,\ldots,m_5$, no irreducible component of the normalized complex central-configuration variety containing a physical point admits a singular sequence approaching AK diagram 1 or diagram 4 after balanced normalization.
\end{proposition}
\begin{proof}
Suppose such a sequence exists, and use the balanced normalization of AK. All limits below are taken as $\epsilon\to0$, and order comparisons refer to moduli. In both diagrams, $12$ is a $z$ edge, there are no $z$ strokes between $S_2$ and $S_3$, and the vertices in $S_3$ have no $z$ circles. Thus
\[
Z_{12}\approx\epsilon^{-2},\qquad
Z_{ij}=o(\epsilon^{-2})\quad(i\in S_2,\ j\in S_3),\qquad
z_j=o(\epsilon^{-2})\quad(j\in S_3).
\]
The original $z$ equations \eqref{eq:2.2} at vertices 1 and 2 give
\[
z_1=m_2Z_{12}+o(\epsilon^{-2}),\qquad
z_2=-m_1Z_{12}+o(\epsilon^{-2}).
\]
Since $m_1,m_2>0$, it follows that
\[
z_{ij}\approx\epsilon^{-2}\quad(i\in S_2,\ j\in S_3),\qquad
z_{12}=(m_1+m_2)Z_{12}+o(\epsilon^{-2}).
\]
For each of the six cross pairs, therefore,
\[
\delta_{ij}=\frac{Z_{ij}}{z_{ij}}=o(1),\qquad
r_{ij}\longrightarrow\infty.
\]
Using the exact identity $r_{ij}^2=z_{ij}w_{ij}$, we obtain
\[
W_{ij}=\frac{w_{ij}}{r_{ij}^3}
=\frac{1}{r_{ij}z_{ij}}=o(\epsilon^2)
\qquad(i\in S_2,\ j\in S_3).
\]

Multiply the $w$ equations for vertices 1 and 2 by their respective masses and add them. The internal terms cancel, yielding
\[
(m_1+m_2)C_{S_2,w}
=\sum\limits_{i\in S_2}\sum\limits_{j\in S_3}m_im_jW_{ij}
=o(\epsilon^2).
\]
Hence $C_{S_2,w}=o(\epsilon^2)$, and total centering gives $C_{S_3,w}=o(\epsilon^2)$ as well. Balanced normalization bounds every $z_i$, $i\in S_5$, by $O(\epsilon^{-2})$. Consequently,
\[
(C_{S_2,z}-C_{S_3,z})(C_{S_2,w}-C_{S_3,w})=o(1).
\]

In both diagrams, the triangle on $S_3$ has all three $w$ edges and no $w$ circles. For distinct $j,k\in S_3$, this gives
\[
W_{jk}\approx\epsilon^{-2},\qquad
w_{jk}=o(\epsilon^{-2}),\qquad
r_{jk}^3=\frac{w_{jk}}{W_{jk}}\longrightarrow0.
\]
Its internal moment of inertia therefore satisfies
\[
I_{S_3}=\frac{1}{e_1(S_3)}
\sum\limits_{\substack{j<k\\j,k\in S_3}}m_jm_kr_{jk}^2
\longrightarrow0.
\]
The binary $z$ equations above also give
\[
r_{12}^3=\frac{z_{12}}{Z_{12}}\longrightarrow m_1+m_2.
\]
Using the relative coordinates defined in \eqref{eq:2.relative_coordinates}, write
$z_i=C_{S,z}+\zeta_{S,i}$ and $w_i=C_{S,w}+\omega_{S,i}$ for $S=S_2,S_3$.
The center-of-mass identities
\[
\sum\limits_{i\in S}m_i\zeta_{S,i}
=\sum\limits_{i\in S}m_i\omega_{S,i}=0
\]
cancel the mixed terms in each cluster. Hence
\begin{equation}\label{eq:4.1}
\begin{aligned}
I&=\sum\limits_{i\in S_2}m_i z_iw_i+\sum\limits_{i\in S_3}m_i z_iw_i\\
&=\sum\limits_{i\in S_2}m_i(C_{S_2,z}+\zeta_{S_2,i})(C_{S_2,w}+\omega_{S_2,i})+\sum\limits_{i\in S_3}m_i(C_{S_3,z}+\zeta_{S_3,i})(C_{S_3,w}+\omega_{S_3,i})\\
&=I_{S_2}+I_{S_3}+(m_1+m_2)C_{S_2,z}C_{S_2,w}+e_1(S_3)C_{S_3,z}C_{S_3,w}\\
&=\frac{m_1m_2}{m_1+m_2}r_{12}^2+I_{S_3}+\frac{(m_1+m_2)e_1(S_3)}{M}(C_{S_2,z}-C_{S_3,z})(C_{S_2,w}-C_{S_3,w}).
\end{aligned}
\end{equation}
The last equality uses $I_{S_2}=m_1m_2r_{12}^2/(m_1+m_2)$ and the total center-of-mass relation
\[
(m_1+m_2)C_{S_2}+e_1(S_3)C_{S_3}=0,
\]
which gives
\[
\begin{aligned}
C_{S_2}&=\frac{e_1(S_3)}{M}(C_{S_2}-C_{S_3}),\\
C_{S_3}&=-\frac{m_1+m_2}{M}(C_{S_2}-C_{S_3}).
\end{aligned}
\]
The preceding estimates and \eqref{eq:4.1} give
\[
I=\frac{m_1m_2}{m_1+m_2}r_{12}^2+o(1),\qquad
|r_{12}|\longrightarrow(m_1+m_2)^{1/3}.
\]
Since $I$ is a positive real constant by Lemma~\ref{lem:2.1}, taking moduli and passing to the limit gives
\[
I=|I|=\frac{m_1m_2}{(m_1+m_2)^{1/3}}
<\frac{e_2(S_5)}{M^{1/3}}.
\]
The strict inequality is the mass inequality in the proof of Lemma~\ref{lem:3.2}, applied to the blocks $S_2,\{3\},\{4\},\{5\}$; its use here does not require the separated Lagrange hypotheses of that lemma. On the other hand, Lemmas~\ref{lem:2.1} and~\ref{lem:3.1} imply that $I$ is a positive real constant on the component and satisfies $I\ge e_2(S_5)/M^{1/3}$, a contradiction.

The reciprocal coordinate changes used for balancing preserve every signed distance and $I$, so applying them separately to the terms of the sequence leaves this contradiction unchanged. The modulus estimate covers every signed-distance branch.
\end{proof}

\section{The exclusion of Diagram 7}\label{sec:8}

The condition that some 4-product tends to infinity simplifies the analysis of diagram 7. All complex distance phases are retained.

The global sets $S_2,T_2$ in \eqref{eq:2.index_sets} are the two outer pairs; $T_3,S_3$ are the vertex sets of the two wings, which share vertex 5. Write $q_{ij}=q_i-q_j$. In the notation of Section~\ref{sec:2}, the contribution of vertices 1 and 2 to the right-hand side of the equation for vertex 5 is
\[
f_5^{T_3}=\sum\limits_{i\in S_2}m_iq_{5i}\delta_{i5}.
\]

\begin{proposition}\label{prop:8.1}
For any five positive masses, a complex normalized component containing a physical point cannot admit a singular sequence approaching AK diagram 7 along which some 4-product tends to infinity.
\end{proposition}

\begin{proof}
Define the coordinate separation scales of the two wings by
\[
\begin{aligned}
\Theta_{T_3,z}&=\max\limits_{\substack{i<j\\i,j\in T_3}}|z_{ij}|,\qquad&
\Theta_{T_3,w}&=\max\limits_{\substack{i<j\\i,j\in T_3}}|w_{ij}|,\\
\Theta_{S_3,z}&=\max\limits_{\substack{i<j\\i,j\in S_3}}|z_{ij}|,&
\Theta_{S_3,w}&=\max\limits_{\substack{i<j\\i,j\in S_3}}|w_{ij}|.
\end{aligned}
\]
The diagram gives
\[
\Theta_{S_3,z}=o(\Theta_{T_3,z}),\qquad
\Theta_{T_3,w}=o(\Theta_{S_3,w}).
\]
Every distance in either wing collapses, so
\[
\Theta_{T_3,z}\Theta_{T_3,w}\to0,\qquad
\Theta_{S_3,z}\Theta_{S_3,w}\to0.
\]
Each of the other four squared distances has modulus comparable to $\Theta_{T_3,z}\Theta_{S_3,w}$.

If $\Theta_{T_3,z}\Theta_{S_3,w}$ stayed bounded, all ten distances would be bounded. The four crossed edges join vertices in $S_2$ to vertices in $T_2$ and do not involve vertex 5, so they cannot form a three-cycle. Hence every three-cycle contains at least one edge within a wing, whose distance collapses. Each 4-product \eqref{eq:2.3} would tend to zero, contradicting the assumption that some 4-product tends to infinity. Consequently $\Theta_{T_3,z}\Theta_{S_3,w}\to\infty$.

With $C_{S_2}$ and $C_{T_2}$ denoting the centers of mass of the two outer pairs, the exact decomposition of the moment of inertia gives
\[
\begin{aligned}
I={}&\frac{m_1m_2}{m_1+m_2}r_{12}^2+\frac{m_3m_4}{m_3+m_4}r_{34}^2\\
&+\frac{(m_1+m_2)e_1(S_3)}{M}q_{S_2,5}^2\\
&+\frac{(m_3+m_4)e_1(T_3)}{M}q_{T_2,5}^2\\
&-\frac{2(m_1+m_2)(m_3+m_4)}{M}q_{S_2,5}\cdot q_{T_2,5}.
\end{aligned}
\]
The first four terms vanish: each is a fixed linear combination of squared distances within a collapsed wing. Also
\[
\begin{aligned}
\omega_{S_2,5}\zeta_{T_2,5}&=O(\Theta_{T_3,w}\Theta_{S_3,z})\\
&=O\left(\frac{(\Theta_{T_3,z}\Theta_{T_3,w})(\Theta_{S_3,z}\Theta_{S_3,w})}
{\Theta_{T_3,z}\Theta_{S_3,w}}\right)\longrightarrow0.
\end{aligned}
\]
It remains to control $\zeta_{S_2,5}\omega_{T_2,5}$.

We use the skew product defined in Section~\ref{sec:2.1}. Multiplying the normalized force equations for $i\in S_2$ by $m_i$ and summing cancels the interaction between vertices 1 and 2. Each cross interaction $q_{ij}\delta_{ij}$ has components $(Z_{ij},W_{ij})$ in the $(z,w)$ coordinates. The resulting vector force equation and its torque identity give
\[
\begin{aligned}
(m_1+m_2)C_{S_2}
&=-m_5f_5^{T_3}+\sum\limits_{\substack{i\in S_2\\j\in T_2}}m_im_jq_{ij}\delta_{ij},\\
0&=m_5[q_5,f_5^{T_3}]+\sum\limits_{\substack{i\in S_2\\j\in T_2}}m_im_j\delta_{ij}[q_i,q_j].
\end{aligned}
\]
Taking the bracket with $q_5$ in the first equation and using the second yields
\[
(m_1+m_2)[q_5,C_{S_2}]
=\sum\limits_{\substack{i\in S_2\\j\in T_2}}m_im_j\delta_{ij}\bigl([q_5,q_{ij}]+[q_i,q_j]\bigr).
\]
For each $i\in S_2$, $j\in T_2$, the definition of the skew product gives the exact identity for the complete summand:
\[
\delta_{ij}\bigl([q_5,q_{ij}]+[q_i,q_j]\bigr)
=w_{i5}Z_{ij}-z_{i5}W_{ij}.
\]
The separation scales of this diagram give
\[
\begin{aligned}
|z_{i5}|&\le\Theta_{T_3,z}\approx z_{ij},\\
|w_{i5}|&\le\Theta_{T_3,w}=o(\Theta_{S_3,w}),\qquad
w_{ij}\approx\Theta_{S_3,w}.
\end{aligned}
\]
Since $w_{ij}Z_{ij}=z_{ij}W_{ij}=1/r_{ij}$, we obtain
\[
\begin{aligned}
|w_{i5}Z_{ij}-z_{i5}W_{ij}|
&\le |w_{i5}|\,|Z_{ij}|+|z_{i5}|\,|W_{ij}|\\
&\preceq |w_{ij}Z_{ij}|+|z_{ij}W_{ij}|
=\frac{2}{|r_{ij}|}\longrightarrow0.
\end{aligned}
\]
All four crossed distances diverge, so substitution into the finite mass-weighted sum yields $[q_5,C_{S_2}]\to0$. The center identity
\[
(m_1+m_2)C_{S_2}+(m_3+m_4)C_{T_2}+m_5q_5=0
\]
gives
\[
[q_{S_2,5},q_{T_2,5}]=\frac{M}{m_3+m_4}[q_5,C_{S_2}]\longrightarrow0.
\]
Since $\omega_{S_2,5}\zeta_{T_2,5}\to0$, it follows that $\zeta_{S_2,5}\omega_{T_2,5}\to0$. The identity for the moment of inertia now yields $I\to0$, contradicting the positive constant inherited from a physical point.
\end{proof}

\section{The exclusion of Diagram 6}\label{sec:9}

We use the global sets in \eqref{eq:2.index_sets}. The two wings are $T_3$ and $S_3$, with outer pairs $S_2$ and $T_2$, respectively. We retain the notation $q_{ij}=q_i-q_j$.

AK's torque equations imply $m_1=m_2$ and $m_3=m_4$. Under these relations, the expansion of the reduced function excludes every singular sequence with Diagram~6.

\begin{theorem}\label{thm:9.1}
For positive masses, no singular sequence of complex, multiplier-one planar central configurations can have Albouy--Kaloshin diagram 6. In particular, diagram 6 cannot occur along the singular sequences supplied by Lemma~\ref{lem:2.3}.
\end{theorem}

We give a local analytic proof, including all signed distance branches. The signed distance variables satisfy $r_{ij}^2=z_{ij}w_{ij}$, with $r_{ij}\ne0$. As usual,
\[
U=\sum\limits_{i<j\atop i,j\in S_5}\frac{m_im_j}{r_{ij}},\qquad
I=\sum\limits_{i\in S_5}m_iz_iw_i,\qquad G=U+\frac12 I.
\]
On the centered configuration space, the critical equations of $G$ are precisely the multiplier-one central-configuration equations. All square roots below are local holomorphic branches on their specified sheets.

\subsection{Coordinates and distance formulas}\label{sec:9.1}

AK's calculation in Section 8.6 forces
\[
m_1=m_2>0,\qquad m_3=m_4>0,\qquad m_5>0.
\]
The centers of mass of the two outer pairs are $C_{S_2}$ and $C_{T_2}$, respectively, with the notation defined in Section~\ref{sec:2}.
The maximal $z$-edge $12$ and maximal $w$-edge $34$ give
\[
\begin{aligned}
Z_{12}&\approx z_{12}\approx\epsilon^{-2},&
W_{34}&\approx w_{34}\approx\epsilon^{-2},\\
r_{12},r_{34}&\approx1,&
w_{12},z_{34}&\approx\epsilon^2,\\
W_{12}&\approx\epsilon^2,&
Z_{34}&\approx\epsilon^2.
\end{aligned}
\]
Indeed, $Z_{12}=z_{12}/r_{12}^3$ and $W_{34}=w_{34}/r_{34}^3$ give the distance bounds, and $r_{ij}^2=z_{ij}w_{ij}$ then gives the complementary coordinate scales. The center estimates are $C_{S_2},C_{T_2},q_5=O(\epsilon^2)$ componentwise. For these estimates, \cite[Section~8.6, p.~576]{AK} obtains $z_5,w_5=O(\epsilon^2)$ from the maximal-wing estimates and the cluster equations (30); the small complementary coordinates in each maximal wing and the total center-of-mass equations then give the two pair-center estimates. Choose a local square root
\[
\mathfrak{a}=\left(\frac{z_{12}w_{34}}4\right)^{-1/2}.
\]
After a complex rotation as in \cite[p.~542]{AK}, with bounded rotation factor and inverse, we have $z_{12}=w_{34}=2/\mathfrak{a}$. For convenience, from this point onward we use the original symbols for the transformed positions, complex coordinates, and centers of mass. All ordered pairs below are written in their $(z,w)$ components. Thus $\mathfrak{a}\approx\epsilon^2$. The hats below denote the positions $q_{S_2,5},q_{T_2,5}$ of body 5 relative to the centers of mass of the two pairs, defined in Section~\ref{sec:2}, divided by $\mathfrak a$. The exact parametrization is
\begin{equation*}
\begin{aligned}
q_{12}&=\left(\frac2{\mathfrak{a}},\frac{\mathfrak{a}r_{12}^2}{2}\right),&
q_{34}&=\left(\frac{\mathfrak{a}r_{34}^2}{2},\frac2{\mathfrak{a}}\right),\\
\widehat{\zeta}_{S_2,5}&=\frac{\zeta_{S_2,5}}{\mathfrak{a}},&
\widehat{\omega}_{S_2,5}&=\frac{\omega_{S_2,5}}{\mathfrak{a}},\\
\widehat{\zeta}_{T_2,5}&=\frac{\zeta_{T_2,5}}{\mathfrak{a}},&
\widehat{\omega}_{T_2,5}&=\frac{\omega_{T_2,5}}{\mathfrak{a}},\\
q_5&=\frac{\mathfrak{a}}{M}\bigl(2m_1(\widehat{\zeta}_{S_2,5},\widehat{\omega}_{S_2,5})+2m_3(\widehat{\zeta}_{T_2,5},\widehat{\omega}_{T_2,5})\bigr).
\end{aligned}
\end{equation*}
The six variables $r_{12},r_{34},\widehat{\zeta}_{S_2,5},\widehat{\omega}_{S_2,5},\widehat{\zeta}_{T_2,5},\widehat{\omega}_{T_2,5}$ are bounded. Passing to a subsequence, they have limits. The six wing distances have finite nonzero limits in diagram 6; we work in a neighborhood of one such limiting point and its fixed signed branches. The maps $r_{12}\mapsto r_{12}^2/4$ and $r_{34}\mapsto r_{34}^2/4$ are locally invertible there, so the signed distances themselves may be used as coordinates.
Choose this neighborhood small enough that all wing squared distances and all scaled crossed squared distances $\mathfrak a^2r_{ij}^2$, $i\in S_2$, $j\in T_2$, stay away from zero; the latter tend to $-(-1)^{i+j}$ by \eqref{eq:9.3}. Their chosen local square roots are then holomorphic on a fixed polydisc about the limiting point.

For $\mathfrak{a}\ne0$, these seven variables are coordinates on the centered configuration space modulo complex rotations. Indeed $\mathfrak{a}$ is rotation invariant and the condition $z_{12}=w_{34}$ is transverse to the rotation orbit: the derivative of $z_{12}-w_{34}$ in the rotation direction at the slice is $4/\mathfrak{a}\ne0$. Rotational invariance of $G$ therefore shows that its seven critical equations in this chart are equivalent to the full centered critical equations.

The remaining wing squared distances are
\begin{equation}
\begin{aligned}
r_{15}^2&=(1-\mathfrak{a}^2\widehat{\zeta}_{S_2,5})\left(\frac{r_{12}^2}{4}-\widehat{\omega}_{S_2,5}\right),&
r_{25}^2&=(1+\mathfrak{a}^2\widehat{\zeta}_{S_2,5})\left(\frac{r_{12}^2}{4}+\widehat{\omega}_{S_2,5}\right),\\
r_{35}^2&=\left(\frac{r_{34}^2}{4}-\widehat{\zeta}_{T_2,5}\right)(1-\mathfrak{a}^2\widehat{\omega}_{T_2,5}),&
r_{45}^2&=\left(\frac{r_{34}^2}{4}+\widehat{\zeta}_{T_2,5}\right)(1+\mathfrak{a}^2\widehat{\omega}_{T_2,5}).
\end{aligned}
\label{eq:9.2}
\end{equation}
For the four crossed distances $r_{ij}$, use the actual body indices
\[
i\in S_2,\qquad
j\in T_2.
\]
The endpoint factors are $(-1)^{i+1}$ and $(-1)^{j+1}$, respectively. The crossed squared distances satisfy
{\footnotesize\begin{equation}
\begin{aligned}r_{ij}^2={}&-\frac{(-1)^{i+j}}{\mathfrak{a}^2}+\frac{r_{12}^2+r_{34}^2}{4}-(-1)^{i+1}(\widehat{\omega}_{S_2,5}-\widehat{\omega}_{T_2,5})+(-1)^{j+1}(\widehat{\zeta}_{S_2,5}-\widehat{\zeta}_{T_2,5})\\&{}+\mathfrak{a}^2\left(\widehat{\zeta}_{S_2,5}-\widehat{\zeta}_{T_2,5}+\frac{(-1)^{j+1} r_{34}^2}{4}\right)\left(\widehat{\omega}_{S_2,5}-\widehat{\omega}_{T_2,5}-\frac{(-1)^{i+1} r_{12}^2}{4}\right).\end{aligned}
\label{eq:9.3}
\end{equation}}
The exact chart and boundedness of the six remaining coordinates give
\[
z_{ij},w_{ij},r_{ij}\approx\mathfrak{a}^{-1}\approx\epsilon^{-2},
\qquad i\in S_2,\quad j\in T_2.
\]
Consequently, the four crossed interactions satisfy
\[
Z_{ij}=\frac{z_{ij}}{r_{ij}^3}\approx\epsilon^4,\qquad
W_{ij}=\frac{w_{ij}}{r_{ij}^3}\approx\epsilon^4,
\qquad i\in S_2,\quad j\in T_2.
\]
The moment of inertia is exactly
{\footnotesize\begin{equation}
\begin{aligned}
I={}&\frac{m_1r_{12}^2+m_3r_{34}^2}{2}
+\mathfrak{a}^2\Biggl[2m_1\widehat{\zeta}_{S_2,5}\widehat{\omega}_{S_2,5}
+2m_3\widehat{\zeta}_{T_2,5}\widehat{\omega}_{T_2,5}\\
&\qquad{}-\frac{(2m_1\widehat{\zeta}_{S_2,5}+2m_3\widehat{\zeta}_{T_2,5})(2m_1\widehat{\omega}_{S_2,5}+2m_3\widehat{\omega}_{T_2,5})}{M}\Biggr].
\end{aligned}
\label{eq:9.4}
\end{equation}}
Consequently $G$ extends holomorphically to $\mathfrak{a}=0$: the wing terms are analytic in $\mathfrak{a}^2$, and every crossed potential is $\mathfrak{a}$ times an analytic function of $\mathfrak{a}^2$ and the six other variables.

\subsection{Expansion of the crossed potential}\label{sec:9.2}

For each integer $i\ge0$, let $G_{(\mathfrak{a},i)}$ denote the coefficient of $\mathfrak{a}^i$ in the Taylor expansion of $G$ at $\mathfrak{a}=0$, with the other six coordinates $r_{12},r_{34},\widehat{\zeta}_{S_2,5},\widehat{\omega}_{S_2,5},\widehat{\zeta}_{T_2,5},\widehat{\omega}_{T_2,5}$ held fixed. Thus $G_{(\mathfrak{a},i)}=[\mathfrak{a}^i]G$ is a function of these six coordinates; in particular, $G_{(\mathfrak{a},0)}$ is the constant Taylor coefficient, not a critical-point evaluation.

Write $r_{ij}\sim\rho_{ij}/\mathfrak{a}$, with $\rho_{ij}^2=-(-1)^{i+j}$. The coefficient of $\mathfrak{a}$ in $G$ is $m_1m_3\sum\limits_{\substack{i\in S_2\\j\in T_2}}\rho_{ij}^{-1}$ and is independent of the other variables. Thus $\partial_{\mathfrak{a}}G=0$ along an approaching sequence forces this coefficient to vanish. The leading roots for $r_{13},r_{24}$ are $\pm\mathrm{i}$ and those for $r_{14},r_{23}$ are $\pm1$, so the real and imaginary parts cancel separately. In this section, let $\sigma_{13},\sigma_{14}\in\{1,-1\}$ label the leading branches of $r_{13}$ and $r_{14}$, respectively. The branches of $r_{24}$ and $r_{23}$ are then determined by the cancellations just obtained. The allowed sheets have the form
\begin{equation}
\rho_{13}=\mathrm{i}\sigma_{13},\quad \rho_{24}=-\mathrm{i}\sigma_{13},\quad
\rho_{14}=\sigma_{14},\quad \rho_{23}=-\sigma_{14},\quad \sigma_{13},\sigma_{14}\in\{1,-1\}.
\label{eq:9.5}
\end{equation}
To expand the crossed potential terms $m_1m_3/r_{ij}$, we take the reciprocal of the signed square root of the entire right-hand side of \eqref{eq:9.3}. Holding the other six coordinates $r_{12},r_{34},\widehat{\zeta}_{S_2,5},\widehat{\omega}_{S_2,5},\widehat{\zeta}_{T_2,5},\widehat{\omega}_{T_2,5}$ fixed, and factoring out $-(-1)^{i+j}/\mathfrak{a}^2$, gives
\[
\begin{aligned}
\frac{1}{r_{ij}}
={}&\frac{\mathfrak{a}}{\rho_{ij}}\Biggl[
1-(-1)^{i+j}\mathfrak{a}^2\Bigl(\frac{r_{12}^2+r_{34}^2}{4}\\
&\qquad{}-(-1)^{i+1}(\widehat{\omega}_{S_2,5}-\widehat{\omega}_{T_2,5})
+(-1)^{j+1}(\widehat{\zeta}_{S_2,5}-\widehat{\zeta}_{T_2,5})\Bigr)
+O(\mathfrak{a}^4)\Biggr]^{-1/2}.
\end{aligned}
\]
The factor inside the brackets tends to $1$. Its reciprocal square root is taken on the local branch with value $1$ at $\mathfrak{a}=0$, consistent with $\mathfrak{a}r_{ij}/\rho_{ij}\to1$. The reciprocal square root has value $1$ and derivative $-1/2$ at $1$, so Taylor expansion yields
\begin{equation}
\begin{aligned}
\frac{1}{r_{ij}}
={}&\frac{\mathfrak{a}}{\rho_{ij}}\\
&+\frac{\mathfrak{a}^3(-1)^{i+j}}{2\rho_{ij}}
\Bigl(\frac{r_{12}^2+r_{34}^2}{4}\\
&\qquad{}-(-1)^{i+1}(\widehat{\omega}_{S_2,5}-\widehat{\omega}_{T_2,5})+(-1)^{j+1}(\widehat{\zeta}_{S_2,5}-\widehat{\zeta}_{T_2,5})\Bigr)
+O(\mathfrak{a}^5).
\end{aligned}
\label{eq:9.7}
\end{equation}
Under \eqref{eq:9.5}, the four signed branches satisfy
\[
\begin{aligned}
\sum\limits_{\substack{i\in S_2\\j\in T_2}}\frac{(-1)^{j+1}}{\rho_{ij}}
&=-2(\sigma_{14}+\mathrm{i}\sigma_{13}),\\
\sum\limits_{\substack{i\in S_2\\j\in T_2}}\frac{(-1)^{i+1}}{\rho_{ij}}
&=2(\sigma_{14}-\mathrm{i}\sigma_{13}),\\
\sum\limits_{\substack{i\in S_2\\j\in T_2}}\frac{(-1)^{i+j}}{\rho_{ij}}
&=0.
\end{aligned}
\]
The wing potentials and the moment of inertia contain only even powers of $\mathfrak{a}$. Therefore, multiplying the cubic coefficient in \eqref{eq:9.7} by $m_1m_3$ and summing over the four crossed pairs $i\in S_2$, $j\in T_2$ gives
\[
\begin{aligned}
G_{(\mathfrak{a},3)}={}&\frac{m_1m_3}{2}
\sum\limits_{\substack{i\in S_2\\j\in T_2}}
\frac{(-1)^{i+j}}{\rho_{ij}}
\Bigl(\frac{r_{12}^2+r_{34}^2}{4}-(-1)^{i+1}(\widehat{\omega}_{S_2,5}-\widehat{\omega}_{T_2,5})\\
&\hspace{6em}{}+(-1)^{j+1}(\widehat{\zeta}_{S_2,5}-\widehat{\zeta}_{T_2,5})\Bigr)\\
={}&m_1m_3
\bigl\{(\sigma_{14}-\mathrm{i}\sigma_{13})(\widehat{\zeta}_{S_2,5}-\widehat{\zeta}_{T_2,5})\\
&\qquad{}+(\sigma_{14}+\mathrm{i}\sigma_{13})(\widehat{\omega}_{S_2,5}-\widehat{\omega}_{T_2,5})\bigr\}.
\end{aligned}
\]

\subsection{Critical points and the Hessian of \texorpdfstring{$G_{(\mathfrak a,0)}$}{G\_(a,0)}}\label{sec:9.3}

Put $\mathbf{y}=(r_{12},\widehat{\omega}_{S_2,5},r_{34},\widehat{\zeta}_{T_2,5})$ and $\mathbf{x}=(\widehat{\zeta}_{S_2,5},\widehat{\omega}_{T_2,5})$, and write the coefficients as $G_{(\mathfrak{a},i)}(\mathbf{y},\mathbf{x})$. Denote the limits of these coordinate groups along the selected subsequence by $\mathbf{y}^0=\lim\limits_{\epsilon\to0}\mathbf{y}$ and $\mathbf{x}^0=\lim\limits_{\epsilon\to0}\mathbf{x}$; along the same subsequence, $\mathfrak{a}\to0$. The constant coefficient $G_{(\mathfrak{a},0)}(\mathbf{y},\mathbf{x})=G(0,\mathbf{y},\mathbf{x})$ is independent of $\mathbf{x}$, and is therefore written as $G_{(\mathfrak{a},0)}(\mathbf{y})$. After \eqref{eq:9.5},
\begin{equation}
G(\mathfrak{a},\mathbf{y},\mathbf{x})=G_{(\mathfrak{a},0)}(\mathbf{y})+\mathfrak{a}^2G_{(\mathfrak{a},2)}(\mathbf{y},\mathbf{x})+\mathfrak{a}^3G_{(\mathfrak{a},3)}(\mathbf{y},\mathbf{x})+O(\mathfrak{a}^4).
\label{eq:9.8}
\end{equation}
In the formulas for $G_{(\mathfrak{a},0)}$ and $G_{(\mathfrak{a},2)}$ below, the wing distances $r_{ij}$ and their inverse cubes $\delta_{ij}$ are evaluated on the chart $\mathfrak{a}=0$, with their chosen signed branches. The contribution from $T_3$ to $G_{(\mathfrak{a},0)}$ is
\[
\frac{m_1^2}{r_{12}}+\frac{m_1m_5}{r_{15}}+\frac{m_1m_5}{r_{25}}+\frac{m_1r_{12}^2}{4},\qquad
r_{15}^2=\frac{r_{12}^2}{4}-\widehat{\omega}_{S_2,5},\qquad r_{25}^2=\frac{r_{12}^2}{4}+\widehat{\omega}_{S_2,5}.
\]
The two limiting critical equations $\partial_{r_{12}}G_{(\mathfrak{a},0)}=\partial_{\widehat{\omega}_{S_2,5}}G_{(\mathfrak{a},0)}=0$ are
\begin{equation}
\delta_{15}=\delta_{25},\qquad 2m_1\delta_{12}+m_5\delta_{15}=1.
\label{eq:9.9}
\end{equation}
For $S_3$, the corresponding formulas hold with $m_1,\allowbreak r_{12},\allowbreak\widehat{\omega}_{S_2,5}$ replaced by $m_3,\allowbreak r_{34},\allowbreak\widehat{\zeta}_{T_2,5}$; its limiting critical equations give $\delta_{35}=\delta_{45}$ and $2m_3\delta_{34}+m_5\delta_{35}=1$. The equality $\delta_{15}=\delta_{25}$ gives $(r_{25}/r_{15})^3=1$. Together with $r_{12}^2=2(r_{15}^2+r_{25}^2)$, this implies
\[
\left(\frac{\delta_{15}}{\delta_{12}}\right)^2
=8\left[1+\left(\frac{r_{25}}{r_{15}}\right)^2\right]^3
\in\{64,-8\}.
\]
The calculation for $S_3$ gives the same two possible squared values. Hence
\begin{equation}
\frac{\delta_{15}}{\delta_{12}},\qquad
\frac{\delta_{35}}{\delta_{34}}
\in\{8,-8,2\sqrt2\,\mathrm{i},-2\sqrt2\,\mathrm{i}\}.
\label{eq:9.11}
\end{equation}
All wing distances are finite and nonzero at an admissible limiting point, so $\delta_{12}\delta_{34}\ne0$. For $T_3$ at the limiting point, $\widehat{\omega}_{S_2,5}=0$ precisely when $r_{15}=r_{25}$, or equivalently $\delta_{15}/\delta_{12}=\pm8$. The remaining two values have $\widehat{\omega}_{S_2,5}\ne0$. Analogous statements hold for $\widehat{\zeta}_{T_2,5}$ and $\delta_{35}/\delta_{34}$. These statements concern only the limiting rescaled coordinates relative to the centers of mass of the corresponding pairs; no symmetry of the actual configurations is being imposed.

At a limiting critical point $\mathbf{y}^0$ of $G_{(\mathfrak{a},0)}$, define the two $2\times2$ Hessian blocks by
\[
\begin{aligned}
H^{0}_{r_{12},\widehat{\omega}_{S_2,5}}&:=\left.D^2_{(r_{12},\widehat{\omega}_{S_2,5})}G_{(\mathfrak{a},0)}\right|_{\mathbf{y}=\mathbf{y}^0},\\
H^{0}_{r_{34},\widehat{\zeta}_{T_2,5}}&:=\left.D^2_{(r_{34},\widehat{\zeta}_{T_2,5})}G_{(\mathfrak{a},0)}\right|_{\mathbf{y}=\mathbf{y}^0}.
\end{aligned}
\]
In particular,
\[
H^{0}_{r_{12},\widehat{\omega}_{S_2,5}}=\left.
\begin{pmatrix}
\partial_{r_{12}}^2G_{(\mathfrak{a},0)} & \partial_{r_{12}}\partial_{\widehat{\omega}_{S_2,5}}G_{(\mathfrak{a},0)}\\
\partial_{\widehat{\omega}_{S_2,5}}\partial_{r_{12}}G_{(\mathfrak{a},0)} & \partial_{\widehat{\omega}_{S_2,5}}^2G_{(\mathfrak{a},0)}
\end{pmatrix}
\right|_{\mathbf{y}=\mathbf{y}^0}.
\]
Since the contribution from $S_3$ to $G_{(\mathfrak{a},0)}$ is independent of $(r_{12},\widehat{\omega}_{S_2,5})$, $H^{0}_{r_{12},\widehat{\omega}_{S_2,5}}$ is exactly the Hessian of the contribution from $T_3$ displayed above, evaluated at $\mathbf{y}^0$. The derivatives are taken before imposing \eqref{eq:9.9}, with $r_{15}$ and $r_{25}$ differentiated through their squared-distance relations on the chosen signed branches. Evaluating at $\mathbf{y}^0$ and using \eqref{eq:9.9} gives the entries
\[
\begin{aligned}
(H^{0}_{r_{12},\widehat{\omega}_{S_2,5}})_{11}&=3m_1^2\delta_{12}+\frac{3m_1m_5r_{12}^2}{16}(r_{15}^{-5}+r_{25}^{-5}),\\
(H^{0}_{r_{12},\widehat{\omega}_{S_2,5}})_{22}&=\frac{3m_1m_5}{4}(r_{15}^{-5}+r_{25}^{-5}),\\
(H^{0}_{r_{12},\widehat{\omega}_{S_2,5}})_{12}=(H^{0}_{r_{12},\widehat{\omega}_{S_2,5}})_{21}
&=-\frac{3m_1m_5r_{12}}{8}(r_{15}^{-5}-r_{25}^{-5}).
\end{aligned}
\]
These formulas include the factors $\partial_{r_{12}}(r_{12}^2/4)=r_{12}/2$ from the signed-distance coordinate. The first derivatives vanish at the critical point. Using $r_{12}^2=2(r_{15}^2+r_{25}^2)$ and \eqref{eq:9.9}, its determinant is
\begin{equation*}
\det H^{0}_{r_{12},\widehat{\omega}_{S_2,5}}=\frac{9m_1^2m_5\delta_{15}r_{12}^2}{16r_{15}^2r_{25}^2}\ne0.
\end{equation*}
Similarly,
\[
\det H^{0}_{r_{34},\widehat{\zeta}_{T_2,5}}=\frac{9m_3^2m_5\delta_{35}r_{34}^2}{16r_{35}^2r_{45}^2}\ne0.
\]
Both blocks are invertible, so the block-diagonal Hessian $\left.D_{\mathbf{y}}^2G_{(\mathfrak{a},0)}\right|_{\mathbf{y}=\mathbf{y}^0}$ is invertible.

\subsection{Eliminating the six bounded variables}\label{sec:9.4}

We use the convergent subsequence and fixed signed neighborhood constructed in Section~\ref{sec:9.1}. The bounded rotation factor and its inverse preserve the coordinate bounds, and all root factors in the analytic extensions below stay nonzero there.

At the limiting critical point $\mathbf{y}^0$ of $G_{(\mathfrak{a},0)}$, define the scalar coefficients $\gamma_{11},\gamma_{22}$ and $\gamma_{12}=\gamma_{21}$ of the bilinear expression \eqref{eq:9.14}. The indices $1$ and $2$ label the pairs $S_2$ and $T_2$, respectively. These coefficients depend on the masses and the chosen point $\mathbf{y}^0$; the values $\delta_{15},\delta_{35}$ in \eqref{eq:9.13} are taken at $\mathbf{y}^0$ on the chart $\mathfrak{a}=0$:
\begin{equation}
\begin{aligned}
\gamma_{11}&=m_1\left(\frac{2m_3+m_5}{M}-m_5\delta_{15}\right),\\
\gamma_{22}&=m_3\left(\frac{2m_1+m_5}{M}-m_5\delta_{35}\right),\qquad \gamma_{12}=-\frac{2m_1m_3}{M}.
\end{aligned}
\label{eq:9.13}
\end{equation}
The associated symmetric coefficient matrix is
\[
\Gamma_2:=\begin{pmatrix}\gamma_{11}&\gamma_{12}\\\gamma_{12}&\gamma_{22}\end{pmatrix}.
\]
Expansion of \eqref{eq:9.2} and \eqref{eq:9.4}, using the critical-point relations $\partial_{\mathbf{y}}G_{(\mathfrak{a},0)}(\mathbf{y}^0)=0$, gives
\begin{equation}
G_{(\mathfrak{a},2)}(\mathbf{y}^0,\mathbf{x})=\gamma_{11}\widehat{\zeta}_{S_2,5}\widehat{\omega}_{S_2,5}^0+\gamma_{12}(\widehat{\zeta}_{S_2,5}\widehat{\omega}_{T_2,5}+\widehat{\zeta}_{T_2,5}^0\widehat{\omega}_{S_2,5}^0)+\gamma_{22}\widehat{\zeta}_{T_2,5}^0\widehat{\omega}_{T_2,5}.
\label{eq:9.14}
\end{equation}
The $\mathfrak{a}^2$ coefficient of $I/2$ in \eqref{eq:9.4} contributes $m_1(2m_3+m_5)/M$ to $\gamma_{11}$, $m_3(2m_1+m_5)/M$ to $\gamma_{22}$, and $-2m_1m_3/M$ to each mixed term. The wing potentials in \eqref{eq:9.2} contribute $-m_1m_5\delta_{15}$ and $-m_3m_5\delta_{35}$ to $\gamma_{11}$ and $\gamma_{22}$, respectively, using $\delta_{15}=\delta_{25}$ and $\delta_{35}=\delta_{45}$ at $\mathbf{y}^0$. The crossed potentials in \eqref{eq:9.7} have no $\mathfrak{a}^2$ term.

Before imposing the critical-point relations $\partial_{\mathbf{y}}G_{(\mathfrak{a},0)}=0$, the contribution from $T_3$ to $G_{(\mathfrak{a},2)}$ is
\[
\frac{m_1m_5}{2}\widehat{\zeta}_{S_2,5}\left\{\left(\frac{r_{12}^2}{4}-\widehat{\omega}_{S_2,5}\right)\delta_{15}+\left(-\frac{r_{12}^2}{4}-\widehat{\omega}_{S_2,5}\right)\delta_{25}\right\};
\]
its reduction to $-m_1m_5\delta_{15}\widehat{\zeta}_{S_2,5}\widehat{\omega}_{S_2,5}$ uses $\delta_{15}=\delta_{25}$ at $\mathbf{y}=\mathbf{y}^0$. Equation~\eqref{eq:9.14} gives values and derivatives with respect to $\mathbf{x}$ at $\mathbf{y}=\mathbf{y}^0$; derivatives with respect to $\mathbf{y}$ must be taken from the unreduced coefficient.

The six equations
\begin{equation}
\partial_{\mathbf{y}}G=0,\qquad \mathfrak{a}^{-2}\partial_{\widehat{\zeta}_{S_2,5}}G=0,\qquad \mathfrak{a}^{-2}\partial_{\widehat{\omega}_{T_2,5}}G=0
\label{eq:9.15}
\end{equation}
extend analytically to $\mathfrak{a}=0$. Since $G_{(\mathfrak{a},0)}$ is independent of $\mathbf{x}=(\widehat{\zeta}_{S_2,5},\widehat{\omega}_{T_2,5})$, \eqref{eq:9.8} gives
\[
\mathfrak{a}^{-2}\partial_{\mathbf{x}}G
=\partial_{\mathbf{x}}G_{(\mathfrak{a},2)}+\mathfrak{a}\,\partial_{\mathbf{x}}G_{(\mathfrak{a},3)}+O(\mathfrak{a}^2).
\]
Taking $\mathfrak{a}\to0$, $\mathbf{y}\to \mathbf{y}^0$ and $\mathbf{x}\to \mathbf{x}^0=(\widehat{\zeta}_{S_2,5}^0,\widehat{\omega}_{T_2,5}^0)$ in the last two equations of \eqref{eq:9.15} therefore gives $\partial_{\mathbf{x}}G_{(\mathfrak{a},2)}(\mathbf{y}^0,\mathbf{x}^0)=0$. Differentiating \eqref{eq:9.14} with respect to $\widehat{\zeta}_{S_2,5}$ and $\widehat{\omega}_{T_2,5}$ yields
\[
\begin{aligned}
\gamma_{11}\widehat{\omega}_{S_2,5}^0+\gamma_{12}\widehat{\omega}_{T_2,5}^0&=0,\\
\gamma_{12}\widehat{\zeta}_{S_2,5}^0+\gamma_{22}\widehat{\zeta}_{T_2,5}^0&=0.
\end{aligned}
\]
Since $\gamma_{12}=-2m_1m_3/M\ne0$, solving these two linear equations gives
\begin{equation}
\widehat{\omega}_{T_2,5}^0=-\frac{\gamma_{11}}{\gamma_{12}}\widehat{\omega}_{S_2,5}^0,\qquad \widehat{\zeta}_{S_2,5}^0=-\frac{\gamma_{22}}{\gamma_{12}}\widehat{\zeta}_{T_2,5}^0.
\label{eq:9.16}
\end{equation}
The Jacobian of \eqref{eq:9.15} with respect to $(\mathbf{y},\mathbf{x})$ at this point is block lower triangular. The block corresponding to $\mathbf{y}$ is the invertible matrix $D_{\mathbf{y}}^2G_{(\mathfrak{a},0)}$, and the block corresponding to $\mathbf{x}$ is
\[
\begin{pmatrix}0&\gamma_{12}\\\gamma_{12}&0\end{pmatrix},
\]
which is invertible because $m_1,m_3>0$. The analytic implicit-function theorem therefore gives a unique analytic solution $(\mathbf{y}(\mathfrak{a}),\mathbf{x}(\mathfrak{a}))$ to these six equations near the limiting point. Every approaching sequence of full central configurations must lie on this solution. No assumption that its coordinates initially have integral Puiseux exponents is required.

Define the reduced function $\tilde{G}(\mathfrak{a})=G(\mathfrak{a},\mathbf{y}(\mathfrak{a}),\mathbf{x}(\mathfrak{a}))$. The remaining full critical equation is $\tilde{G}'(\mathfrak{a})=0$, since the six other derivatives vanish. We have $\mathbf{y}(\mathfrak{a})=\mathbf{y}^0+O(\mathfrak{a}^2)$ and $\mathbf{x}(\mathfrak{a})=\mathbf{x}^0+O(\mathfrak{a})$. Stationarity of $G_{(\mathfrak{a},2)}$ in $\mathbf{x}$ at $(\mathbf{y}^0,\mathbf{x}^0)$ yields
\begin{equation}
\begin{aligned}
\tilde{G}_{(\mathfrak{a},2)}&=-\frac{\det \Gamma_2}{\gamma_{12}}\widehat{\omega}_{S_2,5}^0\widehat{\zeta}_{T_2,5}^0,\\
\tilde{G}_{(\mathfrak{a},3)}&=m_1m_3
\{(\sigma_{14}-\mathrm{i}\sigma_{13})(\widehat{\zeta}_{S_2,5}^0-\widehat{\zeta}_{T_2,5}^0)+(\sigma_{14}+\mathrm{i}\sigma_{13})(\widehat{\omega}_{S_2,5}^0-\widehat{\omega}_{T_2,5}^0)\}.
\end{aligned}
\label{eq:9.17}
\end{equation}
Here $\tilde{G}_{(\mathfrak{a},j)}:=[\mathfrak{a}^j]\tilde{G}$, whereas $G_{(\mathfrak{a},j)}=[\mathfrak{a}^j]G$ is taken before the other coordinates are eliminated. The limiting critical equations for the two wings simplify the determinant of the matrix in \eqref{eq:9.13} to
\begin{equation}
\det \Gamma_2=\frac{4m_1^2m_3^2m_5}{M}\left(\delta_{12}\delta_{34}-\delta_{15}\delta_{34}-\delta_{35}\delta_{12}\right).
\label{eq:9.18}
\end{equation}
Since $\delta_{12}\delta_{34}\ne0$, vanishing of \eqref{eq:9.18} would require $\delta_{15}/\delta_{12}+\delta_{35}/\delta_{34}=1$. No pair of values in \eqref{eq:9.11} has sum one. Hence $\det \Gamma_2\ne0$.

\subsection{Exclusion by \texorpdfstring{$\tilde G'(\mathfrak a)\ne0$}{Gtilde'(a) is nonzero}}\label{sec:9.5}

If both limiting rescaled transverse coordinates are nonzero, then $\widehat{\omega}_{S_2,5}^0\widehat{\zeta}_{T_2,5}^0\ne0$ and \eqref{eq:9.17} gives $\tilde{G}_{(\mathfrak{a},2)}\ne0$. Thus $\tilde{G}'(\mathfrak{a})=2\tilde{G}_{(\mathfrak{a},2)}\mathfrak{a}+O(\mathfrak{a}^2)$ is nonzero for every sufficiently small $\mathfrak{a}\ne0$.

If $\widehat{\omega}_{S_2,5}^0\ne0$ and $\widehat{\zeta}_{T_2,5}^0=0$, then $\widehat{\zeta}_{S_2,5}^0=0$, $\tilde{G}_{(\mathfrak{a},2)}=0$, and
\begin{equation*}
\tilde{G}_{(\mathfrak{a},3)}=m_1m_3(\sigma_{14}+\mathrm{i}\sigma_{13})\widehat{\omega}_{S_2,5}^0\left(1+\frac{\gamma_{11}}{\gamma_{12}}\right).
\end{equation*}
Now $\gamma_{11}+\gamma_{12}=m_1m_5(1/M-\delta_{15})$. If this vanished, $\delta_{15}=1/M$ and \eqref{eq:9.9} would imply $\delta_{15}/\delta_{12}=m_1/(m_1+m_3)$, whereas $\widehat{\omega}_{S_2,5}^0\ne0$ makes this ratio nonzero purely imaginary. Therefore $\tilde{G}_{(\mathfrak{a},3)}\ne0$, and $\tilde{G}'(\mathfrak{a})$ is again nonzero for every sufficiently small $\mathfrak{a}\ne0$. The case $\widehat{\omega}_{S_2,5}^0=0$, $\widehat{\zeta}_{T_2,5}^0\ne0$ is identical, with
\[
\tilde{G}_{(\mathfrak{a},3)}=-m_1m_3(\sigma_{14}-\mathrm{i}\sigma_{13})\widehat{\zeta}_{T_2,5}^0\left(1+\frac{\gamma_{22}}{\gamma_{12}}\right),\qquad \gamma_{22}+\gamma_{12}=m_3m_5(1/M-\delta_{35})\ne0.
\]
It remains to consider $\widehat{\omega}_{S_2,5}^0=\widehat{\zeta}_{T_2,5}^0=0$. Then all four limiting rescaled coordinates of body 5 relative to the centers of mass of the pairs vanish by \eqref{eq:9.16}. At this fixed point the wing terms and the moment of inertia are exactly independent of $\mathfrak{a}$. The crossed squared distances coincide in opposite pairs, whose signed reciprocal distances cancel exactly by \eqref{eq:9.5}. Therefore
\begin{equation}
G(\mathfrak{a},\mathbf{y}^0,\mathbf{x}^0)=G_{(\mathfrak{a},0)}(\mathbf{y}^0)
\label{eq:9.20}
\end{equation}
identically for small $\mathfrak{a}$. Moreover $\partial_{\mathbf{y}}G_{(\mathfrak{a},2)}(\mathbf{y}^0,\mathbf{x}^0)=0$, using the unreduced coefficient of $G_{(\mathfrak{a},2)}$: the contributions from $T_3$ and $S_3$ contain, respectively, the rescaled relative $z$-coordinate $\widehat{\zeta}_{S_2,5}$ and the rescaled relative $w$-coordinate $\widehat{\omega}_{T_2,5}$ as factors. Both factors vanish at $(\mathbf{y}^0,\mathbf{x}^0)$, and the coefficient contributed by the moment of inertia is quadratic in the four rescaled relative coordinates. The solution for $\mathbf{y}$ consequently satisfies $\mathbf{y}(\mathfrak{a})=\mathbf{y}^0+O(\mathfrak{a}^3)$. The last two equations in \eqref{eq:9.15}, at order $\mathfrak{a}$, now give
\begin{equation}
\widehat{\zeta}_{S_2,5}(\mathfrak{a})=-\frac{M(\sigma_{14}+\mathrm{i}\sigma_{13})}{2}\mathfrak{a}+O(\mathfrak{a}^2),\qquad \widehat{\omega}_{T_2,5}(\mathfrak{a})=\frac{M(\sigma_{14}-\mathrm{i}\sigma_{13})}{2}\mathfrak{a}+O(\mathfrak{a}^2).
\label{eq:9.21}
\end{equation}
At order $\mathfrak{a}^4$, the only nonzero terms of the reduced function are
\[
\mathfrak{a}^2\gamma_{12}\widehat{\zeta}_{S_2,5}\widehat{\omega}_{T_2,5}+\mathfrak{a}^3(m_1m_3(\sigma_{14}-\mathrm{i}\sigma_{13})\widehat{\zeta}_{S_2,5}-m_1m_3(\sigma_{14}+\mathrm{i}\sigma_{13})\widehat{\omega}_{T_2,5}).
\]
Corrections to $\mathbf{y}$ change $G_{(\mathfrak{a},0)}$ only at order $\mathfrak{a}^6$; their contributions to $\mathfrak{a}^2G_{(\mathfrak{a},2)}$ and to the cubic term in \eqref{eq:9.8} are also of order at least $\mathfrak{a}^6$. Every explicit higher coefficient evaluated at the limiting point vanishes by \eqref{eq:9.20}. Substituting \eqref{eq:9.21} and $(\sigma_{14}-\mathrm{i}\sigma_{13})(\sigma_{14}+\mathrm{i}\sigma_{13})=2$ gives
\begin{equation*}
\tilde{G}(\mathfrak{a})=G_{(\mathfrak{a},0)}(\mathbf{y}^0)-m_1m_3M\mathfrak{a}^4+O(\mathfrak{a}^5),\qquad \tilde{G}'(\mathfrak{a})=-4m_1m_3M\mathfrak{a}^3+O(\mathfrak{a}^4).
\end{equation*}
Its leading coefficient is nonzero. The reduced critical equation has no sufficiently small nonzero solution in this final case either. These three cases exhaust all limiting wings, and prove Theorem~\ref{thm:9.1}.


\section{The exclusion of Diagram 2}
\label{sec:11}

We use the binary and triangle index sets $S_2$ and $S_3$ from \eqref{eq:2.index_sets}. In every crossed pair, $i\in S_2$ and $j\in S_3$.

The remaining critical equations exclude the non-Lagrange Euler limits. We treat separately the cases in which the partial Hessian is invertible or singular.

\begin{theorem}\label{thm:11.1}
For any five positive masses, no singular sequence on a normalized complex
component containing a physical configuration can have Albouy--Kaloshin diagram 2. More precisely,
no complex normalized singular sequence with positive masses has a diagram-2 limit consisting of
a binary and a non-Lagrange Euler triangle. The Lagrange alternative is excluded by the physical lower bound \eqref{eq:3.1} on the component.
\end{theorem}

\subsection{Equations for the limiting triangle}\label{sec:10.6}

AK Section 5.0 decouples diagram 2 into a normalized binary on $S_2$ and a normalized triangle on $S_3$; all six
cross distances diverge \cite[pp.~556--557,572--575]{AK}. The binary distance satisfies
$r_{12}^3\to m_1+m_2$. We use the following terminology for the internal limiting triangle, following AK \cite[p.~584]{AK}.

\begin{definition}[Lagrange and Euler cases]\label{def:lagrange_euler}
Consider a collision-free, centered complex three-body central configuration on $S_3$, normalized to multiplier one. In this definition, $q_j=(z_j,w_j)$ denotes its coordinates and $r_{jk}$ its nonzero distances on the specified signed branches.
It is in the \emph{Lagrange case} if
\[
r_{34}^3=r_{35}^3=r_{45}^3=e_1(S_3).
\]
It is in the \emph{Euler case} if its three points lie on one complex affine line, equivalently
\[
z_{35}w_{45}-z_{45}w_{35}=0.
\]
A \emph{non-Lagrange Euler configuration} is an Euler configuration that does not satisfy the displayed equal-cube condition.
\end{definition}

These cases are not disjoint over $\mathbb C$: a collinear complex configuration may also satisfy the Lagrange condition \cite[p.~584]{AK}. If the two relative position vectors are independent, the normalized equations, as in AK equation (45), force the Lagrange condition. Thus every normalized complex triangle belongs either to the Lagrange case (including its intersection with the Euler case) or to the non-Lagrange Euler case.
For diagram 2 these definitions apply to the centered, rescaled limiting triangle constructed in Section~\ref{sec:11.1}. Every Lagrange limit is excluded by Lemma~\ref{lem:3.2}, including a Lagrange--Euler intersection. Only non-Lagrange Euler limits remain to be considered.

For a collinear triangle its $w$ coordinates $w_3,w_4,w_5$ are pairwise distinct and nonzero in
diagram 2. Up to a common nonzero scalar, its four necessary signed Euler equations are
\begin{equation}
\mathcal{E}(w_3,w_4,w_5)=w_3w_{35}^2w_{34}^2\pm w_4w_{45}^2w_{34}^2
\pm w_5w_{45}^2w_{35}^2=0.
\label{eq:10.26}
\end{equation}
The two signs are chosen independently; for each choice, $\mathcal{E}$ denotes the corresponding polynomial.
Indeed AK equation (47) is $w_3W_{45}+w_4W_{53}+w_5W_{34}=0$; on a nonisotropic line each
inverse-square factor has an independent sign up to simultaneous reversal. Clearing its three
nonzero squared denominators gives \eqref{eq:10.26}. The masses enter through
$m_3w_3+m_4w_4+m_5w_5=0$. The all-plus polynomial contains the Lagrange--Euler intersection factor;
this is an actual Lagrange locus, not a factor discarded by a genericity assumption \cite[pp.~584--585]{AK}.
For centered real collinear configurations, see also \cite[Proposition~3]{AM}.

\subsection{Coordinates and rescaling}\label{sec:11.1}
\label{subsec:11.1}

Write $C_{S_2},C_{S_3}$ for the centers of mass of the binary and the triangle. In diagram 2, both binary $z$ coordinates and all three triangle $w$ coordinates are circled. The binary edge $12$ is a maximal $z$-edge, and every internal triangle edge is a maximal $w$-edge. The disconnected-block estimates in \cite[Section~5.0, pp.~556--557]{AK}, applied to the diagram described in \cite[Section~8.2, pp.~572--573]{AK}, give
\[
\begin{aligned}
Z_{12}&\approx z_{12}\approx\epsilon^{-2},&
W_{12}&\approx w_{12}\approx\epsilon^2,\\
W_{jk}&\approx w_{jk}\approx\epsilon^{-2},&
Z_{jk}&\approx z_{jk}\approx\epsilon^2,\\
r_{12},r_{jk}&\approx1,&&j,k\in S_3,\quad j\ne k.
\end{aligned}
\]
Every crossed pair satisfies $z_{ij},w_{ij},r_{ij}\approx\epsilon^{-2}$, so its individual force components have the two-sided estimates
\[
Z_{ij}=\frac{z_{ij}}{r_{ij}^3}\approx\epsilon^4,\qquad
W_{ij}=\frac{w_{ij}}{r_{ij}^3}\approx\epsilon^4,
\qquad i\in S_2,\quad j\in S_3.
\]
Summing the exact mass-weighted normalized equations within each cluster cancels every internal interaction and gives
\[
\begin{aligned}
(m_1+m_2)C_{S_2,z}
&=\sum\limits_{\substack{i\in S_2\\j\in S_3}}m_im_jZ_{ij}=O(\epsilon^4),\\
(m_1+m_2)C_{S_2,w}
&=\sum\limits_{\substack{i\in S_2\\j\in S_3}}m_im_jW_{ij}=O(\epsilon^4),\\
e_1(S_3)C_{S_3,z}&=-(m_1+m_2)C_{S_2,z},\\
e_1(S_3)C_{S_3,w}&=-(m_1+m_2)C_{S_2,w}.
\end{aligned}
\]
Consequently,
\begin{equation}
C_{S_2},C_{S_3}=O(\epsilon^4)
\label{eq:11.1}
\end{equation}
componentwise. These are upper bounds: the complex force components can cancel in the weighted sums even for fixed positive masses, so individual two-sided estimates do not give two-sided center estimates.

We use the relative coordinates $\zeta_{S_3,k},\omega_{S_3,k}$, $k\in S_5$, defined in Section~\ref{sec:2}.
The maximal binary edge gives $z_{21}\approx\epsilon^{-2}$. Since vertex 5 is $w$-circled, $w_5\approx\epsilon^{-2}$; combining this point-coordinate estimate with \eqref{eq:11.1} yields $\omega_{S_3,5}=w_5-C_{S_3,w}\sim w_5$, and hence $\omega_{S_3,5}\approx\epsilon^{-2}$.
Choose a local square root $\mathfrak{c}=(z_{21}\omega_{S_3,5})^{-1/2}$ and a complex rotation
to impose $z_{21}=\omega_{S_3,5}=\mathfrak{c}^{-1}$. The rotation factor and its inverse are bounded, and $\mathfrak{c}\approx\epsilon^2$. The slice is
transverse to rotation because the derivative of $z_{21}-\omega_{S_3,5}$ in the rotation direction is $2/\mathfrak{c}\ne0$.

For the triangle, use the rescaled relative coordinates
\[
\widehat{\zeta}_{S_3,j}=\frac{\zeta_{S_3,j}}{\mathfrak{c}},\qquad
\widehat{\omega}_{S_3,j}=\mathfrak{c}\omega_{S_3,j},\qquad j\in S_3.
\]
The normalization above gives $\widehat{\omega}_{S_3,5}=1$. We take $\widehat{\omega}_{S_3,3}$ as the independent triangle $w$ coordinate; $\widehat{\omega}_{S_3,4}$ is determined by the center-of-mass condition.
Use $r_{12}^2=z_{21}w_{21}$ as the binary coordinate, with its signed square root $r_{12}$ fixed locally.
Define the rescaled separation of the centers of mass by
\[
\nu=(\nu_z,\nu_w)=\frac{C_{S_2}-C_{S_3}}{\mathfrak{c}^2}.
\]
Using the original body labels, write the exact chart as
\begin{equation}
\begin{aligned}
\zeta_{S_2,1}&=-\frac{m_2}{(m_1+m_2)\mathfrak{c}},&
\omega_{S_2,1}&=-\frac{m_2\mathfrak{c}r_{12}^2}{m_1+m_2},\\
\zeta_{S_2,2}&=\frac{m_1}{(m_1+m_2)\mathfrak{c}},&
\omega_{S_2,2}&=\frac{m_1\mathfrak{c}r_{12}^2}{m_1+m_2},\\
\zeta_{S_3,j}&=\mathfrak{c}\widehat{\zeta}_{S_3,j},&
\omega_{S_3,j}&=\widehat{\omega}_{S_3,j}/\mathfrak{c} &&(j\in S_3),\\
\widehat{\omega}_{S_3,4}&=-\frac{m_3\widehat{\omega}_{S_3,3}+m_5}{m_4},&
\widehat{\zeta}_{S_3,5}&=-\frac{m_3\widehat{\zeta}_{S_3,3}+m_4\widehat{\zeta}_{S_3,4}}{m_5},\\
C_{S_2}&=\frac{e_1(S_3)}{M}\mathfrak{c}^2(\nu_z,\nu_w),&
C_{S_3}&=-\frac{m_1+m_2}{M}\mathfrak{c}^2(\nu_z,\nu_w).
\end{aligned}
\label{eq:11.2}
\end{equation}
The triangle coordinates retain their indices in $S_3$; every crossed pair has $i\in S_2$ and $j\in S_3$.
Write $\widehat{\zeta}_{S_3,T_2}=(\widehat{\zeta}_{S_3,3},\widehat{\zeta}_{S_3,4})$.
Derivative subscripts $r_{12}^2$, $\widehat{\zeta}_{S_3,T_2}$, and $\widehat{\omega}_{S_3,3}$ refer to the squared-distance coordinate, the two independent triangle $z$ coordinates, and the independent triangle $w$ coordinate, respectively. Primes denote differentiation with respect to $\widehat{\omega}_{S_3,3}$. There are seven quotient coordinates $(\mathfrak{c},r_{12}^2,\widehat{\omega}_{S_3,3},\widehat{\zeta}_{S_3,3},\widehat{\zeta}_{S_3,4},\nu_z,\nu_w)$,
as required after translation and rotation. They invert to the actual positions for $\mathfrak{c}\ne0$. All six
coordinates other than $\mathfrak{c}$ are bounded by the diagram estimates and \eqref{eq:11.1}. Pass to a convergent
subsequence of them. Every internal wing squared distance has a nonzero finite limit; each $\widehat{\omega}_{S_3,j}$ and
each $\widehat{\omega}_{S_3,j}-\widehat{\omega}_{S_3,k}$ likewise has a nonzero limit, by the circles and maximal triangle edges. Thus the
chart discards no diagram-2 leading point. On a sufficiently small polydisc about this limit, the internal squared distances and the scaled crossed squared distances $\mathfrak c^2r_{ij}^2$ stay away from zero. Fix their local holomorphic square roots on the inherited signed sheet. All subsequent Taylor estimates and local inverse or implicit functions are taken in this same neighborhood.

The reciprocal rescaling preserves the internal squared distances, and the crossed force terms tend to zero in the rescaled equations. Hence the limiting triangle is a collision-free normalized three-body configuration. With superscript $0$ denoting the limit along this subsequence as $\epsilon\to0$, Definition~\ref{def:lagrange_euler} gives the following concrete tests:
\[
\begin{aligned}
\text{Lagrange case:}\quad &(r_{34}^0)^3=(r_{35}^0)^3=(r_{45}^0)^3=e_1(S_3),\\
\text{Euler case:}\quad &\widehat{\zeta}_{S_3,j}^0=\widehat{\zeta}_{S_3,5}^0\widehat{\omega}_{S_3,j}^0\qquad(j\in S_3).
\end{aligned}
\]
The second test uses the center-of-mass condition and $\widehat{\omega}_{S_3,5}^0=1$; collision-freeness then gives $\widehat{\zeta}_{S_3,5}^0\ne0$. Both tests may hold simultaneously. The classification concerns the limiting triangle; the internal triangle at a nonzero $\mathfrak c$ is still perturbed by the other cluster.

On a fixed signed-distance sheet, define the holomorphic function
\[
G=U+\frac12I.
\]
Its critical equations on the centered space are the multiplier-one equations. Rotation invariance and
the transverse slice show that its seven partial derivatives in \eqref{eq:11.2} are exactly the required critical
equations. No independent internal phase is identified with another: the binary distance, the three
triangle distances, and the six crossed distances retain their chosen local signed square roots.

\subsection{Expansion of \texorpdfstring{$G$}{G} at \texorpdfstring{$\mathfrak c=0$}{c=0}}\label{sec:11.2}
\label{subsec:11.2}
The moment of inertia is exactly
\begin{equation*}
I=\frac{m_1m_2}{m_1+m_2}r_{12}^2+\sum\limits_{j\in S_3}m_j\widehat{\zeta}_{S_3,j}\widehat{\omega}_{S_3,j}
+\frac{(m_1+m_2)e_1(S_3)}{M}\mathfrak{c}^4\nu_z\nu_w.
\end{equation*}
Let $G_{(\mathfrak{c},0)}$ be the sum of the contributions to $G$ from the binary and the triangle:
\begin{equation}
\begin{aligned}
G_{(\mathfrak{c},0)}={}&\frac{m_1m_2}{r_{12}}+\frac{m_1m_2}{2(m_1+m_2)}r_{12}^2\\
&+\sum\limits_{\substack{j,k\in S_3\\j<k}}\frac{m_jm_k}{\sqrt{(\widehat{\zeta}_{S_3,j}-\widehat{\zeta}_{S_3,k})(\widehat{\omega}_{S_3,j}-\widehat{\omega}_{S_3,k})}}\\
&+\frac12\sum\limits_{j\in S_3}m_j\widehat{\zeta}_{S_3,j}\widehat{\omega}_{S_3,j}.
\end{aligned}
\label{eq:11.4}
\end{equation}
Every square root in this formula is the locally specified branch. For a crossed pair the exact squared
distance is
\begin{equation}
\begin{aligned}
r_{1j}^2&=\left(-\frac{m_2}{(m_1+m_2)\mathfrak{c}}-\mathfrak{c}\widehat{\zeta}_{S_3,j}+\mathfrak{c}^2\nu_z\right)
\left(-\frac{m_2\mathfrak{c}r_{12}^2}{m_1+m_2}-\frac{\widehat{\omega}_{S_3,j}}{\mathfrak{c}}+\mathfrak{c}^2\nu_w\right),\\
r_{2j}^2&=\left(\frac{m_1}{(m_1+m_2)\mathfrak{c}}-\mathfrak{c}\widehat{\zeta}_{S_3,j}+\mathfrak{c}^2\nu_z\right)
\left(\frac{m_1\mathfrak{c}r_{12}^2}{m_1+m_2}-\frac{\widehat{\omega}_{S_3,j}}{\mathfrak{c}}+\mathfrak{c}^2\nu_w\right),\qquad j\in S_3.
\end{aligned}
\label{eq:11.5}
\end{equation}
The coefficients of $\mathfrak{c}^{-2}$ are $m_2\widehat{\omega}_{S_3,j}/(m_1+m_2)$ and $-m_1\widehat{\omega}_{S_3,j}/(m_1+m_2)$, respectively, and are nonzero. Factoring $\mathfrak{c}^{-2}$ in \eqref{eq:11.5} shows that each reciprocal distance
is $\mathfrak{c}$ times a holomorphic function whose first correction has order $\mathfrak{c}^2$.
For each integer $i\ge0$, let $G_{(\mathfrak{c},i)}$ denote the coefficient of $\mathfrak{c}^i$ in the Taylor expansion of $G$ at $\mathfrak{c}=0$, with the other six coordinates $r_{12}^2,\widehat{\zeta}_{S_3,3},\widehat{\zeta}_{S_3,4},\widehat{\omega}_{S_3,3},\nu_z,\nu_w$ held fixed. Thus $G_{(\mathfrak{c},i)}=[\mathfrak{c}^i]G$; the subscript records both the local parameter and its power, as in Section~\ref{sec:9}. This gives
\begin{equation}
\begin{aligned}
G={}&G_{(\mathfrak{c},0)}(r_{12}^2,\widehat{\zeta}_{S_3,T_2},\widehat{\omega}_{S_3,3})+\mathfrak{c}G_{(\mathfrak{c},1)}(\widehat{\omega}_{S_3,3})\\
&+\mathfrak{c}^3G_{(\mathfrak{c},3)}(r_{12}^2,\widehat{\zeta}_{S_3,T_2},\widehat{\omega}_{S_3,3})\\
&+\mathfrak{c}^4G_{(\mathfrak{c},4)}(r_{12}^2,\widehat{\zeta}_{S_3,T_2},\widehat{\omega}_{S_3,3},\nu)+O(\mathfrak{c}^5).
\end{aligned}
\label{eq:11.6}
\end{equation}
In particular, $G_{(\mathfrak{c},2)}=0$. The expansion is holomorphic in all seven coordinates on
a fixed neighborhood of the limiting point. Its remainder estimates are uniform there and remain
valid after taking the partial derivatives used below. Bounded center variables therefore cause no
hidden error term of lower order.

Choose fixed square roots of $m_2/(m_1+m_2)$ and $-m_1/(m_1+m_2)$, and local square roots of $\widehat{\omega}_{S_3,j}$. For cross signs $\sigma_{ij}\in\{1,-1\}$, $i\in S_2$, $j\in S_3$, the
leading coefficient is
\begin{equation}
G_{(\mathfrak{c},1)}(\widehat{\omega}_{S_3,3})=\sum\limits_{j\in S_3}\frac{m_j\gamma_j}{\sqrt{\widehat{\omega}_{S_3,j}}},\qquad
\gamma_j=\frac{m_1\sigma_{1j}}{\sqrt{m_2/(m_1+m_2)}}+\frac{m_2\sigma_{2j}}{\sqrt{-m_1/(m_1+m_2)}}.
\label{eq:11.7}
\end{equation}
Define the positive real constants
\[
\gamma_{\mathrm{R}}=m_1\sqrt{\frac{m_1+m_2}{m_2}},\qquad
\gamma_{\mathrm{I}}=m_2\sqrt{\frac{m_1+m_2}{m_1}}.
\]
They are the absolute values of the real and imaginary parts of each $\gamma_j$, respectively. Thus each $\gamma_j$
belongs to $\{\gamma_{\mathrm{R}}+\mathrm{i}\gamma_{\mathrm{I}},\gamma_{\mathrm{R}}-\mathrm{i}\gamma_{\mathrm{I}},-\gamma_{\mathrm{R}}+\mathrm{i}\gamma_{\mathrm{I}},-\gamma_{\mathrm{R}}-\mathrm{i}\gamma_{\mathrm{I}}\}$. In particular
\begin{equation}
|\gamma_3|=|\gamma_4|=|\gamma_5|>0,\qquad
\gamma_j^2\in\{(\gamma_{\mathrm{R}}+\mathrm{i}\gamma_{\mathrm{I}})^2,(\gamma_{\mathrm{R}}-\mathrm{i}\gamma_{\mathrm{I}})^2\}.
\label{eq:11.8}
\end{equation}
This observation includes all 64 crossed sign choices. Write $\widehat{\omega}_{S_3,3}^0=\lim\limits_{\epsilon\to0}\widehat{\omega}_{S_3,3}$. The equation $\partial_{\mathfrak{c}}G=0$ requires $G_{(\mathfrak{c},1)}(\widehat{\omega}_{S_3,3}^0)=0$ at
the limit; this is the limiting constraint on the crossed terms, corresponding to the constraint in AK Section 5.0 \cite[p.~557]{AK}, on the specified radical sheet.

\subsection{Simple zeros of \texorpdfstring{$G_{(\mathfrak{c},1)}$}{G(c,1)} at non-Lagrange limits}\label{sec:11.3}
\label{subsec:11.3}
A multiple zero means that both $G_{(\mathfrak c,1)}$ and its derivative with respect to $\widehat{\omega}_{S_3,3}$ vanish. The following lemma shows that such a zero can occur at an Euler limit only in the Lagrange--Euler intersection of Definition~\ref{def:lagrange_euler}.

\begin{lemma}\label{lem:11.2}
Consider a non-Lagrange signed Euler-triangle limit with positive masses. If $G_{(\mathfrak{c},1)}(\widehat{\omega}_{S_3,3}^0)=0$, then $G'_{(\mathfrak{c},1)}(\widehat{\omega}_{S_3,3}^0)\ne0$.
\end{lemma}
\begin{proof}
Suppose both vanish. Differentiating \eqref{eq:11.7} and using $\widehat{\omega}_{S_3,3}'=1$, $\widehat{\omega}_{S_3,4}'=-m_3/m_4$ gives
\[
\frac{\gamma_3}{\widehat{\omega}_{S_3,3}^{3/2}}=\frac{\gamma_4}{\widehat{\omega}_{S_3,4}^{3/2}}.
\]
Inserting this relation in $G_{(\mathfrak{c},1)}=0$ and using $m_3\widehat{\omega}_{S_3,3}+m_4\widehat{\omega}_{S_3,4}+m_5=0$ gives $\gamma_3/\widehat{\omega}_{S_3,3}^{3/2}=\gamma_5$, since $\widehat{\omega}_{S_3,5}=1$. Thus
\[
\gamma_j^2=\gamma_5^2\widehat{\omega}_{S_3,j}^3\qquad(j\in S_3).
\]
Equation \eqref{eq:11.8} gives $\gamma_5\ne0$; together with the preceding identity, it implies $|\widehat{\omega}_{S_3,j}|=1$ for every $j\in S_3$. Two of the $\gamma_j^2$ coincide.
Choose distinct indices $j,k,l$ with $\{j,k,l\}=S_3$ and $\gamma_j^2=\gamma_k^2$. The ratio
$\widehat{\omega}_{S_3,k}/\widehat{\omega}_{S_3,j}$ has cube one and is not one, since the triangle coordinates are distinct. Hence
\[
\left(\frac{\widehat{\omega}_{S_3,k}}{\widehat{\omega}_{S_3,j}}\right)^2
+\frac{\widehat{\omega}_{S_3,k}}{\widehat{\omega}_{S_3,j}}+1=0,\qquad
\left|\frac{\widehat{\omega}_{S_3,l}}{\widehat{\omega}_{S_3,j}}\right|=1.
\]
The coordinate ratios are used only to exploit the degree-five homogeneity of the auxiliary polynomial $\mathcal E$ below. The normalized triangle coordinates are unchanged; in particular, the original condition $\widehat{\omega}_{S_3,5}=1$ remains in force. Permuting the order $j,k,l$
merely permutes the four signed Euler equations, all of which are retained below.

We use the signed Euler polynomial $\mathcal{E}$ defined in \eqref{eq:10.26}, with
$(w_3,w_4,w_5)$ replaced by $(\widehat{\omega}_{S_3,j},\widehat{\omega}_{S_3,k},\widehat{\omega}_{S_3,l})$ in its three argument slots.
The two signs are chosen according to the signed sheet of the limiting triangle, so that
$\mathcal{E}(\widehat{\omega}_{S_3,j},\widehat{\omega}_{S_3,k},\widehat{\omega}_{S_3,l})=0$.
Each of the four polynomials is homogeneous of degree five.

For each fixed sign choice, the one-argument notation below denotes the univariate specialization of this same polynomial:
\[
\mathcal{E}\left(\frac{\widehat{\omega}_{S_3,l}}{\widehat{\omega}_{S_3,j}}\right)
=\mathcal{E}\left(1,\frac{\widehat{\omega}_{S_3,k}}{\widehat{\omega}_{S_3,j}},
\frac{\widehat{\omega}_{S_3,l}}{\widehat{\omega}_{S_3,j}}\right).
\]
Here $\widehat{\omega}_{S_3,k}/\widehat{\omega}_{S_3,j}$ is fixed, and $\widehat{\omega}_{S_3,l}/\widehat{\omega}_{S_3,j}$ is treated as the polynomial variable,
then evaluated at its actual coordinate value. Homogeneity gives zero at that value.
Complex conjugation sends $\widehat{\omega}_{S_3,k}/\widehat{\omega}_{S_3,j}$ to its square.
Since the original three-variable polynomial has integer coefficients, the conjugate reciprocal of this specialization is
\[
\mathcal{E}^*\left(\frac{\widehat{\omega}_{S_3,l}}{\widehat{\omega}_{S_3,j}}\right)
=\left(\frac{\widehat{\omega}_{S_3,l}}{\widehat{\omega}_{S_3,j}}\right)^5
\mathcal{E}\left(1,\left(\frac{\widehat{\omega}_{S_3,k}}{\widehat{\omega}_{S_3,j}}\right)^2,
\frac{\widehat{\omega}_{S_3,j}}{\widehat{\omega}_{S_3,l}}\right).
\]
At the actual ratio of modulus one, $\mathcal{E}=0$ implies $\mathcal{E}^*=0$.
Exact Euclidean arithmetic in the polynomial ring
\[
\mathbb{Q}\left(\frac{\widehat{\omega}_{S_3,k}}{\widehat{\omega}_{S_3,j}}\right)
\left[\frac{\widehat{\omega}_{S_3,l}}{\widehat{\omega}_{S_3,j}}\right]
\]
gives
\begin{equation*}
d:=\gcd(\mathcal{E},\mathcal{E}^*)=
\begin{cases}
\displaystyle\frac{\widehat{\omega}_{S_3,l}}{\widehat{\omega}_{S_3,j}}
-\left(\frac{\widehat{\omega}_{S_3,k}}{\widehat{\omega}_{S_3,j}}\right)^2,&\text{for two plus signs},\\
1,&\text{for the other three choices}.
\end{cases}
\end{equation*}
Each $\mathcal{E}$ has degree five with leading coefficient $\pm1$, determined by the last sign in
\eqref{eq:10.26}; the conjugate reciprocal is formed with that exact degree, even if its own leading coefficient vanishes.
For each sign choice, the gcd calculation consists of an extended Euclidean identity
$P\mathcal{E}+V\mathcal{E}^*=d$, with $P,V$ in this polynomial ring, together with divisibility of both arguments by the displayed monic polynomial $d$.
The only possible shape therefore satisfies
$\widehat{\omega}_{S_3,j}\widehat{\omega}_{S_3,l}=\widehat{\omega}_{S_3,k}^2$ on the all-plus sheet. Together with the first ratio being a nontrivial cube root of one,
this makes the three cyclic differences have equal cubes, so the all-plus signed distances also have equal cubes. The
normalized triangle equations make their common cube $e_1(S_3)$. This is the Lagrange--Euler
intersection in Definition~\ref{def:lagrange_euler}, contrary to the hypothesis.
\end{proof}

\subsection{The partial Hessian}\label{sec:11.4}
\label{subsec:11.4}
At a collinear limiting triangle, the normalization $\widehat{\omega}_{S_3,5}=1$ gives
\[
\widehat{\zeta}_{S_3,j}=\widehat{\zeta}_{S_3,5}\widehat{\omega}_{S_3,j},\qquad j\in S_3,\qquad \widehat{\zeta}_{S_3,5}\ne0.
\]
Indeed a nonisotropic line gives such a proportionality, and an isotropic line would make every
squared distance zero, contrary to the diagram estimates. If the two centered coordinate vectors are
independent instead, the triangle is Lagrange: the normalized force operator acts as the identity on
the whole centered plane and all three inverse-cube weights equal $1/e_1(S_3)$. The collinear and independent-vector alternatives are exhaustive; the collinear alternative may also be Lagrange, as specified in Definition~\ref{def:lagrange_euler}.

On a collinear triangle define the signed force matrix $L$, with rows and columns ordered as $3,4,5$, by
\begin{equation}
L_{jk}=-m_k\delta_{jk}\quad(j\ne k),\qquad
L_{jj}=\sum\limits_{\substack{k\in S_3\\k\ne j}}m_k\delta_{jk},\qquad j,k\in S_3.
\label{eq:11.10}
\end{equation}
It preserves the two-dimensional centered plane and is self-adjoint for the complex bilinear mass form.
The normalized triangle coordinates satisfy
\[
\sum\limits_{\substack{k\in S_3\\k\ne j}}m_k(\widehat{\omega}_{S_3,j}-\widehat{\omega}_{S_3,k})\delta_{jk}
=\widehat{\omega}_{S_3,j},\qquad j\in S_3.
\]
Thus one of its centered eigenvalues is one. Its other centered eigenvalue, counted algebraically, is
\begin{equation*}
\tau=\sum\limits_{\substack{j,k\in S_3\\j<k}}(m_j+m_k)r_{jk}^{-3}-1.
\end{equation*}
This definition assumes neither diagonalizability nor a nonzero moment of inertia.

Eliminate $\widehat{\zeta}_{S_3,5}$ using the centering relation in \eqref{eq:11.2}.
Differentiate $G_{(\mathfrak{c},0)}$ in \eqref{eq:11.4} with respect to the two independent coordinates
$\widehat{\zeta}_{S_3,3},\widehat{\zeta}_{S_3,4}$, keeping $\widehat{\omega}_{S_3,j}$ fixed for every $j\in S_3$,
and then evaluate at the selected collinear limit. The notation $\left.\cdot\right|_0$ means substitution of $\widehat{\zeta}_{S_3,T_2}=\widehat{\zeta}_{S_3,T_2}^{0}$ and $\widehat{\omega}_{S_3,3}=\widehat{\omega}_{S_3,3}^0$ after differentiation; the dependent coordinates, distances, $L$, and $\tau$ are evaluated at that same point. Define
\[
H^{0}_{\widehat{\zeta}_{S_3,T_2}\widehat{\zeta}_{S_3,T_2}}
:=\left.D^2_{(\widehat{\zeta}_{S_3,3},\widehat{\zeta}_{S_3,4})}G_{(\mathfrak{c},0)}\right|_0.
\]
Using $m_jL_{jk}=m_kL_{kj}$ gives
{\small\setlength{\arraycolsep}{4pt}\begin{equation}
\begin{aligned}
&H^{0}_{\widehat{\zeta}_{S_3,T_2}\widehat{\zeta}_{S_3,T_2}}\\
&\quad=\left.\frac{3}{4\widehat{\zeta}_{S_3,5}}
\begin{pmatrix}
m_3L_{33}+\frac{m_3^2}{m_5}(L_{55}-2L_{35})
&
m_3L_{34}+\frac{m_3m_4}{m_5}(L_{55}-L_{35}-L_{45})
\\[4pt]
m_3L_{34}+\frac{m_3m_4}{m_5}(L_{55}-L_{35}-L_{45})
&
m_4L_{44}+\frac{m_4^2}{m_5}(L_{55}-2L_{45})
\end{pmatrix}\right|_0,\\[6pt]
&\det H^{0}_{\widehat{\zeta}_{S_3,T_2}\widehat{\zeta}_{S_3,T_2}}
=\left.\frac{9m_3m_4e_1(S_3)}{16m_5\widehat{\zeta}_{S_3,5}^2}\tau\right|_0.
\end{aligned}
\label{eq:11.12}
\end{equation}}
For each edge, the coefficient of the squared coordinate difference in the second variation is
$\frac34m_jm_k(\widehat{\omega}_{S_3,j}-\widehat{\omega}_{S_3,k})^2r_{jk}^{-5}$.
Using $r_{jk}^2=\widehat{\zeta}_{S_3,5}(\widehat{\omega}_{S_3,j}-\widehat{\omega}_{S_3,k})^2$
and eliminating the fifth coordinate gives the displayed matrix.
The determinant formula follows from the centered eigenvalues $1,\tau$ of $L$
and the determinant $m_3m_4e_1(S_3)/m_5$ of the mass form in the two independent coordinates.
Thus $H^{0}_{\widehat{\zeta}_{S_3,T_2}\widehat{\zeta}_{S_3,T_2}}$ is invertible whenever $\tau$ is nonzero at the selected limit. This does not require the full Euler shape Hessian
to be invertible; in particular no exceptional value such as $\tau=-1/2$ or $\tau=1$ is discarded, and there
is no division by the triangle's moment of inertia.

\subsection{Exclusion when the partial Hessian is invertible}\label{sec:11.5}
\label{subsec:11.5}
\begin{lemma}\label{lem:11.3}
A positive-mass complex normalized sequence cannot have a non-Lagrange Euler
diagram-2 limit with $\tau\ne0$.
\end{lemma}
\begin{proof}
Let $((r_{12}^0)^2,\widehat{\zeta}_{S_3,T_2}^{0},\widehat{\omega}_{S_3,3}^0,\nu^0)$ be the bounded-coordinate limit in the fixed signed chart of Section~\ref{sec:11.1},
where $r_{12}^0=\lim\limits_{\epsilon\to0}r_{12}\ne0$ and
$\widehat{\zeta}_{S_3,T_2}^{0}=\lim\limits_{\epsilon\to0}\widehat{\zeta}_{S_3,T_2}$ and $\nu^0=\lim\limits_{\epsilon\to0}\nu$. The limiting internal
equations make it a critical point of $G_{(\mathfrak{c},0)}$. Lemma~\ref{lem:11.2} gives $G'_{(\mathfrak{c},1)}(\widehat{\omega}_{S_3,3}^0)\ne0$. From \eqref{eq:11.6}, uniformly near
this limit,
\begin{equation}
\begin{aligned}
G_{\mathfrak{c}}&=G_{(\mathfrak{c},1)}(\widehat{\omega}_{S_3,3})+O(\mathfrak{c}^2),& G_{r_{12}^2}&=(G_{(\mathfrak{c},0)})_{r_{12}^2}+O(\mathfrak{c}^3),\\
G_{\widehat{\zeta}_{S_3,T_2}}&=(G_{(\mathfrak{c},0)})_{\widehat{\zeta}_{S_3,T_2}}+O(\mathfrak{c}^3),& G_{\widehat{\omega}_{S_3,3}}&=(G_{(\mathfrak{c},0)})_{\widehat{\omega}_{S_3,3}}+\mathfrak{c}G'_{(\mathfrak{c},1)}(\widehat{\omega}_{S_3,3})+O(\mathfrak{c}^3).
\end{aligned}
\label{eq:11.13}
\end{equation}
The first equation and the local inverse-function estimate for $G_{(\mathfrak{c},1)}$ give $\widehat{\omega}_{S_3,3}-\widehat{\omega}_{S_3,3}^0=O(\mathfrak{c}^2)$. The contribution of the binary to $G_{(\mathfrak{c},0)}$ has, with respect to $r_{12}^2$, nonzero second derivative $3m_1m_2/(4r_{12}^5)$, independently of its signed phase, so its equation
gives $r_{12}^2-(r_{12}^0)^2=O(\mathfrak{c}^3)$.

Use only the two equations $G_{\widehat{\zeta}_{S_3,T_2}}=0$ to solve for $\widehat{\zeta}_{S_3,T_2}$ as an analytic function of $(\widehat{\omega}_{S_3,3},r_{12}^2,\mathfrak{c},\nu)$ near the
limit. At that limit, their Jacobian with respect to $\widehat{\zeta}_{S_3,T_2}$ is $H^{0}_{\widehat{\zeta}_{S_3,T_2}\widehat{\zeta}_{S_3,T_2}}$ in \eqref{eq:11.12}, which is invertible. Since the leading contribution from the triangle is independent of $r_{12}^2,\nu$ and the perturbation in these equations begins at order $\mathfrak{c}^3$, their
solution satisfies
\[
\widehat{\zeta}_{S_3,T_2}-\widehat{\zeta}_{S_3,T_2}^{0}=O(\widehat{\omega}_{S_3,3}-\widehat{\omega}_{S_3,3}^0)+O(\mathfrak{c}^3)=O(\mathfrak{c}^2).
\]
These estimates are uniform for $\nu$ near $\nu^0$; solving a center equation or assuming an analytic
approaching arc is unnecessary. Because $(G_{(\mathfrak{c},0)})_{\widehat{\omega}_{S_3,3}}=0$ at the leading point, the final equation in \eqref{eq:11.13}
now reads
\[
0=G_{\widehat{\omega}_{S_3,3}}=\mathfrak{c}G'_{(\mathfrak{c},1)}(\widehat{\omega}_{S_3,3}^0)+O(\mathfrak{c}^2).
\]
This is impossible for sufficiently small nonzero complex $\mathfrak{c}$. The contradiction applies directly to every
convergent subsequence of bounded chart variables, including one initially described by fractional
Puiseux powers.
\end{proof}

\subsection{Exclusion when the partial Hessian is singular}\label{sec:11.6}
\label{subsec:11.6}
\begin{lemma}\label{lem:11.4}
Let a limiting signed Euler triangle in diagram 2 have positive masses, and let $\tau$ be
the second eigenvalue of its normalized force operator on the centered two-dimensional space, the
first eigenvalue being one. If $\tau=0$, then the triangle masses are proportional to $(1,1,5)$, and no
positive binary masses can satisfy the limiting crossed-term constraint $G_{(\mathfrak c,1)}=0$.
\end{lemma}
\begin{proof}
Use centered collinear coordinates $\widehat{\omega}_{S_3,j}$, $j\in S_3$, with every coordinate and every pairwise difference
nonzero, as required in diagram 2. Use the force matrix $L$ defined in \eqref{eq:11.10}.
The force operator is self-adjoint for the mass bilinear form and satisfies $L1=0$.
Here $1=(1,1,1)^T$ is the constant vector indexed by $S_3$; $L1=0$ follows because each row of $L$ sums to zero.
This vector lies outside the centered plane, since $e_1(S_3)>0$.
The normalized coordinate equations make the column vector with components $\widehat{\omega}_{S_3,j}$, $j\in S_3$, an eigenvector with eigenvalue one.
The centered mass form is nondegenerate because the positive total mass is nonzero. If $\tau=0$, its two
eigenvalues one and zero are distinct, and their eigenlines are orthogonal. Hence this coordinate eigenvector satisfies
$\sum\limits_{l\in S_3}m_l\widehat{\omega}_{S_3,l}^2\ne0$: otherwise its orthogonal line in the centered plane would
coincide with itself and could not contain the independent zero eigenvector. Consequently $L$ is the
mass-orthogonal rank-one projection onto the line spanned by this coordinate vector, and its off-diagonal entries give
\begin{equation}
\delta_{jk}=-\frac{\widehat{\omega}_{S_3,j}\widehat{\omega}_{S_3,k}}{\displaystyle\sum\limits_{l\in S_3}m_l\widehat{\omega}_{S_3,l}^2},\qquad j,k\in S_3,\quad j\ne k.
\label{eq:11.14}
\end{equation}
Take uncentered line coordinates $(0,1,z)$ for bodies $3,4,5$.
In this subsection, $\sigma_{35},\sigma_{45}\in\{1,-1\}$ denote the signs of the $35$ and $45$ distance branches relative to the $34$ branch:
$r_{35}=\sigma_{35}zr_{34}$ and $r_{45}=\sigma_{45}(z-1)r_{34}$. Thus
\[
\frac{\delta_{35}}{\delta_{34}}=\frac{\sigma_{35}}{z^3},\qquad
\frac{\delta_{45}}{\delta_{34}}=\frac{\sigma_{45}}{(z-1)^3}.
\]
Equation \eqref{eq:11.14} implies
\[
(\widehat{\omega}_{S_3,3},\widehat{\omega}_{S_3,4},\widehat{\omega}_{S_3,5})\text{ is proportional to }(\sigma_{45}(z-1)^3,\sigma_{35} z^3,1).
\]
The affine ratio $\widehat{\omega}_{S_3,5}-\widehat{\omega}_{S_3,3}=z(\widehat{\omega}_{S_3,4}-\widehat{\omega}_{S_3,3})$ therefore forces
\begin{equation}
\sigma_{35} z^4-\sigma_{45}(z-1)^4=1.
\label{eq:11.15}
\end{equation}
For the four choices of $(\sigma_{35},\sigma_{45})$, the left side of \eqref{eq:11.15} minus one factors as
\[
\begin{aligned}
(1,1) &: \quad 2(z-1)(2z^2-z+1),\\
(1,-1) &: \quad 2z(z-1)(z^2-z+2),\\
(-1,1) &: \quad -2(z^2-z+1)^2,\\
(-1,-1) &: \quad -2z(2z^2-3z+2).
\end{aligned}
\]
The factors $z,z-1$ represent collisions. The remaining quadratics all have negative discriminant.

The center-of-mass condition is $m_3\sigma_{45}(z-1)^3+m_4\sigma_{35} z^3+m_5=0$.
Reducing it modulo the four quadratics gives, respectively,
\[
\begin{aligned}
&(m_3-m_4+4m_5)+(5m_3-m_4)z=0,\\
&(-3m_3-2m_4+m_5)+(m_3-m_4)z=0,\\
&e_1(S_3)=0,\\
&(-2m_3+6m_4+4m_5)+(m_3-5m_4)z=0.
\end{aligned}
\]
Because $z$ is nonreal and the masses are real, the first, second and fourth cases give ratios $(1,5,1)$,
$(1,1,5)$, and $(5,1,1)$. The third is impossible for positive masses. This exhausts all signed zero-eigenvalue triangles.

After relabeling within $S_3$ and a common mass normalization, normalize the centered coordinates by $\widehat{\omega}_{S_3,5}=1$ to obtain
\begin{equation}
(m_3,m_4,m_5)=(1,1,5),\qquad
\widehat{\omega}_{S_3,4}=-\widehat{\omega}_{S_3,3}-5,\qquad \widehat{\omega}_{S_3,3}^2+5\widehat{\omega}_{S_3,3}+8=0.
\label{eq:11.16}
\end{equation}
Both conjugate roots are retained. In the case $(\sigma_{35},\sigma_{45})=(1,-1)$ above, the centered vector is
$(z-3,-z-2,1)$ and $z^2-z+2=0$; putting $\widehat{\omega}_{S_3,3}=z-3$ gives \eqref{eq:11.16}. The other two cases give the
same description after permuting the equal-mass endpoints and the mass-five body.

It remains to show that the limiting crossed-term constraint $G_{(\mathfrak c,1)}=0$ fails on every signed sheet. Take $\widehat{\omega}_{S_3,3}=
(-5+\mathrm{i}\sqrt7)/2$, so $\widehat{\omega}_{S_3,4}=\overline{\widehat{\omega}_{S_3,3}}$ and $\widehat{\omega}_{S_3,5}=1$. Choose
\[
\begin{aligned}
\frac1{\sqrt{\widehat{\omega}_{S_3,3}}}
&=\sqrt{\frac{4\sqrt2-5}{32}}-\mathrm{i}\sqrt{\frac{4\sqrt2+5}{32}},\\
\frac1{\sqrt{\widehat{\omega}_{S_3,4}}}
&=\sqrt{\frac{4\sqrt2-5}{32}}+\mathrm{i}\sqrt{\frac{4\sqrt2+5}{32}}.
\end{aligned}
\]
The square roots on the right are positive real, and we take $\sqrt{\widehat{\omega}_{S_3,5}}=1$.
For the mass square roots in \eqref{eq:11.7}, choose the positive real root of $m_2/(m_1+m_2)$ and set
$\sqrt{-m_1/(m_1+m_2)}=\mathrm{i}\sqrt{m_1/(m_1+m_2)}$, with the root on the right positive.
Other choices of these square roots are absorbed into the signs $\sigma_{ij}$.
By \eqref{eq:11.7}--\eqref{eq:11.8}, the coefficient $G_{(\mathfrak c,1)}$ is
\[
G_{(\mathfrak{c},1)}=\gamma_{\mathrm{R}}\sum\limits_{j\in S_3}\frac{m_j\sigma_{1j}}{\sqrt{\widehat{\omega}_{S_3,j}}}
-\mathrm{i}\gamma_{\mathrm{I}}\sum\limits_{j\in S_3}\frac{m_j\sigma_{2j}}{\sqrt{\widehat{\omega}_{S_3,j}}},\qquad
\gamma_{\mathrm{R}},\gamma_{\mathrm{I}}>0.
\]
Each of the two sums is either a nonzero real number $\pm5\pm\sqrt{(4\sqrt2-5)/8}$ or a nonreal number
$\pm5\pm \mathrm{i}\sqrt{(4\sqrt2+5)/8}$, since
\[
0<\frac{4\sqrt2-5}{8}<\frac{4\sqrt2+5}{8}<2.
\]
For every choice of the six cross signs $\sigma_{ij}$, the product of the first sum and the complex conjugate of the second has nonzero real part.
Indeed, when both sums are nonreal, its real part belongs to $\{\pm25\pm(4\sqrt2+5)/8\}$, whose elements are all nonzero;
when exactly one sum is real, the real part is $\pm5$ times that nonzero real sum;
and when both sums are real, it is their nonzero product. But $G_{(\mathfrak{c},1)}=0$ would give
\[
\left(\sum\limits_{j\in S_3}\frac{m_j\sigma_{1j}}{\sqrt{\widehat{\omega}_{S_3,j}}}\right)
\overline{\left(\sum\limits_{k\in S_3}\frac{m_k\sigma_{2k}}{\sqrt{\widehat{\omega}_{S_3,k}}}\right)}
=\mathrm{i}\frac{\gamma_{\mathrm{I}}}{\gamma_{\mathrm{R}}}
\left|\sum\limits_{k\in S_3}\frac{m_k\sigma_{2k}}{\sqrt{\widehat{\omega}_{S_3,k}}}\right|^2.
\]
The right side has zero real part, a contradiction. This excludes every positive binary ratio and every crossed sign choice.
The conjugate root exchanges the first two coordinates, and relabeling preserves the equation $G_{(\mathfrak c,1)}=0$, so all canonical cases
are excluded.
\end{proof}

\begin{proof}[Proof of Theorem~\ref{thm:11.1}]
By Definition~\ref{def:lagrange_euler} and the normalized triangle equations, the internal limiting triangle is either Lagrange (including the Lagrange--Euler intersection) or a non-Lagrange Euler configuration.
In the Lagrange case, its three distance cubes tend to $e_1(S_3)$; the binary distance cube tends to $m_1+m_2$, and all
crossed distances diverge. Lemma~\ref{lem:3.2} excludes this separated-block limit on a component containing
a physical point, for every choice of distance phases.

For a non-Lagrange Euler triangle, Lemma~\ref{lem:11.3} excludes the case $\tau\ne0$. The remaining case
$\tau=0$ is excluded by Lemma~\ref{lem:11.4}: the only possible positive triangle masses are proportional to
$(1,1,5)$, and their limiting crossed-term constraint $G_{(\mathfrak c,1)}=0$ is inconsistent for every positive binary mass ratio. All
alternatives are exhausted, proving the theorem.
\end{proof}

\section{The exclusion of Diagram 5}
\label{sec:5}
\label{sec:6}
\label{sec:7}

The global sets $S_2,T_2,S_3,S_4,T_4,S_5$ in \eqref{eq:2.index_sets} index the exterior pair, the two triangle endpoints, the isolated $zw$-triangle, the subset excluding body 5, the subset excluding body 2, and the whole system, respectively. A cross pair has one index in $S_2$ and the other in $S_3$.

\begin{theorem}\label{thm:7.1}
For arbitrary positive five-body masses, a normalized complex solution component containing a physical point cannot admit a singular sequence approaching AK diagram 5 along which some 4-product tends to infinity.
\end{theorem}

Suppose, for a contradiction, that such a singular sequence exists. After relabeling the internal vertices so that the smallest internal mass is at vertex 5, the masses of the isolated $zw$-triangle satisfy the necessary relation given in \cite[Section~5.3]{AK}:
\[
\frac1{\sqrt{m_5}}=\frac1{\sqrt{m_3}}+\frac1{\sqrt{m_4}}.
\]
Masses for which no internal relabeling satisfies this relation already exclude Diagram 5. We first establish the common identities and coordinate reduction in Section~\ref{sec:6.1} and the Schur complement estimates in Section~\ref{sec:5.2}, then treat $m_3=m_4$ in Section~\ref{sec:diagram5:equal} and $m_3\ne m_4$ in Section~\ref{sec:diagram5:unequal}. Each branch uses its own common positive mass scaling and the corresponding coordinate scaling that keeps the multiplier equal to one. Such a scaling multiplies $U$, $I$, and the physical lower bound \eqref{eq:3.1} by the same positive factor, and does not change any order comparison.

\subsection{Common identities and coordinate reduction}\label{sec:6.1}
\label{subsec:6.1}
\label{sec:7.2}\label{subsec:7.2}

We use the cluster potential $U_S$ and the contributions $f_i^S$ to the right-hand sides of the central configuration equations, as defined in Section~\ref{sec:2}.

\begin{lemma}[Identities for the cluster potential]\label{lem:5.2}
For every nonempty $S\subseteq S_5$ and any $k\in S$,
\begin{equation*}
U_S=\sum\limits_{i\in S}m_i(z_i-z_k)f_{i,w}^S
=\sum\limits_{i\in S}m_i(w_i-w_k)f_{i,z}^S.
\end{equation*}
The identities are translation invariant: positions may be measured relative to any body $k\in S$, and the term with $i=k$ vanishes identically.
\end{lemma}
\begin{proof}
In the first sum, the two summands for an unordered pair $i,j$ combine to $m_im_jz_{ij}W_{ij}=m_im_j/r_{ij}$. The terms containing $z_k$ cancel because $\sum\limits_{i\in S}m_if_i^S=0$. The second identity follows by exchanging the coordinates. Each sum equals $U_S$, so neither identity has a factor $1/2$.
\end{proof}

The following reduction uses only positivity of the fixed masses and the diagram-5 estimates; no common mass normalization is imposed here. Write $M=e_1(S_5)$. The physical component gives $U=I\ge e_2(S_5)/M^{1/3}>0$, and
\begin{equation*}
\begin{aligned}
e_2(S_5)&=m_1m_2+e_1(S_3)(m_1+m_2)+e_2(S_3),\\
e_2(S_3)&=m_3m_4+m_3m_5+m_4m_5>0.
\end{aligned}
\end{equation*}
By the AK estimates in Section~\ref{sec:2.2}, all internal differences $z_{jk},w_{jk},r_{jk}$, for distinct $j,k\in S_3$, have order $\epsilon$, every cross distance is much larger than $\epsilon$, and all absolute coordinates are $o(\epsilon^{-2})$. In particular,
\[
Z_{jk}=\frac{z_{jk}}{r_{jk}^3}\approx\epsilon^{-2},\qquad
W_{jk}=\frac{w_{jk}}{r_{jk}^3}\approx\epsilon^{-2},
\qquad j,k\in S_3,\quad j\ne k.
\]

With the global index sets, $C_{S_2}$ and $C_{S_3}$ are the centers of mass of the exterior pair and the internal triangle, respectively. Total centering gives
\[
(m_1+m_2)C_{S_2}+e_1(S_3) C_{S_3}=0.
\]
For $i\in S_2$, use the relative coordinates $\zeta_{S_3,i},\omega_{S_3,i}$ defined in Section~\ref{sec:2}; in particular, $\zeta_{S_3,1}-\zeta_{S_3,2}=z_{12}$. The exact identities for the center of mass and the moment of inertia give
\begin{equation}
C_{S_3,z}=-\frac{m_1\zeta_{S_3,1}+m_2\zeta_{S_3,2}}{M},\qquad C_{S_3,w}=-\frac{m_1\omega_{S_3,1}+m_2\omega_{S_3,2}}{M},
\label{eq:6.2}
\end{equation}
\begin{equation}
\begin{aligned}
I={}&m_1\zeta_{S_3,1}\omega_{S_3,1}+m_2\zeta_{S_3,2}\omega_{S_3,2}\\
&-\frac{(m_1\zeta_{S_3,1}+m_2\zeta_{S_3,2})(m_1\omega_{S_3,1}+m_2\omega_{S_3,2})}{M}+I_{S_3},
\qquad I_{S_3}=O(\epsilon^2).
\end{aligned}
\label{eq:6.3}
\end{equation}
Initially let
\[
\Theta_{\zeta}=\max(\epsilon,|\zeta_{S_3,1}|,|\zeta_{S_3,2}|),\qquad
\Theta_{\omega}=\max(\epsilon,|\omega_{S_3,1}|,|\omega_{S_3,2}|).
\]
We pass to subsequences as needed so that the distances and the normalized ratios under consideration have limits as $\epsilon\to0$. A cross distance $r_{ij}$, $i\in S_2$, $j\in S_3$, is called an \emph{attachment} if its limit is finite, possibly zero.
Define the minimum modulus of these six cross distances by
\[
\rho=\min_{\substack{i\in S_2\\j\in S_3}}|r_{ij}|
=\min\{|r_{13}|,|r_{14}|,|r_{15}|,|r_{23}|,|r_{24}|,|r_{25}|\}.
\]

If some 4-product tends to infinity, at least one exterior distance must diverge, since every such product contains an internal edge. If all six cross distances diverge, then $\rho\to\infty$ and $\delta_{ij}\to0$ for $i\in S_2$, $j\in S_3$. The mass-weighted cluster identities in Section~\ref{sec:2.1} give
\[
C_{S_3,z}=o\!\left(\max\limits_{i\in S_5}|z_i|\right),\qquad
C_{S_3,w}=o\!\left(\max\limits_{i\in S_5}|w_i|\right).
\]
After exchanging coordinates if necessary, choose $\Theta_w=\max\limits_{i\in S_5}|w_i|\succ\epsilon$. Since the internal $w$ differences are $O(\epsilon)$, all triangle coordinates satisfy $w_j=o(\Theta_w)$ for $j\in S_3$. Total centering then gives $w_i\approx\Theta_w$ for $i\in S_2$, with $w_1/w_2\to-m_2/m_1$. Thus $w_{12}\approx\Theta_w$ and $w_{ij}\approx\Theta_w$ for $i\in S_2$, $j\in S_3$. The pair-component identities now give
\[
Z_{ij}=\frac{1}{r_{ij}w_{ij}}=O\!\left(\frac1{\rho\Theta_w}\right),\qquad
W_{ij}=w_{ij}\delta_{ij}=o(\Theta_w),\qquad i\in S_2,\quad j\in S_3.
\]
Subtracting the $w$ equations for bodies 1 and 2 yields
\[
\begin{aligned}
w_{12}&=(m_1+m_2)W_{12}+\sum\limits_{j\in S_3}m_j(W_{1j}-W_{2j})\\
&=(m_1+m_2)W_{12}+o(\Theta_w).
\end{aligned}
\]
Since $w_{12}\approx\Theta_w$, it follows that $r_{12}^3=w_{12}/W_{12}\to m_1+m_2$. Applying the mass-weighted cluster identity to the $Z_{ij}$ estimates gives $C_{S_3,z}=O((\rho\Theta_w)^{-1})$. Together with $C_{S_3,w}=o(\Theta_w)$, this implies $C_{S_3,z}C_{S_3,w}=o(\rho^{-1})\to0$. The cluster decomposition of the moment of inertia gives $I=m_1m_2r_{12}^2/(m_1+m_2)+o(1)$. Since $I>0$ is constant, taking moduli as at the end of Proposition~\ref{prop:4.1} yields $I=m_1m_2/(m_1+m_2)^{1/3}$ on every distance branch. This is below the physical lower bound because concavity and $m_1m_2\le (m_1+m_2)^2/4$ give
\begin{equation}
\begin{aligned}
m_1m_2(M/(m_1+m_2))^{1/3}
&\le m_1m_2+\frac{e_1(S_3)}3\frac{m_1m_2}{(m_1+m_2)}\\
&\le m_1m_2+e_1(S_3)(m_1+m_2)/12<e_2(S_5).
\end{aligned}
\label{eq:6.4}
\end{equation}
Consequently at least one cross distance is bounded.

Exterior divergence gives $\Theta_{\zeta} \Theta_{\omega}\to\infty$. If $\Theta_{\zeta},\Theta_{\omega}\succ\epsilon$, pass to limiting vectors
\[
\widehat{\zeta}_{S_3,S_2}^{0}=\lim\limits_{\epsilon\to0}\frac{(\zeta_{S_3,1},\zeta_{S_3,2})^T}{\Theta_{\zeta}},\qquad
\widehat{\omega}_{S_3,S_2}^{0}=\lim\limits_{\epsilon\to0}\frac{(\omega_{S_3,1},\omega_{S_3,2})^T}{\Theta_{\omega}},
\]
each with maximum component modulus one. Denote their components by $\widehat{\zeta}_{S_3,i}^0$ and $\widehat{\omega}_{S_3,i}^0$, respectively, to distinguish them from the unscaled relative coordinates $\zeta_{S_3,i},\omega_{S_3,i}$. A bounded attachment at body 1 gives $\widehat{\zeta}_{S_3,1}^0\widehat{\omega}_{S_3,1}^0=0$. Exchange coordinates if needed so $\widehat{\zeta}_{S_3,1}^0=0$; then $\widehat{\zeta}_{S_3,2}^0\ne0$. Boundedness of the moment of inertia now gives $m_1\widehat{\omega}_{S_3,1}^0=(m_1+e_1(S_3))\widehat{\omega}_{S_3,2}^0$, so $\widehat{\omega}_{S_3,2}^0\ne0$ and $\widehat{\omega}_{S_3,1}^0-\widehat{\omega}_{S_3,2}^0=(e_1(S_3)/m_1)\widehat{\omega}_{S_3,2}^0\ne0$. All distances incident to body 2 diverge. For every $j\in S_5\setminus\{2\}$, we have $z_{2j}/\Theta_{\zeta}=O(1)$ and $\delta_{2j}\to0$, hence
\[
\frac{Z_{2j}}{\Theta_{\zeta}}
=\frac{z_{2j}}{\Theta_{\zeta}}\delta_{2j}\longrightarrow0.
\]
The normalized $z$ equation for body 2 therefore gives
\[
\frac{z_2}{\Theta_{\zeta}}
=\sum\limits_{j\in S_5\setminus\{2\}}m_j\frac{Z_{2j}}{\Theta_{\zeta}}
\longrightarrow0.
\]
On the other hand, \eqref{eq:6.2} gives
\[
\frac{z_2}{\Theta_{\zeta}}
=\frac{m_1+e_1(S_3)}{M}\frac{\zeta_{S_3,2}}{\Theta_{\zeta}}
-\frac{m_1}{M}\frac{\zeta_{S_3,1}}{\Theta_{\zeta}}
\longrightarrow\frac{m_1+e_1(S_3)}{M}\widehat{\zeta}_{S_3,2}^0\ne0,
\]
a contradiction. Hence, after exchanging coordinates,
\begin{equation}
w_i=O(\epsilon)\quad(i\in S_5),\qquad \epsilon \Theta_{\zeta}\to\infty,\qquad \Theta_{\zeta}=o(\epsilon^{-2}),\qquad \zeta_{S_3,1},\zeta_{S_3,2}\succ\epsilon,
\label{eq:6.5}
\end{equation}
where $\Theta_{\zeta}=\max(|\zeta_{S_3,1}|,|\zeta_{S_3,2}|)$. The bound $\Theta_{\zeta}=o(\epsilon^{-2})$ follows because every vertex is uncircled. If $\zeta_{S_3,i}=O(\epsilon)$ for some $i\in S_2$, every cross squared distance from that body would be $O(\epsilon^2)$, contradicting the absence of cross $zw$ edges. Thus $\zeta_{S_3,i}\succ\epsilon$. Every exterior body with $\zeta_{S_3,i}\approx\Theta_{\zeta}$ has at most one attachment: a bounded distance has $w$-separation $O(\Theta_{\zeta}^{-1})=o(\epsilon)$, whereas distinct internal $w$ levels are separated at scale $\epsilon$.

We record a common estimate for the branches in which all six cross inverse cubes are bounded. For $i\in S_2$ and $j\in S_3$, \eqref{eq:6.5} gives $w_{ji}=O(\epsilon)$, hence $W_{ji}=-W_{ij}=w_{ji}\delta_{ij}=O(\epsilon)$. The individual $w$ equations therefore give
\[
f_{j,w}^{S_3}=w_j-\sum\limits_{i\in S_2}m_iW_{ji}=O(\epsilon),\qquad j\in S_3.
\]
Using the first identity in Lemma~\ref{lem:5.2} with $k=5$ and the AK bounds $z_{j5}=O(\epsilon)$ for $j\in T_2$, we obtain
\begin{equation}
U_{S_3}=\sum\limits_{j\in T_2}m_jz_{j5}f_{j,w}^{S_3}=O(\epsilon^2)\longrightarrow0.
\label{eq:6.6}
\end{equation}
This estimate does not require a bound on $\delta_{12}$.

After a further subsequence and, in the last case, an exchange of the two exterior labels, the relative $z$-coordinates $\zeta_{S_3,1},\zeta_{S_3,2}$ and their $z$-separation $z_{12}$ satisfy exactly one of the following cases:
\begin{enumerate}[label=\textup{\Alph*.},ref=\Alph*]
\item\label{case:6.A} $\zeta_{S_3,1}\approx\zeta_{S_3,2}\approx\zeta_{S_3,1}-\zeta_{S_3,2}\approx\Theta_{\zeta}$.
For $m_3=m_4$, Section~\ref{subsec:6.2} derives the remaining three-point limit, which Section~\ref{subsec:6.3} excludes by its potential. For $m_3\ne m_4$, Section~\ref{subsec:7.3} excludes this case directly.
\item\label{case:6.B} $\zeta_{S_3,1}\approx\zeta_{S_3,2}\approx\Theta_{\zeta}$ and $\zeta_{S_3,1}-\zeta_{S_3,2}=o(\Theta_{\zeta})$.
This case is excluded in Section~\ref{subsec:6.4} when $m_3=m_4$, and in Section~\ref{subsec:7.3} when $m_3\ne m_4$.
\item\label{case:6.C} $\epsilon\prec\zeta_{S_3,1}\prec|\zeta_{S_3,2}|=\Theta_{\zeta}$.
This case is excluded in Section~\ref{subsec:6.5} when $m_3=m_4$, and in Sections~\ref{subsec:7.7}--\ref{subsec:7.6} when $m_3\ne m_4$.
\end{enumerate}

\subsection{Schur complement estimates}\label{sec:5.2}

The following definitions and estimates apply to both mass cases. No special normalization of the positive masses is imposed in this subsection. Define the symmetric matrix $\Omega$ associated with the centered central configuration equations by
\begin{equation*}
\Omega_{ij}=-m_im_j(\delta_{ij}-1/M)\quad(i\ne j),\qquad
\Omega_{ii}=\sum\limits_{j\ne i}m_im_j(\delta_{ij}-1/M).
\end{equation*}
Let $\mathcal{L}$ be the exact internal triangle Laplacian, with off-diagonal entries $\mathcal{L}_{ij}=-m_im_j\delta_{ij}$ for distinct $i,j\in S_3$, zero entries involving the exterior bodies, and row sums zero. Define $\mathcal{R}=\Omega-\mathcal{L}$. Then
\[
\Omega=\mathcal{L}+\mathcal{R},\qquad
\mathcal{L}_{55}=m_5(m_3\delta_{35}+m_4\delta_{45}),\qquad
\Omega_{55}=\mathcal{L}_{55}+\mathcal{R}_{55}.
\]
The row sums of $\Omega$ vanish. The force and center equations give $\Omega z=\Omega w=0$. Thus, whenever $1,z,w$ are linearly independent over $\mathbb{C}$,
\begin{equation}\label{eq:5.11}
\operatorname{rank}\Omega\le2.
\end{equation}
The common reduction gives $\Theta_{\zeta}\to\infty$. Whenever a bounded attachment has $z$-separation comparable to $\Theta_{\zeta}$, it establishes the required independence: the ratio of its $z$- and $w$-separations has modulus at least a constant times $\Theta_{\zeta}^2$, whereas an internal edge has a ratio bounded above and away from zero. A relation $z=aw+b1$ cannot hold for both. We use this criterion in the branches where such an attachment has been established.

The Schur complement is the coefficient matrix obtained by eliminating variables using an invertible diagonal block. This operation requires neither symmetry nor invertibility of the full matrix. Here, whenever $\Omega_{55}\ne0$, the fifth row of $\Omega z=0$ gives
\[
z_5=-\frac{1}{\Omega_{55}}\sum\limits_{l\in S_4}\Omega_{5l}z_l.
\]
Substituting into the other four rows gives $\sum\limits_{l\in S_4}\mathcal{S}_{kl}z_l=0$ for $k\in S_4$, where the $4\times4$ matrix $\mathcal{S}$ is defined by
\[
\mathcal{S}_{kl}=\Omega_{kl}-\frac{\Omega_{k5}\Omega_{5l}}{\Omega_{55}},
\qquad k,l\in S_4.
\]
This is the Schur complement of $\Omega_{55}$ in $\Omega$. The same elimination applies to $\Omega w=0$. Invertible row and column operations reduce $\Omega$ to the block diagonal matrix with blocks $\mathcal{S}$ and $\Omega_{55}$, so
\[
\operatorname{rank}\Omega=1+\operatorname{rank}\mathcal{S}.
\]
Consequently, whenever \eqref{eq:5.11} holds, $\operatorname{rank}\mathcal{S}\le1$ and every $2\times2$ minor of $\mathcal{S}$ vanishes. When $\mathcal{L}_{55}\ne0$, define the scalar coefficient $\mathcal{K}$ of the internal Schur complement by
\begin{equation}\label{eq:diagram5:schur-scalar}
\mathcal{K}=-\mathcal{L}_{34}-\frac{\mathcal{L}_{35}\mathcal{L}_{45}}{\mathcal{L}_{35}+\mathcal{L}_{45}}
=m_3m_4\delta_{34}+\frac{m_3m_4m_5\delta_{35}\delta_{45}}{m_3\delta_{35}+m_4\delta_{45}}.
\end{equation}

For the seven pairs not wholly contained in the internal triangle $S_3$, set
\[
\Theta_{\delta}:=\max\{1,|\delta_{12}|,|\delta_{13}|,|\delta_{14}|,|\delta_{15}|,|\delta_{23}|,|\delta_{24}|,|\delta_{25}|\}.
\]
The moduli make this a positive real scale, and the entry $1$ also bounds the constant terms in $\mathcal{R}$.

\begin{lemma}[Schur complement estimates]\label{lem:5.4}
For fixed positive masses, suppose
\[
\mathcal{L}_{55}\approx\epsilon^{-3},\qquad
\frac{\mathcal{L}_{k5}}{\mathcal{L}_{55}}=O(1)\quad(k\in T_2).
\]
If $\epsilon^3\Theta_{\delta}\to0$, then, for sufficiently small $\epsilon$, the scalar pivot $\Omega_{55}$ is nonzero, and $\mathcal{S}$ satisfies
\begin{equation}\label{eq:5.12}
\begin{gathered}
\mathcal{S}_{13}=O(\Theta_{\delta}),\quad \mathcal{S}_{23}=O(\Theta_{\delta}),\quad \mathcal{S}_{33}=\mathcal{K}+O(\Theta_{\delta}),\\
\mathcal{S}_{ij}=\Omega_{ij}+O(\epsilon^3\Theta_{\delta}^2)\quad(\text{for all }i,j\in S_2).
\end{gathered}
\end{equation}
Alternatively, retain the assumptions on $\mathcal{L}_{55}$ and the ratios $\mathcal{L}_{35}/\mathcal{L}_{55}$ and $\mathcal{L}_{45}/\mathcal{L}_{55}$, and suppose that all six cross coefficients $\delta_{ij}$, $i\in S_2$, $j\in S_3$, are bounded. The binary coefficient $\delta_{12}$ may be arbitrary; this variant does not require $\epsilon^3\Theta_{\delta}\to0$. The pivot $\Omega_{55}$ is again nonzero for sufficiently small $\epsilon$, and
\begin{equation}\label{eq:5.13}
\mathcal{S}_{13},\mathcal{S}_{23}=O(1),\qquad \mathcal{S}_{33}=\mathcal{K}+O(1),\qquad
\mathcal{S}_{11}=m_1m_2\delta_{12}+O(1).
\end{equation}
In either case, when \eqref{eq:5.11} holds, $\operatorname{rank}\mathcal{S}\le1$. The constants may depend on the fixed masses and the bounds in the hypotheses.
\end{lemma}
\begin{proof}
The Schur complement of $\mathcal{L}$ alone, on indices $1,2,3,4$, is exactly the $4\times4$ matrix
\[
\mathcal{K}
\begin{pmatrix}0\\0\\1\\-1\end{pmatrix}
\begin{pmatrix}0&0&1&-1\end{pmatrix}
=\mathcal{K}
\begin{pmatrix}
0&0&0&0\\
0&0&0&0\\
0&0&1&-1\\
0&0&-1&1
\end{pmatrix}.
\]
Direct elimination of the $55$ entry of $\mathcal{L}$ gives \eqref{eq:diagram5:schur-scalar}. For retained indices $k,l$, the exact difference between the two Schur complements is
\begin{equation}\label{eq:5.14}
\mathcal{R}_{kl}+\mathcal{L}_{k5}\mathcal{L}_{5l}(\mathcal{L}_{55}^{-1}-\Omega_{55}^{-1})
-\frac{\mathcal{L}_{k5}\mathcal{R}_{5l}+\mathcal{R}_{k5}\mathcal{L}_{5l}+\mathcal{R}_{k5}\mathcal{R}_{5l}}{\Omega_{55}}.
\end{equation}
All entries of $\mathcal{R}$ are $O(\Theta_{\delta})$. Since $\mathcal{R}_{55}/\mathcal{L}_{55}\to0$, the pivot is nonzero, and its internal ratios remain bounded. The terms in \eqref{eq:5.14} are $O(\Theta_{\delta})+O(\epsilon^3\Theta_{\delta}^2)=O(\Theta_{\delta})$. For two indices in $S_2$ the internal pivot entries vanish, giving the sharper last estimate in \eqref{eq:5.12}. The rank conclusion follows from the Schur complement rank identity above.

If all six cross coefficients $\delta_{ij}$, $i\in S_2$, $j\in S_3$, are bounded, then $\mathcal{R}_{55}=O(1)$, so $\Omega_{55}\sim\mathcal{L}_{55}\approx\epsilon^{-3}$ independently of $\delta_{12}$. The only entries of $\mathcal{R}$ containing $\delta_{12}$ are its $1,2$ principal block. This coefficient enters neither the pivot nor any entry used to obtain $\mathcal{S}_{13},\mathcal{S}_{23},\mathcal{S}_{33}$. Their remainders in \eqref{eq:5.14} are therefore bounded. The formula for $\mathcal{S}_{11}$, with its exterior mass factor $m_1m_2$, follows from the bounded cross entries and $\Omega_{55}^{-1}=O(\epsilon^3)$, independently of the size of $\delta_{12}$.
\end{proof}

\subsection{The case \texorpdfstring{$m_3=m_4$}{m3=m4}}\label{sec:diagram5:equal}

Assume $m_3=m_4>0$. The necessary mass relation gives $m_5=m_3/4$, while $m_1,m_2>0$ remain unrestricted.
Common mass scaling and the corresponding multiplier normalization reduce the proof to the masses
\begin{equation}
\begin{gathered}
m_3=m_4=1,\qquad m_5=1/4,\qquad M=m_1+m_2+9/4,\\
e_2(S_5)=m_1m_2+\tfrac94(m_1+m_2)+\tfrac32,
\end{gathered}
\label{eq:6.1}
\end{equation}
with arbitrary positive exterior masses $m_1,m_2$ and the internal labels $3,4,5$. We retain the same symbols for the rescaled masses. The total mass of the internal triangle is $e_1(S_3)=9/4$.
The physical lower bound inherited by a normalized component containing a physical point is \eqref{eq:3.1}.
The identities below use this equal-endpoint normalization. The matrix definitions of Section~\ref{sec:5.2} remain in force with these masses.

\subsubsection{Exact identities for the isolated \texorpdfstring{$zw$}{zw}-triangle}\label{sec:diagram5:equal-identities}
\label{sec:5.1}

\begin{lemma}[Coordinates for masses $(1,1,1/4)$]\label{lem:5.1}
Suppose the isolated $zw$-triangle of a diagram-5 limit has masses $(1,1,1/4)$ at vertices $3,4,5$, with arbitrary positive exterior masses. With the relative coordinates defined in Section~\ref{sec:2}, the exact positions satisfy
\begin{equation}\label{eq:5.1}
\begin{aligned}
q_{S_3,3}&=\frac12q_{34}-\frac19q_{T_2,5},\\
q_{S_3,4}&=-\frac12q_{34}-\frac19q_{T_2,5},\\
q_{S_3,5}&=\frac89q_{T_2,5}.
\end{aligned}
\end{equation}
where
\[
r_{34}^2=z_{34}w_{34},\quad r_{34},z_{34},w_{34}\approx\epsilon,\quad
\zeta_{T_2,5},\omega_{T_2,5}=o(\epsilon),\quad r_{35},r_{45}\sim-r_{34}/2.
\]
\end{lemma}
\begin{proof}
The dominant internal force terms at each vertex of the isolated $zw$-triangle in Diagram~5 sum to zero; see \cite[Sections~8.4--8.5]{AK}. The leading midpoint and signs follow directly from these equations. These equations make the leading $z,w$ differences between vertices in $S_3$ proportional, with nonisotropic differences by \eqref{eq:2.4}. In a complex scalar coordinate local to this proof, place vertices 3 and 4, each of mass one, at $1,-1$, and body 5, of mass $1/4$, at $\xi\in\mathbb C$; no reality assumption on $\xi$ is made. The signed distances have the form $2k\sigma_{34}$, $k(1-\xi)\sigma_{35}$, $k(-1-\xi)\sigma_{45}$, with $k\ne0$ and $\sigma_{ij}=\pm1$. The leading internal force equations at vertices 3 and 4 give
\[
(1-\xi)^2=-\sigma_{35}/\sigma_{34},\qquad (1+\xi)^2=\sigma_{45}/\sigma_{34}.
\]
Thus $1-\xi$ and $1+\xi$ lie in $\{1,-1,\mathrm{i},-\mathrm{i}\}$ and sum to two. Both are one, so $\xi=0$ and $r_{35},r_{45}\sim-r_{34}/2$. Using $m_3=m_4=1$ and the definition of $q_{T_2,5}$ gives \eqref{eq:5.1} exactly; the midpoint limit gives the stated smallness of $q_{T_2,5}$.
\end{proof}

In this normalization, Lemma~\ref{lem:5.1} gives
\[
\mathcal{L}_{55}=\frac{\delta_{35}+\delta_{45}}4\sim-4\delta_{34}\approx\epsilon^{-3},\qquad
\frac{\mathcal{L}_{35}}{\mathcal{L}_{55}},\frac{\mathcal{L}_{45}}{\mathcal{L}_{55}}\longrightarrow-\frac12.
\]
Thus the hypotheses of Lemma~\ref{lem:5.4} on $\mathcal{L}_{55}$ and the ratios $\mathcal{L}_{35}/\mathcal{L}_{55}$ and $\mathcal{L}_{45}/\mathcal{L}_{55}$ hold. In each application, we will check either $\epsilon^3\Theta_{\delta}\to0$ or boundedness of all six cross coefficients $\delta_{ij}$, $i\in S_2$, $j\in S_3$.

Put
\[
\phi_z=-\frac{2\zeta_{T_2,5}}{z_{34}},\qquad\phi_w=-\frac{2\omega_{T_2,5}}{w_{34}}.
\]
Thus $\phi_z,\phi_w\to0$ by Lemma~\ref{lem:5.1}. Under \eqref{eq:6.1}, the general scalar \eqref{eq:diagram5:schur-scalar} becomes
\begin{equation}\label{eq:5.3}
\begin{aligned}
\mathcal{K}&=\delta_{34}+\frac{\delta_{35}\delta_{45}}{4(\delta_{35}+\delta_{45})}\\
&=\delta_{34}\,\frac{\bigl((1+\phi_z)(1+\phi_w)\bigr)^{3/2}+\bigl((1-\phi_z)(1-\phi_w)\bigr)^{3/2}-2}
{\bigl((1+\phi_z)(1+\phi_w)\bigr)^{3/2}+\bigl((1-\phi_z)(1-\phi_w)\bigr)^{3/2}}.
\end{aligned}
\end{equation}
Its denominator is nonzero for sufficiently small $\epsilon$, since $\delta_{35},\delta_{45}\sim-8\delta_{34}$.

On the signed branches for the isolated $zw$-triangle in Lemma~\ref{lem:5.1}, the sums of the two dominant terms in the equations for vertex 5 are
\[
f_{5,w}^{S_3}=m_3W_{53}+m_4W_{54}=4W_{34}\,\hat f_{5,w},\qquad
f_{5,z}^{S_3}=m_3Z_{53}+m_4Z_{54}=4Z_{34}\,\hat f_{5,z}.
\]
The two summands in each component have order $\epsilon^{-2}$, while their sum is $o(\epsilon^{-2})$ because the leading terms cancel. Here $\hat f_{5,w},\hat f_{5,z}$ denote the normalized sums of these dominant terms, given exactly by
\begin{equation}\label{eq:5.4}
\hat f_{5,w}=(1+\phi_z)^{-3/2}(1+\phi_w)^{-1/2}
-(1-\phi_z)^{-3/2}(1-\phi_w)^{-1/2},
\end{equation}
\begin{equation}\label{eq:5.5}
\hat f_{5,z}=(1+\phi_z)^{-1/2}(1+\phi_w)^{-3/2}
-(1-\phi_z)^{-1/2}(1-\phi_w)^{-3/2}.
\end{equation}

Let $F$ be the analytic map near $(0,0)\in\mathbb{C}^2$ defined by \eqref{eq:5.4}--\eqref{eq:5.5}:
\[
F(\phi_z,\phi_w)
=\bigl(\hat f_{5,w}(\phi_z,\phi_w),\hat f_{5,z}(\phi_z,\phi_w)\bigr).
\]

\begin{lemma}[Local inverse and asymptotic estimates]\label{lem:5.3}
All powers use their analytic branches near one. The odd analytic map $F$ has an odd analytic local inverse $F^{-1}$ near zero, satisfying $F^{-1}(\hat f_{5,w},\hat f_{5,z})=(\phi_z,\phi_w)$. Its component expansions are
\begin{equation}\label{eq:5.6}
\phi_z=\frac{-3\hat f_{5,w}+\hat f_{5,z}}{8}+O((|\hat f_{5,w}|+|\hat f_{5,z}|)^3),
\end{equation}
\begin{equation}\label{eq:5.7}
\phi_w=\frac{\hat f_{5,w}-3\hat f_{5,z}}{8}+O((|\hat f_{5,w}|+|\hat f_{5,z}|)^3).
\end{equation}
The map has no parameters, so the neighbourhood of zero and the constant in the remainder are universal.
Moreover,
\begin{equation}\label{eq:5.8}
\begin{gathered}
\bigl((1+\phi_z)(1+\phi_w)\bigr)^{3/2}+\bigl((1-\phi_z)(1-\phi_w)\bigr)^{3/2}-2\\
=\tfrac34(\phi_z^2+\phi_w^2)+\tfrac92\phi_z\phi_w+O\bigl((|\phi_z|+|\phi_w|)^4\bigr).
\end{gathered}
\end{equation}
The following consequences hold.
\begin{enumerate}
\item If $\Theta_{f_{5,z}}=o(\epsilon^{-2})$, $f_{5,z}^{S_3}\approx \Theta_{f_{5,z}}$, and $f_{5,w}^{S_3}=o(\Theta_{f_{5,z}})$, then
\begin{equation}\label{eq:5.9}
\phi_z\approx\epsilon^2\Theta_{f_{5,z}},\qquad\phi_w/\phi_z\to-3,\qquad
\mathcal{K}\approx\epsilon \Theta_{f_{5,z}}^2.
\end{equation}
\item If $f_{5,w}^{S_3}=O(\epsilon)$ and $\mathcal{K}=O(1)$, then $|\phi_z|+|\phi_w|=O(\epsilon^{3/2})$.
\end{enumerate}
\end{lemma}
\begin{proof}
The exact positions relative to vertex 5 are $q_{35}=q_{34}/2-q_{T_2,5}$ and $q_{45}=-q_{34}/2-q_{T_2,5}$. Therefore, on the signed branches of Lemma~\ref{lem:5.1},
\[
\begin{aligned}
r_{35}&=-\frac{r_{34}}{2}\sqrt{(1+\phi_z)(1+\phi_w)},\\
r_{45}&=-\frac{r_{34}}{2}\sqrt{(1-\phi_z)(1-\phi_w)},
\end{aligned}
\]
where both square roots use the analytic branch near one with $\sqrt{1}=1$. The contribution of vertices 3 and 4 to the right-hand side of the equation for vertex 5 is
\[
f_5^{S_3}=\frac12(\delta_{45}-\delta_{35})q_{34}+(\delta_{35}+\delta_{45})q_{T_2,5}.
\]
Substitution gives \eqref{eq:5.4}--\eqref{eq:5.5}. Their linear terms are $-3\phi_z-\phi_w$ and $-\phi_z-3\phi_w$. The derivative matrix $\left(\begin{smallmatrix}-3&-1\\-1&-3\end{smallmatrix}\right)$ has determinant eight. The analytic inverse theorem gives the local inverse $F^{-1}$. Since $F(-\phi_z,-\phi_w)=-F(\phi_z,\phi_w)$, uniqueness of the inverse implies that $F^{-1}$ is odd. Its expansion has no even terms, giving \eqref{eq:5.6}--\eqref{eq:5.7}. Substituting $\delta_{35}/\delta_{34}=-8\bigl((1+\phi_z)(1+\phi_w)\bigr)^{-3/2}$ and $\delta_{45}/\delta_{34}=-8\bigl((1-\phi_z)(1-\phi_w)\bigr)^{-3/2}$ into \eqref{eq:5.3} proves its exact expression. Taylor expansion gives the displayed quadratic part; the left side of \eqref{eq:5.8} is invariant under $\phi\mapsto-\phi$, so its odd-order terms vanish.

Since $Z_{34},W_{34}\approx\epsilon^{-2}$, the normalization identities $\hat f_{5,z}=f_{5,z}^{S_3}/(4Z_{34})$ and $\hat f_{5,w}=f_{5,w}^{S_3}/(4W_{34})$ give, under the first set of hypotheses, $\hat f_{5,z}\approx\epsilon^2\Theta_{f_{5,z}}\to0$ and $\hat f_{5,w}/\hat f_{5,z}\to0$. The cubic inverse remainders divided by $\hat f_{5,z}$ tend to zero. Thus $\phi_z\sim \hat f_{5,z}/8$ and $\phi_w/\phi_z\to-3$. The left side of \eqref{eq:5.8} is consequently $-6\phi_z^2(1+o(1))$, and \eqref{eq:5.3} gives $\mathcal{K}\sim-3\phi_z^2\delta_{34}$, so \eqref{eq:5.9} follows.

For the second assertion, $\hat f_{5,w}=f_{5,w}^{S_3}/(4W_{34})=O(\epsilon^3)$. If $(|\phi_z|+|\phi_w|)/\epsilon^{3/2}$ were unbounded on a subsequence, $\epsilon^3=o(|\phi_z|+|\phi_w|)$ and \eqref{eq:5.4} would give $\phi_w+3\phi_z=o(|\phi_z|+|\phi_w|)$, hence $\phi_w/\phi_z\to-3$ and $\phi_z\approx|\phi_z|+|\phi_w|$. Equations \eqref{eq:5.3} and \eqref{eq:5.8} would then imply $\mathcal{K}\approx(|\phi_z|+|\phi_w|)^2/\epsilon^3\to\infty$, contradicting its boundedness.
\end{proof}

\subsubsection{Reduction to a three-point limit}\label{sec:6.2}
\label{subsec:6.2}

We consider Case~\ref{case:6.A} of the classification in Section~\ref{sec:6.1}, so $\zeta_{S_3,1}\approx\zeta_{S_3,2}\approx\zeta_{S_3,1}-\zeta_{S_3,2}\approx\Theta_{\zeta}$. By \eqref{eq:6.5}, $\epsilon\Theta_{\zeta}\to\infty$ and each exterior body has at most one attachment. All cross $z$-separations and $z_{12}$ have order $\Theta_{\zeta}$, while $z_1,z_2=O(\Theta_{\zeta})$. If at least one attachment collapses, apply the equations for the exterior bodies
\[
z_i=m_{3-i}z_{i,3-i}\delta_{12}+\sum\limits_{k\in S_3}m_kz_{ik}\delta_{ik},\qquad i\in S_2.
\]
There are three alternatives to two collapsing attachments at the same vertex.
\begin{enumerate}[label=\textup{(\roman*)},leftmargin=*]
\item If only one exterior body has an attachment, the other equation has cross sum $o(\Theta_{\zeta})$, hence $\delta_{12}=O(1)$. The equation for the attached exterior body then bounds its sole attachment coefficient, contradicting collapse.
\item If the two attachments meet distinct vertices $j,k\in S_3$, their $w$-separations are $O(\Theta_{\zeta}^{-1})=o(\epsilon)$. Since $w_{jk}\approx\epsilon$, we get $w_{12}\approx\epsilon$ and $|r_{12}|^2\approx\Theta_{\zeta}\epsilon\to\infty$. Thus $\delta_{12}\to0$, and the equation for each exterior body again bounds its attachment coefficient, contradicting collapse.
\item If both attachments meet the same vertex and exactly one has a finite nonzero limit, the equation for that exterior body first gives $\delta_{12}=O(1)$. The other equation then bounds the collapsing attachment coefficient, a contradiction.
\end{enumerate}
The only remaining possibility is two collapsing attachments to the same vertex $j\in S_3$. Apply Lemma~\ref{lem:5.2} to both clusters with reference body $j$. The four other cross distances diverge, so the individual equations give
\[
f_{k,w}^{S_3}=w_k-\sum\limits_{i\in S_2}m_iW_{ki}=O(\epsilon),\qquad k\in S_3\setminus\{j\}.
\]
Since $z_{kj}=O(\epsilon)$ for these two vertices $k\in S_3\setminus\{j\}$, Lemma~\ref{lem:5.2} yields
\[
U_{S_3}=\sum\limits_{k\in S_3\setminus\{j\}}m_kz_{kj}f_{k,w}^{S_3}=O(\epsilon^2).
\]
For the cluster $S_2\cup\{j\}$, the equations for the exterior bodies similarly give
\[
f_{i,z}^{S_2\cup\{j\}}=z_i-\sum\limits_{k\in S_3\setminus\{j\}}m_kZ_{ik}=O(\Theta_{\zeta}),\qquad i\in S_2.
\]
The collapsing attachments have $w_{ij}=r_{ij}^2/z_{ij}=o(1/\Theta_{\zeta})$, hence the second identity in Lemma~\ref{lem:5.2} gives
\[
U_{S_2\cup\{j\}}=\sum\limits_{i\in S_2}m_iw_{ij}f_{i,z}^{S_2\cup\{j\}}=o(1).
\]
The terms corresponding to vertex $j$ vanish exactly because their coordinate differences are zero; no bound on $f_{j,w}^{S_3}$ or $f_{j,z}^{S_2\cup\{j\}}$ is needed. All other cross potentials vanish, and hence $U\to0$, contradicting \eqref{eq:3.1}.

If no attachment collapses, all cross coefficients are bounded. The exact $z$ equations for the exterior bodies yield
\begin{equation}
\begin{pmatrix}
\begin{gathered}1-m_1/M-m_2\delta_{12}\\[-2pt]-\sum\limits_{j\in S_3}m_j\delta_{1j}\end{gathered}&m_2(\delta_{12}-1/M)\\[4pt]
m_1(\delta_{12}-1/M)&\begin{gathered}1-m_2/M-m_1\delta_{12}\\[-2pt]-\sum\limits_{j\in S_3}m_j\delta_{2j}\end{gathered}
\end{pmatrix}\begin{pmatrix}\zeta_{S_3,1}\\\zeta_{S_3,2}\end{pmatrix}=O(\epsilon),
\label{eq:6.7}
\end{equation}
Since $z_{12}\approx\Theta_{\zeta}$, either original $z$ equation for the exterior bodies bounds $\delta_{12}$. All seven exterior or cross distances $r_{12},r_{13},r_{14},r_{15},r_{23},r_{24},r_{25}$ are therefore bounded away from zero, and their $z$-separations have modulus comparable to $\Theta_{\zeta}$. Thus the pair-component identity in Section~\ref{sec:2.2} gives
\[
W_{ij}=O(\Theta_{\zeta}^{-1}),\qquad i,j\in S_5,\quad i<j,\quad \{i,j\}\cap S_2\ne\varnothing.
\]
The individual $w$ equations for the exterior bodies and total centering then imply
\[
w_1,w_2,C_{S_3,w}=O(\Theta_{\zeta}^{-1})=o(\epsilon),
\]
since $\Theta_{\zeta}\epsilon\to\infty$. Combining this estimate with Lemma~\ref{lem:5.1}, we obtain, for $i\in S_2$ and $j\in T_2$,
\[
z_{ij}\approx\Theta_{\zeta},\qquad w_{ij}\approx\epsilon,\qquad
r_{ij}^2=z_{ij}w_{ij}\approx\Theta_{\zeta}\epsilon\longrightarrow\infty.
\]
All attachments therefore go to vertex 5. For later use, the identities for $f_5^{S_3}$ from Section~\ref{sec:2.1} take the exact form
\[
\begin{aligned}
f_{5,z}^{S_3}&=z_5+m_1Z_{15}+m_2Z_{25},\\
f_{5,w}^{S_3}&=w_5+m_1W_{15}+m_2W_{25}.
\end{aligned}
\]
Here the plus signs follow from $Z_{5i}=-Z_{i5}$ and $W_{5i}=-W_{i5}$ for $i\in S_2$.

All six cross coefficients $\delta_{ij}$, $i\in S_2$, $j\in S_3$, are now bounded, and a finite attachment has $z$-separation comparable to $\Theta_{\zeta}$. Hence \eqref{eq:5.11} and the estimates \eqref{eq:5.13} of Lemma~\ref{lem:5.4} apply in both cases below.

If $\delta_{12}\to0$, exactly one attachment exists; relabel it $1$--$5$. The limiting kernel of \eqref{eq:6.7} gives
\[
\begin{gathered}
\sum\limits_{j\in S_3}m_j\delta_{1j}\to\frac{e_1(S_3)}{e_1(T_4)},\qquad \delta_{15}\to\frac{4e_1(S_3)}{e_1(T_4)},\\
\frac{\zeta_{S_3,2}}{\zeta_{S_3,1}}\to\frac{m_1}{e_1(T_4)},\qquad\frac{C_{S_3,z}}{\zeta_{S_3,1}}\to-\frac{m_1}{e_1(T_4)}.
\end{gathered}
\]
Consequently, the identities for $f_5^{S_3}$ above yield
\[
f^{S_3}_{5,z}/\zeta_{S_3,1}\to\frac{m_1(4e_1(S_3)-1)}{e_1(T_4)}=\frac{8m_1}{e_1(T_4)}\ne0,
\qquad f^{S_3}_{5,w}=O(\epsilon).
\]
The diagonal entry $\mathcal{S}_{11}$ has the limit
\[
\mathcal{S}_{11}\longrightarrow m_1\left(\frac{e_1(S_3)}{e_1(T_4)}-\frac{M-m_1}{M}\right)
=-\frac{m_1^2m_2}{Me_1(T_4)}\ne0.
\]
The Schur rank-one minor bounds $\mathcal{K}$, while Lemma~\ref{lem:5.3}(1) makes its modulus comparable to $\epsilon \Theta_{\zeta}^2\to\infty$, a contradiction.

If $\delta_{12}$ tends to a nonzero limit, the binary has a finite nonzero distance. At least one attachment goes to vertex 5; its $z$-separation has order $\Theta_{\zeta}$, so $w_1,w_2=O(1/\Theta_{\zeta})$ imply $w_5=O(1/\Theta_{\zeta})$. The distances $r_{15},r_{25}$ are consequently bounded and, since no attachment collapses, have nonzero limits. Equation~\eqref{eq:6.6} gives $U_{S_3}\to0$. Write $\delta_{ij}^0=\lim\limits_{\epsilon\to0}\delta_{ij}$ for $ij\in\{12,15,25\}$. For $k,l\in S_2$, write $\mathcal{S}_{kl}^0=\lim\limits_{\epsilon\to0}\mathcal{S}_{kl}$. The limiting upper Schur block has entries
\[
\begin{aligned}
\mathcal{S}_{11}^0&=m_1m_2\delta_{12}^0+\tfrac{m_1}4\delta_{15}^0-\frac{m_1(M-m_1)}{M},\\
\mathcal{S}_{22}^0&=m_1m_2\delta_{12}^0+\tfrac{m_2}4\delta_{25}^0-\frac{m_2(M-m_2)}{M},\\
\mathcal{S}_{12}^0&=-m_1m_2\left(\delta_{12}^0-1/M\right).
\end{aligned}
\]
If this block is zero, then $\delta_{12}^0=1/M$ and $\delta_{15}^0=\delta_{25}^0=4e_1(S_3)/M$. The four potential terms with $i\in S_2$, $j\in T_2$ tend to zero, and $U_{S_3}\to0$, so the remaining three potential terms satisfy
\[
M^{1/3}U\le m_1m_2+(m_1+m_2)\frac{e_1(S_3)^{1/3}}{4^{2/3}}<e_2(S_5).
\]
Here $e_1(S_3)^{1/3}/4^{2/3}<e_1(S_3)$. This contradicts \eqref{eq:3.1} and excludes the zero-block branch directly.

Suppose the block is nonzero. It has rank at most one, so one diagonal entry is nonzero; its Schur minor then bounds $\mathcal{K}$. Since $f_{5,w}^{S_3}=O(\epsilon)$, Lemma~\ref{lem:5.3}(2) gives
\[
\Theta_{\zeta} \omega_{T_2,5}\longrightarrow0,\qquad \zeta_{T_2,5}/\Theta_{\zeta}\longrightarrow0.
\]
We now derive the limiting equations from the original five-body equations. For $\{i,j\}=S_2$, multiplication of the equations for the exterior bodies by $m_i$ gives
\[
\begin{aligned}
m_iz_i&=m_1m_2Z_{ij}+m_i\bigl(Z_{i3}+Z_{i4}+\tfrac14Z_{i5}\bigr),\\
m_iw_i&=m_1m_2W_{ij}+m_i\bigl(W_{i3}+W_{i4}+\tfrac14W_{i5}\bigr).
\end{aligned}
\]
The mass-weighted cluster identities in Section~\ref{sec:2.1} give
\[
\begin{aligned}
e_1(S_3)C_{S_3,z}&=-\sum\limits_{i\in S_2}m_i\bigl(Z_{i3}+Z_{i4}+\tfrac14Z_{i5}\bigr),\\
e_1(S_3)C_{S_3,w}&=-\sum\limits_{i\in S_2}m_i\bigl(W_{i3}+W_{i4}+\tfrac14W_{i5}\bigr).
\end{aligned}
\]
The terms for each pair of vertices in $S_3$ cancel exactly; the minus signs arise from reversing the indices of the cross components.
For $i\in S_5$, set \[\widehat{z}_i:=\frac{z_i}{\Theta_{\zeta}},\qquad \widehat{w}_i:=\Theta_{\zeta}w_i.\] Let $\widehat{q}_i$ be the position with these $(z,w)$ components and $\widehat{C}_{S_3}$ the center of mass of the positions indexed by $S_3$. Write $\widehat{z}_{ij}=\widehat{z}_i-\widehat{z}_j$ and similarly for $\widehat{w}_{ij}$. The signed distances and the form of the normalized equations are unchanged. The exterior positions $\widehat{q}_1,\widehat{q}_2$ and the center of mass $\widehat{C}_{S_3}$ are bounded. For $i\in S_2$ and $j\in T_2$, the previously established bounds $z_{ij}\approx\Theta_{\zeta}$ and $r_{ij}\to\infty$ give
\[
\begin{aligned}
\frac{\widehat{z}_{ij}}{r_{ij}^3}
&=\frac{Z_{ij}}{\Theta_{\zeta}}
=\frac{z_{ij}}{\Theta_{\zeta}r_{ij}^3}
=O(r_{ij}^{-3})\longrightarrow0,\\
\frac{\widehat{w}_{ij}}{r_{ij}^3}
&=\Theta_{\zeta}W_{ij}
=\frac{\Theta_{\zeta}}{r_{ij}z_{ij}}
=O(r_{ij}^{-1})\longrightarrow0.
\end{aligned}
\]
Thus the four interactions with vertices 3 and 4 disappear from the transformed limiting equations. Moreover, $q_{S_3,5}=8q_{T_2,5}/9$ and the transformed components of $q_{T_2,5}$ tend to zero. After taking a subsequence, write
\[
q_i^0=\lim\limits_{\epsilon\to0} \widehat{q}_i\quad(i\in S_2),\qquad
q_0^0=\lim\limits_{\epsilon\to0} \widehat{C}_{S_3}=\lim\limits_{\epsilon\to0} \widehat{q}_5,
\]
and denote the finite nonzero signed distances by
\[
r_{12}^0=\lim\limits_{\epsilon\to0} r_{12},\qquad
r_{10}^0=\lim\limits_{\epsilon\to0} r_{15},\qquad
r_{20}^0=\lim\limits_{\epsilon\to0} r_{25}.
\]
Writing $z_k^0,w_k^0$ for the components of $q_k^0$, these distances satisfy
\[
(r_{kl}^0)^2=(z_k^0-z_l^0)(w_k^0-w_l^0),\qquad
r_{kl}^0=r_{lk}^0\ne0
\quad(k,l\in\{0,1,2\},\ k\ne l).
\]
Taking limits in these mass-weighted equations gives
\begin{equation}
\begin{aligned}
m_1q_1^0&=m_1m_2\frac{q_1^0-q_2^0}{(r_{12}^0)^3}
+\frac{m_1}{4}\frac{q_1^0-q_0^0}{(r_{10}^0)^3},\\
m_2q_2^0&=m_1m_2\frac{q_2^0-q_1^0}{(r_{12}^0)^3}
+\frac{m_2}{4}\frac{q_2^0-q_0^0}{(r_{20}^0)^3},\\
e_1(S_3)q_0^0&=\frac{m_1}{4}\frac{q_0^0-q_1^0}{(r_{10}^0)^3}
+\frac{m_2}{4}\frac{q_0^0-q_2^0}{(r_{20}^0)^3}.
\end{aligned}
\label{eq:6.8}
\end{equation}
Their sum gives $m_1q_1^0+m_2q_2^0+e_1(S_3)q_0^0=0$, and their multiplier remains one. The three interaction coefficients are the surviving original mass products $m_1m_2$, $m_1m_5=m_1/4$, and $m_2m_5=m_2/4$; the coefficient $e_1(S_3)=9/4$ on the left comes from summing the three masses indexed by $S_3$.

The original five-body potential has the exact decomposition
\[
U=\frac{m_1m_2}{r_{12}}
+\sum\limits_{i\in S_2}\left(\frac{m_i}{r_{i3}}+\frac{m_i}{r_{i4}}+\frac{m_i}{4r_{i5}}\right)+U_{S_3}.
\]
The four potential terms with $i\in S_2$, $j\in T_2$ vanish by distance divergence, and \eqref{eq:6.6} gives $U_{S_3}\to0$. Define the potential of the limiting points by
\begin{equation}
U^0=\frac{m_1m_2}{r_{12}^0}+\frac{m_1}{4r_{10}^0}+\frac{m_2}{4r_{20}^0}.
\label{eq:6.9}
\end{equation}
The denominators are signed complex distances. For the limit just constructed, $U^0$ equals the original constant $U$ and must therefore be positive real with $U^0\ge e_2(S_5)/M^{1/3}$. Here $M$ and $e_2(S_5)$ still refer to the five masses in \eqref{eq:6.1}. The next subsection excludes this remaining possibility and completes Case~\ref{case:6.A}.

\subsubsection{The case when \texorpdfstring{$\zeta_{S_3,1}\approx\zeta_{S_3,2}\approx z_{12}$}{zeta(S3,1) comparable to zeta(S3,2) comparable to z12}}\label{sec:6.3}
\label{subsec:6.3}

The preceding reduction specifies the three-point equations and their potential explicitly. We prove that every complex solution of that system, on every signed-distance branch, is incompatible with the physical lower bound \eqref{eq:3.1}. Here the superscript $0$ labels quantities in the limiting three-point system; an arbitrary solution of this system need not arise from a five-body degeneration.

\begin{lemma}\label{lem:6.2}
Let $m_1,m_2>0$ and $e_1(S_3)=9/4$, with $M$ and $e_2(S_5)$ as in \eqref{eq:6.1}. Let three complex points $q_1^0,q_2^0,q_0^0$ satisfy \eqref{eq:6.8}, with symmetric nonzero signed distances satisfying $(r_{kl}^0)^2=(z_k^0-z_l^0)(w_k^0-w_l^0)$. Their potential $U^0$, defined in \eqref{eq:6.9}, cannot be a positive real number with $U^0\ge e_2(S_5)/M^{1/3}$.
\end{lemma}

\begin{proof}
For a noncollinear configuration, independence of the two relative vectors at each point gives
\[
(r_{12}^0)^{-3}=1/M,\qquad (r_{10}^0)^{-3}=(r_{20}^0)^{-3}=4e_1(S_3)/M.
\]
Consequently
\[
M^{1/3}|U^0|\le m_1m_2+(m_1+m_2)\frac{e_1(S_3)^{1/3}}{4^{2/3}}<e_2(S_5),
\]
since $e_1(S_3)^{1/3}/4^{2/3}<e_1(S_3)$.

For a collinear configuration, use the components $z_i^0,w_i^0$ of the normalized limiting points $q_i^0$ in \eqref{eq:6.8}, with $i\in\{0,1,2\}$. Since $(r_{12}^0)^2=(z_2^0-z_1^0)(w_2^0-w_1^0)\ne0$, both coordinate differences are nonzero. Collinearity therefore allows us to define
\[
\xi=\frac{z_0^0-z_1^0}{z_2^0-z_1^0}
 =\frac{w_0^0-w_1^0}{w_2^0-w_1^0},
\qquad \xi\in\mathbb C\setminus\{0,1\}.
\]
Here $\xi$ is a dimensionless shape parameter, while the nonzero signed distance $r_{12}^0$ retains the overall distance scale. The signed-distance identities give
\[
r_{10}^0=\sigma_{10}\xi r_{12}^0,\qquad
r_{20}^0=\sigma_{20}(\xi-1)r_{12}^0,\qquad
\sigma_{10},\sigma_{20}\in\{1,-1\}.
\]
The two signs specify the distance branches relative to $r_{12}^0$; all signed-distance choices are included. Substituting these identities into the $z$ components of \eqref{eq:6.8} gives
\[
\begin{aligned}
z_1^0&=\frac{z_2^0-z_1^0}{(r_{12}^0)^3}
\left(-m_2-\frac{\sigma_{10}}{4\xi^2}\right),\\
z_2^0&=\frac{z_2^0-z_1^0}{(r_{12}^0)^3}
\left(m_1-\frac{\sigma_{20}}{4(\xi-1)^2}\right),\\
z_0^0&=\frac{z_2^0-z_1^0}{4e_1(S_3)(r_{12}^0)^3}
\left(\frac{\sigma_{10}m_1}{\xi^2}+\frac{\sigma_{20}m_2}{(\xi-1)^2}\right).
\end{aligned}
\]
Subtracting the equation for point $1$ from that for point $2$ and cancelling $z_2^0-z_1^0\ne0$ yields
\[
(r_{12}^0)^3=m_1+m_2+\frac{\sigma_{10}}{4\xi^2}
-\frac{\sigma_{20}}{4(\xi-1)^2}\ne0.
\]
Subtracting the equation for point $1$ from that for point $0$, using $z_0^0-z_1^0=\xi(z_2^0-z_1^0)$, and eliminating $(r_{12}^0)^3$ by the preceding identity gives, with $e_1(S_3)=9/4$, the quintic equation
\begin{equation}
\begin{aligned}
g_5(\xi)={}&-36(m_1+m_2)\xi^5\\
&+(72m_1+108m_2)\xi^4\\
&+\bigl(-36m_1-108m_2-9\sigma_{10}+9\sigma_{20}\bigr)\xi^3\\
&+\bigl[36m_2+(27+4m_1)\sigma_{10}+4m_2\sigma_{20}\bigr]\xi^2\\
&-(27+8m_1)\sigma_{10}\xi\\
&+(9+4m_1)\sigma_{10}=0.
\end{aligned}
\label{eq:6.10}
\end{equation}
Its leading coefficient is $-36(m_1+m_2)\ne0$; $g_5(0)=\sigma_{10}(9+4m_1)\ne0$ and $g_5(1)=\sigma_{20}(9+4m_2)\ne0$. Define
\begin{equation}
\begin{aligned}
g_2(\xi)={}&4(m_1+m_2)\xi^4-8(m_1+m_2)\xi^3\\
&+\bigl[4(m_1+m_2)+\sigma_{10}-\sigma_{20}\bigr]\xi^2\\
&-2\sigma_{10}\xi+\sigma_{10},
\end{aligned}
\label{eq:6.11}
\end{equation}
\begin{equation}
\begin{aligned}
g_3(\xi)={}&4m_1m_2\xi^2\\
&+\bigl(-4m_1m_2+\sigma_{10}m_1+\sigma_{20}m_2\bigr)\xi
-\sigma_{10}m_1.
\end{aligned}
\label{eq:6.12}
\end{equation}
The scale identity above and the potential definition \eqref{eq:6.9} give
\[
(r_{12}^0)^3=\frac{g_2(\xi)}{4\xi^2(\xi-1)^2}\ne0,\qquad
U^0=\frac{g_3(\xi)}{4\xi(\xi-1)r_{12}^0}.
\]
Cubing the potential and substituting the scale identity, we obtain
\begin{equation}
(U^0)^3=\frac{g_3(\xi)^3}{16\xi(\xi-1)g_2(\xi)}.
\label{eq:6.13}
\end{equation}
This expression is valid on every root of \eqref{eq:6.10} with $g_2\ne0$, whether or not the shape is real.

Suppose, for a contradiction, that the potential is positive real and satisfies
$U^0\ge e_2(S_5)/M^{1/3}$. Set
\[
\Delta_U=(U^0)^3-\frac{e_2(S_5)^3}{M}\ge0.
\]
For the algebraic calculation, regard $\Delta_U$ as an independent variable and define
\begin{equation}
g_6(\xi)=\{64M\Delta_U+64e_2(S_5)^3\}\xi(\xi-1)g_2(\xi)-4M g_3(\xi)^3.
\label{eq:6.14}
\end{equation}
At the value of $\Delta_U$ just specified, \eqref{eq:6.13} gives $g_6(\xi)=0$. Thus the shape parameter is a common root of $g_5$ and $g_6$.

For these polynomials, Proposition~\ref{prop:app-energy} in Appendix~\ref{app:energy} gives a polynomial
$\mathcal{D}_{\sigma_{10},\sigma_{20}}\in\mathbb Z[m_1,m_2,\Delta_U]$, normalized by
\eqref{eq:A.4}, such that the common-root condition implies
\begin{equation}
\mathcal{G}_{\sigma_{10},\sigma_{20}}(m_1,m_2)
\mathcal{D}_{\sigma_{10},\sigma_{20}}(m_1,m_2,\Delta_U)=0,
\label{eq:6.elimination}
\end{equation}
where
\[
\mathcal{G}_{\sigma_{10},\sigma_{20}}(m_1,m_2)
=(-\sigma_{10}m_1+\sigma_{20}m_2-4m_1m_2)^2
+4\sigma_{10}\sigma_{20}m_1m_2.
\]
The common-root implication is valid for positive $m_1,m_2$, every complex $\Delta_U$,
and all four sign combinations.
The appendix obtains it by eliminating the shape variable with a B\'ezout matrix and factoring
its determinant. Exact coefficient verification also establishes that every nonzero coefficient
of $\mathcal{D}_{\sigma_{10},\sigma_{20}}$ is positive and that its constant term in $\Delta_U$ is nonzero.
Consequently,
\begin{equation}
\mathcal{D}_{\sigma_{10},\sigma_{20}}(m_1,m_2,\Delta_U)>0
\quad\text{for }m_1,m_2>0,\ \Delta_U\ge0.
\label{eq:6.17}
\end{equation}

If $\mathcal{G}_{\sigma_{10},\sigma_{20}}(m_1,m_2)\ne0$, \eqref{eq:6.elimination} gives
$\mathcal{D}_{\sigma_{10},\sigma_{20}}(m_1,m_2,\Delta_U)=0$, contradicting \eqref{eq:6.17}.
It remains to treat positive parameters with $\mathcal{G}_{\sigma_{10},\sigma_{20}}=0$.
Fix the corresponding shape root $\check{\xi}$, with $\check{\xi}\ne0,1$ and $g_2(\check{\xi})\ne0$.
Since $\mathcal{G}_{\sigma_{10},\sigma_{20}}$ is a nonzero polynomial in the masses, choose positive parameter pairs
$(m_{1,N},m_{2,N})$, indexed by $N=1,2,\ldots$, tending to $(m_1,m_2)$ as $N\to\infty$ and satisfying
$\mathcal{G}_{\sigma_{10},\sigma_{20}}(m_{1,N},m_{2,N})\ne0$.
The quintic $g_5$ has fixed degree and a nonzero leading coefficient near the limiting masses.
Its roots can therefore be chosen $\check{\xi}_N\to \check{\xi}$ as $N\to\infty$, for example by applying the argument principle on a small circle about $\check{\xi}$.
For sufficiently large $N$, the roots still satisfy $\check{\xi}_N\ne0,1$ and $g_2(\check{\xi}_N)\ne0$, with the polynomials evaluated at $(m_{1,N},m_{2,N})$.
The values $(U_N^0)^3$ given by \eqref{eq:6.13} at these roots and masses converge to $(U^0)^3$. Here the superscript $0$ continues to label the three-point limiting system; $N$ indexes a separate perturbation of its masses and shape, and does not assert the existence of a five-body degeneration for each perturbed system. Only the cubed potentials are needed. Define
\[
\Delta_{U,N}=(U_N^0)^3-\frac{\left(m_{1,N}m_{2,N}+\frac94(m_{1,N}+m_{2,N})+\frac32\right)^3}{m_{1,N}+m_{2,N}+\frac94}.
\]
Then $\Delta_{U,N}\to \Delta_U$ as $N\to\infty$, and the corresponding polynomials $g_5,g_6$ have the common root $\check{\xi}_N$.
Applying \eqref{eq:6.elimination} to these perturbed parameters, with the same fixed signs, and using
$\mathcal{G}_{\sigma_{10},\sigma_{20}}(m_{1,N},m_{2,N})\ne0$ yields
\[
\mathcal{D}_{\sigma_{10},\sigma_{20}}(m_{1,N},m_{2,N},\Delta_{U,N})=0.
\]
The intermediate values $\Delta_{U,N}$ may be complex. Polynomial continuity nevertheless gives
$\mathcal{D}_{\sigma_{10},\sigma_{20}}(m_1,m_2,\Delta_U)=0$ in the limit, again contradicting \eqref{eq:6.17}, since the original $\Delta_U$ is real and nonnegative.
This proves the lemma for all positive $m_1,m_2$ and all four sign combinations.
\end{proof}

The limit obtained in Section~\ref{subsec:6.2} satisfies the equations of Lemma~\ref{lem:6.2} and has $U^0=U\ge e_2(S_5)/M^{1/3}$. The lemma gives a contradiction, completing Case~\ref{case:6.A}.

\subsubsection{The case when \texorpdfstring{$\zeta_{S_3,1}\sim\zeta_{S_3,2}\succ z_{12}$}{zeta(S3,1) asymptotic to zeta(S3,2), with smaller z12}}\label{sec:6.4}
\label{subsec:6.4}

We now consider Case~\ref{case:6.B}: $\zeta_{S_3,1}\approx\zeta_{S_3,2}\approx\Theta_{\zeta}$ and $\zeta_{S_3,1}-\zeta_{S_3,2}=o(\Theta_{\zeta})$. Then $z_2\sim z_1$, $C_{S_3,z}\sim-(m_1+m_2)z_1/e_1(S_3)$, and $z_{ij}\sim Mz_1/e_1(S_3)$ for $i\in S_2$ and $j\in S_3$. The weighted sum of the $z$ equations for bodies 1 and 2 cancels the binary force exactly:
\[
m_1z_1+m_2z_2
=\sum\limits_{\substack{i\in S_2\\l\in S_3}}m_im_lZ_{il}.
\]
At least one attachment exists. Exactly one collapsing attachment cannot be balanced by finite terms. For two collapsing attachments $1$--$j$, $2$--$k$ with $j,k\in S_3$, every other cross term is $O(\Theta_{\zeta})$, whereas
\[
\frac{Z_{1j}}{\Theta_{\zeta}}
=\frac{z_{1j}}{\Theta_{\zeta}}\delta_{1j}
\longrightarrow\infty.
\]
Dividing the weighted identity by $Z_{1j}$, and then applying the same normalization to the body-1 equation
\[
z_1=m_2Z_{12}+m_jZ_{1j}
+\sum\limits_{l\in S_3\setminus\{j\}}m_lZ_{1l},
\]
we obtain, respectively,
\[
\frac{Z_{2k}}{Z_{1j}}\longrightarrow-\frac{m_1m_j}{m_2m_k},\qquad
\frac{Z_{12}}{Z_{1j}}\longrightarrow-\frac{m_j}{m_2}.
\]
Since $z_{2k}/z_{1j}\to1$ and $z_{12}=o(\Theta_{\zeta})$, these limits also give
\[
\frac{\delta_{2k}}{\delta_{1j}}
=\frac{Z_{2k}}{Z_{1j}}\frac{z_{1j}}{z_{2k}}
\longrightarrow-\frac{m_1m_j}{m_2m_k},\qquad
\frac{\delta_{12}}{\delta_{1j}}
=\frac{Z_{12}}{Z_{1j}}
 \frac{z_{1j}}{z_{12}}
\longrightarrow\infty.
\]

Suppose $j=k$. The preceding limit $\delta_{2j}/\delta_{1j}\to-m_1/m_2$ gives $r_{2j}\approx r_{1j}$. Together with $z_{2j}\approx z_{1j}\approx\Theta_{\zeta}$ and $r_{ij}^2=z_{ij}w_{ij}$, this implies $w_{2j}\approx w_{1j}$, hence
\[
|w_{12}|=|w_{1j}-w_{2j}|
\le |w_{1j}|+|w_{2j}|\preceq w_{1j}.
\]
On the other hand, $Z_{ij}=z_{ij}/r_{ij}^3$ and $r_{ij}^2=z_{ij}w_{ij}$ give the exact squared-ratio identity
\[
\left(\frac{Z_{12}}{Z_{1j}}\right)^2
=\frac{z_{1j}}{z_{12}}
 \left(\frac{w_{1j}}{w_{12}}\right)^3
\succeq\frac{\Theta_{\zeta}}{z_{12}}
\longrightarrow\infty,
\]
since $z_{12}=o(\Theta_{\zeta})$. This contradicts the finite limit $Z_{12}/Z_{1j}\to-m_j/m_2$ established above.

If $j\ne k$, the four remaining cross distances diverge, so their $W$ terms are $o(\epsilon)$ by \eqref{eq:6.5}. The $w$-difference equation is therefore
\[
w_{12}-(m_1+m_2)W_{12}=m_jW_{1j}-m_kW_{2k}+o(\epsilon).
\]
The collapsing attachments and $z_{1j}\approx z_{2k}\approx\Theta_{\zeta}$ give $w_{1j},w_{2k}=o(\Theta_{\zeta}^{-1})=o(\epsilon)$. By Lemma~\ref{lem:5.1}, $w_{jk}\approx\epsilon$, hence $w_{12}=w_{jk}+w_{1j}-w_{2k}\sim w_{jk}$. The preceding inverse-cube ratios imply $\delta_{1j}/\delta_{12}\to0$ and $\delta_{2k}/\delta_{12}\to0$, so
\[
\frac{W_{1j}}{W_{12}}
=\frac{\delta_{1j}}{\delta_{12}}\frac{w_{1j}}{w_{12}}
\longrightarrow0,\qquad
\frac{W_{2k}}{W_{12}}
=\frac{\delta_{2k}}{\delta_{12}}\frac{w_{2k}}{w_{12}}
\longrightarrow0.
\]
Moreover, $\delta_{12}\to\infty$ and $W_{12}=\delta_{12}w_{12}\succ\epsilon$, so the remainder divided by $W_{12}$ also tends to zero. Dividing the difference equation by $W_{12}$ now gives $1/\delta_{12}-(m_1+m_2)\to-(m_1+m_2)\ne0$ on the left and zero on the right, a contradiction.

If all attachments in Case~\ref{case:6.B} have finite nonzero limits, all cross coefficients are bounded. Equation~\eqref{eq:6.6} gives $U_{S_3}=O(\epsilon^2)\to0$ independently of the binary scale. A single finite attachment forces $r_{12}\to0$ by the equation for the unattached exterior body; two finite attachments at the same vertex also force $r_{12}\to0$ by geometry. The term $m_1m_2/r_{12}$ then makes $U$ unbounded, since all other external potentials are bounded. For two finite attachments at distinct vertices the $w$-difference and individual $z$ equations give, as $\epsilon\to0$,
\[
\delta_{12}\to1/(m_1+m_2),\qquad m_j\delta_{1j},m_k\delta_{2k}\to e_1(S_3)/M.
\]
Thus
\[
\begin{gathered}
M^{1/3}U\le m_1m_2(M/(m_1+m_2))^{1/3}+e_1(S_3)^{1/3}(m_1m_j^{2/3}+m_2m_k^{2/3})\\
\le m_1m_2+(e_1(S_3)/12+e_1(S_3)^{1/3})(m_1+m_2)<e_2(S_5).
\end{gathered}
\]
Here $m_j,m_k\le1$, \eqref{eq:6.4} supplies the first correction, and $e_1(S_3)/12+e_1(S_3)^{1/3}<e_1(S_3)$ for $e_1(S_3)=9/4$. This completes Case~\ref{case:6.B} without a restriction on the ratio $m_1/m_2$.

\subsubsection{The case when \texorpdfstring{$\epsilon\prec\zeta_{S_3,1}\prec\zeta_{S_3,2}$}{epsilon smaller than zeta(S3,1), smaller than zeta(S3,2)}}\label{sec:6.5}
\label{subsec:6.5}

In Case~\ref{case:6.C}, after relabeling the exterior bodies, $\epsilon\prec \zeta_{S_3,1}\prec\Theta_{\zeta}=|\zeta_{S_3,2}|$. Body 2 has the larger relative $z$ scale and at most one attachment. We first treat its attachment to vertex 3 or 4 using the first identity in Lemma~\ref{lem:5.2} for the cluster $T_4$, with arbitrary positive exterior masses.

If $r_{23}$ is bounded, the identity for the moment of inertia in \eqref{eq:6.3} gives
\[
\omega_{S_3,2}\sim \frac12w_{34},\qquad \omega_{S_3,1}\sim(1+e_1(S_3)/m_1)\frac12w_{34}.
\]
The distances $r_{12},r_{24},r_{25}$ diverge, and the $z$ equation for body 2 gives $\delta_{23}\to(m_1+e_1(S_3))/M$. The $w$-separations from body 1 to vertices 3, 4, and 5 are asymptotic to $(e_1(S_3)/m_1)\frac12w_{34}$, $(2+e_1(S_3)/m_1)\frac12w_{34}$, $(1+e_1(S_3)/m_1)\frac12w_{34}$, respectively. The three $z$-separations from body 1 to $S_3$ are asymptotic to $\zeta_{S_3,1}$, so $r_{1j}^2\approx\epsilon\zeta_{S_3,1}$ for each $j\in S_3$. If $\epsilon\zeta_{S_3,1}$ were bounded away from zero on a subsequence, the coefficients $\delta_{1j}$, $j\in S_3$, would all be bounded, giving $Z_{1j}=O(\zeta_{S_3,1})$ for $j\in S_3$. The divergent binary distance also gives $Z_{12}=o(\Theta_{\zeta})$. The body-1 equation would then imply
\[
z_1=m_2Z_{12}+\sum\limits_{j\in S_3}m_jZ_{1j}
=o(\Theta_{\zeta}),
\]
contradicting $z_1\sim-(m_2/M)\zeta_{S_3,2}$ and $|\zeta_{S_3,2}|=\Theta_{\zeta}$. Hence
\begin{equation*}
\epsilon\zeta_{S_3,1}\longrightarrow0.
\end{equation*}
Apply the first identity of Lemma~\ref{lem:5.2} to $T_4$ with $k=3$. For $i\in T_4\setminus\{3\}$, the only external body is 2, and $r_{i2}$ diverges. Since $w_i,w_{i2}=O(\epsilon)$, we have $W_{i2}=o(\epsilon)$ and
\[
f_{i,w}^{T_4}=w_i-m_2W_{i2}=O(\epsilon),
\qquad i\in T_4\setminus\{3\}.
\]
The $z$-separations satisfy $z_{13}=O(\zeta_{S_3,1})$ and $z_{43},z_{53}=O(\epsilon)$. The term with $i=3$ vanishes exactly because $z_{33}=0$, irrespective of $f_{3,w}^{T_4}$. Thus
\[
U_{T_4}=\sum\limits_{i\in T_4\setminus\{3\}}m_iz_{i3}f_{i,w}^{T_4}
=O(\epsilon\zeta_{S_3,1})+O(\epsilon^2)\longrightarrow0.
\]
The only remaining nonvanishing potential is $m_2/r_{23}$, and
\[
M^{1/3}U=m_2(m_1+e_1(S_3))^{1/3}<m_2(m_1+e_1(S_3))<e_2(S_5),
\]
since $m_1+e_1(S_3)>1$. This contradicts \eqref{eq:3.1}. Vertex 4 is identical after relabeling. In particular no exceptional signed-force cancellation or exterior mass ratio has been excluded from this argument.

Now suppose body 2 is attached to vertex 5. Constancy of the moment of inertia gives $\omega_{S_3,1},\omega_{S_3,2}=o(\epsilon)$. Then $z_{15}\approx \zeta_{S_3,1}$, $z_{25}\approx \Theta_{\zeta}$, and $w_{15}=o(\epsilon)$. The exact center is
\[
z_5=-\frac{m_1z_{15}+m_2z_{25}-2\zeta_{T_2,5}}{M}.
\]
The exact squared distances satisfy
\begin{equation*}
r_{15}^2=z_{15}w_{15},\qquad r_{25}^2=z_{25}w_{25},\qquad r_{12}^2=z_{12}w_{12}.
\end{equation*}
If $\delta_{12}\to\infty$, the $z$ equation for body 2 gives $\delta_{25}/\delta_{12}\to-4m_1$. Thus $r_{25}\approx r_{12}\to0$, and $w_{15}=w_{12}+w_{25}$ gives $w_{15}\preceq w_{12}\approx r_{12}^2/\Theta_{\zeta}$. If $\delta_{12}$ has a finite nonzero limit, the finite nonzero limit of $r_{12}$ and boundedness of $r_{25}$ give the same bound. In both cases $r_{12}=O(1)$, $z_{12}\approx\Theta_{\zeta}$, and $z_{15}\approx z_{1j}\approx\zeta_{S_3,1}=o(\Theta_{\zeta})$ for $j\in T_2$. Since $w_{15}=o(\epsilon)$ and $w_{1j}\approx\epsilon$, the squared-ratio identities give
\[
\begin{aligned}
\left(\frac{Z_{12}}{Z_{15}}\right)^2
&=\frac{z_{15}}{z_{12}}
  \left(\frac{w_{15}}{w_{12}}\right)^3\longrightarrow0,\\
\left(\frac{Z_{1j}}{Z_{15}}\right)^2
&=\frac{z_{15}}{z_{1j}}
  \left(\frac{w_{15}}{w_{1j}}\right)^3
=O\!\left(\frac{w_{15}^3}{\epsilon^3}\right)\longrightarrow0,
\qquad j\in T_2.
\end{aligned}
\]
The position term is negligible as well, because
\[
\frac{Z_{15}^2}{\Theta_{\zeta}^2}
=\frac{1}{\Theta_{\zeta}^2z_{15}\,w_{15}^3}
\succeq\frac{\Theta_{\zeta}}{\zeta_{S_3,1}\,r_{12}^6}
\longrightarrow\infty.
\]
Since $z_1=O(\Theta_{\zeta})$, dividing the exact equation
\[
z_1=m_2Z_{12}+Z_{13}+Z_{14}+\tfrac14Z_{15}
\]
by $Z_{15}$ makes the left side and the first three terms on the right tend to zero, leaving the contradiction $0=1/4$.

If $\delta_{12}\to0$, the force equations give
\[
\delta_{25}\longrightarrow\frac{4(m_1+e_1(S_3))}{M},\qquad
\frac{Z_{15}}{z_{25}}\longrightarrow-\frac{4m_2}{M}.
\]
Consequently,
\[
Z_{15}\approx\Theta_{\zeta},\qquad
\delta_{15}=\frac{Z_{15}}{z_{15}},\qquad
\delta_{15}\approx\frac{\Theta_{\zeta}}{\zeta_{S_3,1}}.
\]
The exact identity $Z_{15}^2z_{15}w_{15}^3=1$ now yields
\[
w_{15}^3\,\zeta_{S_3,1}\,\Theta_{\zeta}^2\approx1,
\qquad
\epsilon^3\zeta_{S_3,1}\Theta_{\zeta}^2\longrightarrow\infty,
\]
where the second conclusion uses $w_{15}=o(\epsilon)$. The identities for $f_5^{S_3}$ from Section~\ref{subsec:6.2} give
\[
f_{5,z}^{S_3}=z_5+m_1Z_{15}+m_2Z_{25},\qquad
\frac{f_{5,z}^{S_3}}{z_{25}}\longrightarrow
\frac{m_2(4e_1(S_3)-1)}{M}=\frac{8m_2}{M}.
\]
For the transverse component, we have
\[
\begin{aligned}
\frac{W_{15}}{\Theta_{\zeta}}
&=\frac{Z_{15}}{\Theta_{\zeta}}\frac{w_{15}}{z_{15}}
=O\!\left(\frac{w_{15}}{\zeta_{S_3,1}}\right)\longrightarrow0,\\
\frac{W_{25}}{\Theta_{\zeta}}&=O(\Theta_{\zeta}^{-2}),\qquad
\frac{w_5}{\Theta_{\zeta}}=O(\epsilon/\Theta_{\zeta}).
\end{aligned}
\]
Thus $f_{5,w}^{S_3}=w_5+m_1W_{15}+m_2W_{25}=o(\Theta_{\zeta})$. Together with $f_{5,z}^{S_3}\approx\Theta_{\zeta}$, Lemma~\ref{lem:5.3}(1) implies $\mathcal{K}\approx\epsilon\Theta_{\zeta}^2$.
For $j\in T_2$, the bounds $z_{1j}\approx z_{15}$, $w_{1j}\approx\epsilon$, and $w_{15}=o(\epsilon)$ give $Z_{1j}/Z_{15}\to0$ by the same squared-ratio identity. Hence
\[
\frac{\delta_{1j}}{\delta_{15}}
=\frac{Z_{1j}}{Z_{15}}\frac{z_{15}}{z_{1j}}\longrightarrow0,
\qquad j\in T_2.
\]
Every exterior coefficient is $O(\delta_{15})$, and $\epsilon^3\delta_{15}\to0$. Lemma~\ref{lem:5.4} therefore gives
\[
\mathcal{S}_{11}=\tfrac{m_1}4\delta_{15}(1+o(1)),\qquad \mathcal{S}_{13}=O(\delta_{15}),\qquad \mathcal{S}_{33}=\mathcal{K}+O(\delta_{15}).
\]
The bounded attachment of body 2 ensures noncollinearity; the rank-one minor forces $\mathcal{K}=O(\Theta_{\zeta}/\zeta_{S_3,1})$. But the cubic relation above and $\epsilon^2\Theta_{\zeta}\to0$ from \eqref{eq:6.5} imply
\[
\epsilon\Theta_{\zeta}\zeta_{S_3,1}
=\frac{\epsilon^3\zeta_{S_3,1}\Theta_{\zeta}^2}{\epsilon^2\Theta_{\zeta}}
\longrightarrow\infty,
\]
contradicting that bound.

Finally, if body 2 has no attachment, its equation gives
\[
\delta_{12}\longrightarrow\frac{m_1+e_1(S_3)}{m_1M}\ne0.
\]
The distance $r_{12}$ has a finite nonzero limit, so $\omega_{S_3,1}-\omega_{S_3,2}=O(1/\Theta_{\zeta})$. In \eqref{eq:6.3} the coefficient of a common transverse level is $(e_1(S_3)/M)(m_1\zeta_{S_3,1}+m_2\zeta_{S_3,2})$, comparable to $\Theta_{\zeta}$; thus $\omega_{S_3,1},\omega_{S_3,2}=O(1/\Theta_{\zeta})$. Lemma~\ref{lem:5.1} then gives $w_{15}=o(\epsilon)$. Since $r_{25}\to\infty$, $z_{25}\approx\Theta_{\zeta}$, and $w_{15}-w_{25}=w_{12}=O(1/\Theta_{\zeta})$, we also have $\Theta_{\zeta}w_{15}\to\infty$. For $j\in T_2$, the bounds $z_{1j}\approx z_{15}$, $w_{1j}\approx\epsilon$, and $w_{15}=o(\epsilon)$ imply $Z_{1j}/Z_{15}\to0$ by the squared-ratio identity in Section~\ref{sec:2.2}. Retaining the binary term in the exact equation
\[
z_1-m_2Z_{12}=Z_{13}+Z_{14}+\tfrac14Z_{15}
\]
and using $z_{15}\sim\zeta_{S_3,1}$ therefore yields
\[
\frac{Z_{15}}{4\zeta_{S_3,2}}\longrightarrow\frac{m_2e_1(S_3)}{m_1M}\ne0,\qquad
\delta_{15}\approx\frac{\Theta_{\zeta}}{\zeta_{S_3,1}}\longrightarrow\infty.
\]
For $j\in T_2$, we have $z_{1j}\approx z_{15}$ and $w_{1j}\approx\epsilon$. The squared-ratio identity therefore gives
\[
\left(\frac{W_{1j}}{W_{15}}\right)^2
=\left(\frac{z_{15}}{z_{1j}}\right)^3
 \frac{w_{15}}{w_{1j}}
=O\!\left(\frac{w_{15}}{\epsilon}\right)\longrightarrow0,
\qquad j\in T_2.
\]
Total centering gives $w_1=O(1/\Theta_{\zeta})$, and $W_{12}=\delta_{12}w_{12}=O(1/\Theta_{\zeta})$. The exact body-1 equation $w_1=m_2W_{12}+W_{13}+W_{14}+\tfrac14W_{15}$ consequently implies
\[
\bigl(\tfrac14+o(1)\bigr)W_{15}
=w_1-m_2W_{12}=O(1/\Theta_{\zeta}),
\]
so $W_{15}=O(1/\Theta_{\zeta})$. On the other hand,
\[
\Theta_{\zeta}W_{15}
=\delta_{15}\,\Theta_{\zeta}w_{15}\longrightarrow\infty,
\]
because both $\delta_{15}$ and $\Theta_{\zeta}w_{15}$ tend to infinity. This contradicts the preceding bound.

These alternatives exhaust the possible attachments of body 2 and complete Case~\ref{case:6.C}. Together with Cases~\ref{case:6.A} and~\ref{case:6.B}, proved in Sections~\ref{subsec:6.2}--\ref{subsec:6.4}, they prove Theorem~\ref{thm:7.1} under the assumption $m_3=m_4$.

\subsection{The case \texorpdfstring{$m_3\ne m_4$}{m3 unequal to m4}}\label{sec:diagram5:unequal}

Assume now $m_3\ne m_4$. The coefficients $\eta_j$, $j\in S_3$, in \eqref{eq:7.3} are then nonzero, as shown below. This removes the finite three-point limit that occurred in the equal-endpoint case; the remaining transverse alignment is controlled by the first nonzero term in the expansion of $\mathcal{K}$.

\subsubsection{Mass normalization and coordinates relative to \texorpdfstring{$C_{S_3}$}{C(S3)}}\label{sec:7.1}
\label{subsec:7.1}

\begin{lemma}\label{lem:7.2}
In an AK diagram-5 limit with positive masses, relabeling the vertices of the isolated $zw$-triangle and applying a common positive mass scaling gives masses satisfying
\begin{equation}
\begin{gathered}
m_3,m_4,m_5>0,\qquad
\frac1{\sqrt{m_5}}=\frac1{\sqrt{m_3}}+\frac1{\sqrt{m_4}},\\
\sqrt{m_3}+\sqrt{m_4}=1,\qquad m_5=m_3m_4.
\end{gathered}
\label{eq:7.1}
\end{equation}
Its total mass is $e_1(S_3)$. Define the exact displacement $\mu$ by
\begin{equation}
q_5=q_3+\sqrt{m_3}(q_4-q_3)+\mu.
\label{eq:7.2}
\end{equation}
Then $r_{34}^2=z_{34}w_{34}$, $z_{34},w_{34},r_{34}\approx\epsilon$, $\mu_z,\mu_w=o(\epsilon)$, and the signed branches satisfy $r_{35}\sim-\sqrt{m_3}\,r_{34}$ and $r_{45}\sim-\sqrt{m_4}\,r_{34}$. Writing $C_{S_3}$ for the center of mass of the isolated $zw$-triangle, its leading centered coordinates are
\begin{equation}
\begin{aligned}
q_{S_3,j}&=\eta_j(q_4-q_3)+o(\epsilon),\qquad j\in S_3,\\
\eta_3&=-\frac{m_4(1+\sqrt{m_3})}{1-\sqrt{m_3m_4}},\qquad
\eta_4=\frac{m_3(1+\sqrt{m_4})}{1-\sqrt{m_3m_4}},\\
\eta_5&=\frac{\sqrt{m_3}-\sqrt{m_4}}{1-\sqrt{m_3m_4}}.
\end{aligned}
\label{eq:7.3}
\end{equation}
In particular $\eta_3<\eta_5<\eta_4$, $\eta_5-\eta_3=\sqrt{m_3}$, and $\eta_4-\eta_5=\sqrt{m_4}$. If $m_3\ne m_4$, none of the three $\eta_j$ is zero.
\end{lemma}

\begin{proof}
For positive triangle masses, the necessary relation from \cite[Section~5.3]{AK} is $m_5^{-1/2}=m_3^{-1/2}+m_4^{-1/2}$ after placing the smallest mass at vertex 5. A common positive mass scaling makes $\sqrt{m_3}+\sqrt{m_4}=1$; retaining the same symbols for the scaled masses, this relation then gives $m_5=m_3m_4$.

For these general masses, the equations expressing cancellation of the leading internal force terms again make the coordinate differences between vertices in $S_3$ proportional in $z,w$. We derive the coefficients $\eta_j$ from the endpoint force equations, retaining every distance sign. These equations imply
\[
\begin{aligned}
\lim\limits_{\epsilon\to0}\frac{z_5-z_3}{z_4-z_3}
&\in\{\pm\sqrt{m_3},\pm \mathrm{i}\sqrt{m_3}\},\\
\lim\limits_{\epsilon\to0}\frac{z_4-z_5}{z_4-z_3}
&\in\{\pm\sqrt{m_4},\pm \mathrm{i}\sqrt{m_4}\}.
\end{aligned}
\]
The two ratios sum to one, and the real parts of their limits are at most $\sqrt{m_3}$ and $\sqrt{m_4}$, respectively. Since $\sqrt{m_3}+\sqrt{m_4}=1$, both bounds are attained: the limits are $\sqrt{m_3}$ and $\sqrt{m_4}$, respectively. The same limits hold with $z$ replaced by $w$. The endpoint equations then require the short signed distances to have the opposite sign to their common long-distance branch, giving $r_{35}\sim-\sqrt{m_3}\,r_{34}$ and $r_{45}\sim-\sqrt{m_4}\,r_{34}$. The AK $zw$-edge estimates give the stated coordinate sizes and $\mu=o(\epsilon)$.

Substituting \eqref{eq:7.2} into the definition of $C_{S_3}$ and using \eqref{eq:7.1} gives the centered-coordinate formulas \eqref{eq:7.3}. The ordering and differences follow from the leading positions $0,\sqrt{m_3},1$. The coefficients satisfy $\eta_3<0<\eta_4$, and $\eta_5=0$ if and only if $m_3=m_4$.
\end{proof}

For the rest of Section~\ref{sec:diagram5:unequal}, retain these normalized masses and the same notation for the correspondingly scaled coordinates. The exterior masses remain arbitrary and positive. The common identities, scale bounds, and Cases~\ref{case:6.A}--\ref{case:6.C} from Section~\ref{sec:6.1} apply unchanged.

\subsubsection{The case when \texorpdfstring{$\zeta_{S_3,1}\approx\zeta_{S_3,2}$}{zeta(S3,1) comparable to zeta(S3,2)}}\label{sec:7.3}
\label{subsec:7.3}

We treat Cases~\ref{case:6.A} and~\ref{case:6.B}, beginning with the latter. First suppose $\zeta_{S_3,1},\zeta_{S_3,2}\approx\Theta_{\zeta}$ and $\zeta_{S_3,1}-\zeta_{S_3,2}=o(\Theta_{\zeta})$. Then $z_1\sim z_2$, $C_{S_3,z}\sim-(m_1+m_2)z_1/e_1(S_3)$, and every cross $z$-separation is asymptotic to $Mz_1/e_1(S_3)$. The mass-weighted equations for the exterior bodies cancel the binary interaction exactly:
\[
m_1z_1+m_2z_2
=\sum\limits_{\substack{i\in S_2\\l\in S_3}}m_im_lZ_{il}.
\]
At least one attachment exists. A lone collapsing attachment, or a collapsing attachment paired with a finite one, would be the unique dominant term and is also impossible. If two attachments $1$--$j$, $2$--$k$ with $j,k\in S_3$ collapse, then $Z_{1j}/\Theta_{\zeta}\to\infty$, while every other cross term is $O(\Theta_{\zeta})$. Dividing the weighted sum and the body-1 equation by $Z_{1j}$ gives
\[
\frac{Z_{2k}}{Z_{1j}}\longrightarrow-\frac{m_1m_j}{m_2m_k},\qquad
\frac{Z_{12}}{Z_{1j}}\longrightarrow-\frac{m_j}{m_2}.
\]
These are the balances used in Section~\ref{subsec:6.4}; they require only positive masses. For $j=k$, its squared-ratio argument applies: $w_{12}\preceq w_{1j}$ and
\[
\left(\frac{Z_{12}}{Z_{1j}}\right)^2
=\frac{z_{1j}}{z_{12}}
 \left(\frac{w_{1j}}{w_{12}}\right)^3
\succeq\frac{\Theta_{\zeta}}{z_{12}}\longrightarrow\infty,
\]
contradicting the finite limit above. For $j\ne k$, the $w$-separations satisfy $w_{1j},w_{2k}=o(\epsilon)$, so $w_{12}\sim w_{jk}$ and $w_{12}\approx\epsilon$ by Lemma~\ref{lem:7.2}. As in Section~\ref{subsec:6.4}, the ratios of the $z$-separations give $\delta_{1j}/\delta_{12},\delta_{2k}/\delta_{12}\to0$ and $\delta_{12}\to\infty$. The difference equation is
\[
w_{12}-(m_1+m_2)W_{12}
=m_jW_{1j}-m_kW_{2k}+o(\epsilon).
\]
Here
\[
\frac{W_{1j}}{W_{12}}
=\frac{\delta_{1j}}{\delta_{12}}\frac{w_{1j}}{w_{12}}\longrightarrow0,\qquad
\frac{W_{2k}}{W_{12}}
=\frac{\delta_{2k}}{\delta_{12}}\frac{w_{2k}}{w_{12}}\longrightarrow0.
\]
Moreover, $W_{12}\succ\epsilon$, so the remainder divided by $W_{12}$ tends to zero. Division by $W_{12}$ gives $1/\delta_{12}-(m_1+m_2)\to-(m_1+m_2)\ne0$ on the left and zero on the right. This excludes all collapsing attachments in Case~\ref{case:6.B}.

With only finite nonzero attachments, all cross inverse cubes are bounded. Equation~\eqref{eq:6.5} gives $W_{ki}=w_{ki}\delta_{ki}=O(\epsilon)$ for $k\in S_3$, $i\in S_2$, hence
\[
f_{k,w}^{S_3}=w_k-\sum\limits_{i\in S_2}m_iW_{ki}=O(\epsilon),
\qquad k\in S_3.
\]
Lemma~\ref{lem:5.2} with reference body 5 therefore gives
\[
U_{S_3}=\sum\limits_{k\in T_2}m_kz_{k5}f_{k,w}^{S_3}
=O(\epsilon^2)\longrightarrow0,
\]
since $z_{k5}=O(\epsilon)$ for $k\in T_2$. A single finite attachment forces a collapsing binary by the equation for the unattached exterior body. Two finite attachments at the same vertex force the same conclusion by their $O(1/\Theta_{\zeta})$ $w$-separations and $z_{12}=o(\Theta_{\zeta})$. In either case $m_1m_2/r_{12}$ is the only unbounded exterior potential, contradicting finite constant $U$. For two finite attachments at distinct vertices $j,k\in S_3$, the force equations give, as $\epsilon\to0$,
\[
\delta_{12}\to1/(m_1+m_2),\qquad m_j\delta_{1j},m_k\delta_{2k}\to e_1(S_3)/M.
\]
Therefore
\begin{equation}
\begin{aligned}
M^{1/3}U&\le m_1m_2(M/(m_1+m_2))^{1/3}\\
&\quad +e_1(S_3)^{1/3}(m_1m_j^{2/3}+m_2m_k^{2/3}).
\end{aligned}
\label{eq:7.9}
\end{equation}
To compare this with $e_2(S_5)$, concavity and $m_1m_2/(m_1+m_2)\le\min(m_1,m_2)$ give
\[
\begin{aligned}
m_1m_2\{(M/(m_1+m_2))^{1/3}-1\}
&\le\frac{e_1(S_3)}{3}\frac{m_1m_2}{m_1+m_2}\\
&\le\frac{e_1(S_3)}{3}\min(m_1,m_2).
\end{aligned}
\]
For each $i\in S_3$, the difference $e_1(S_3)-e_1(S_3)^{1/3}m_i^{2/3}$ is positive. Since $j\ne k$ and $m_j+m_k<e_1(S_3)$, concavity also gives
\[
m_j^{2/3}+m_k^{2/3}
\le 2^{1/3}(m_j+m_k)^{2/3}
<2^{1/3}e_1(S_3)^{2/3}.
\]
Consequently,
\[
\begin{aligned}
&(m_1+m_2)e_1(S_3)-e_1(S_3)^{1/3}(m_1m_j^{2/3}+m_2m_k^{2/3})\\
&\ge \min(m_1,m_2)\bigl[2e_1(S_3)-e_1(S_3)^{1/3}(m_j^{2/3}+m_k^{2/3})\bigr]\\
&>(2-2^{1/3})e_1(S_3)\min(m_1,m_2)\\
&>\frac{e_1(S_3)}{3}\min(m_1,m_2).
\end{aligned}
\]
This exceeds the binary correction above. Equation \eqref{eq:7.9} therefore implies $M^{1/3}U<m_1m_2+e_1(S_3)(m_1+m_2)<e_2(S_5)$, contradicting \eqref{eq:3.1}.

Now suppose $\zeta_{S_3,1},\zeta_{S_3,2},\zeta_{S_3,1}-\zeta_{S_3,2}\approx \Theta_{\zeta}$. If an attachment collapses, there are four possible patterns. One total attachment makes the unattached equation bound $\delta_{12}$, then the attached equation bounds the collapsing coefficient. Two distinct attachment vertices give $w_{12}\approx\epsilon$ and $r_{12}\to\infty$, so both individual equations bound their attachment coefficients. A shared vertex with one finite attachment first bounds $\delta_{12}$ using the equation for that exterior body, then bounds the other attachment coefficient.

The only remaining pattern is two collapsing attachments to a shared vertex $j\in S_3$. As in Section~\ref{subsec:6.2}, we apply the two identities in Lemma~\ref{lem:5.2} to the clusters $S_3$ and $S_2\cup\{j\}$, respectively. The four other cross distances diverge, and the individual equations give
\[
\begin{aligned}
f_{k,w}^{S_3}&=w_k-\sum\limits_{i\in S_2}m_iW_{ki}=O(\epsilon),
&&k\in S_3\setminus\{j\},\\
f_{i,z}^{S_2\cup\{j\}}&=z_i-\sum\limits_{l\in S_3\setminus\{j\}}m_lZ_{il}=O(\Theta_{\zeta}),
&&i\in S_2.
\end{aligned}
\]
Since $z_{kj}=O(\epsilon)$ for $k\in S_3\setminus\{j\}$ and $w_{ij}=o(\Theta_{\zeta}^{-1})$ for $i\in S_2$, applying Lemma~\ref{lem:5.2} to both clusters with reference body $j$ yields
\[
\begin{aligned}
U_{S_3}&=\sum\limits_{k\in S_3\setminus\{j\}}m_kz_{kj}f_{k,w}^{S_3}=O(\epsilon^2),\\
U_{S_2\cup\{j\}}&=\sum\limits_{i\in S_2}m_iw_{ij}f_{i,z}^{S_2\cup\{j\}}=o(1).
\end{aligned}
\]
The terms corresponding to vertex $j$ vanish exactly, and all four other cross potentials vanish. Hence $U\to0$, which is impossible.

If no attachment collapses, all cross coefficients are bounded. The $z$-separation of vertices 1 and 2 bounds the binary coefficient through either $z$ equation for the exterior bodies. Thus all seven exterior or cross distances are bounded away from zero, and their $z$-separations have modulus comparable to $\Theta_{\zeta}$. The identity $W_{ij}=1/(r_{ij}z_{ij})$ gives, as in Section~\ref{subsec:6.2},
\[
W_{12},W_{ij}=O(\Theta_{\zeta}^{-1})=o(\epsilon),
\qquad i\in S_2,\quad j\in S_3.
\]
The individual equations for the exterior bodies and total centering then give
\[
w_1,w_2,C_{S_3,w}=O(\Theta_{\zeta}^{-1}),
\]
so $\omega_{S_3,1},\omega_{S_3,2}=o(\epsilon)$. But every centered $w$ coordinate of a vertex $j\in S_3$ is $\eta_jw_{43}+o(\epsilon)$ with $\eta_j\ne0$, by the inequality $m_3\ne m_4$. No attachment can have its required $O(1/\Theta_{\zeta})$ $w$-separation. This contradicts the bounded attachment established before \eqref{eq:6.5} and completes all comparable-scale cases.

\subsubsection{The case when \texorpdfstring{$\epsilon\prec\zeta_{S_3,1}\prec\zeta_{S_3,2}$}{epsilon smaller than zeta(S3,1), smaller than zeta(S3,2)} and \texorpdfstring{$r_{2j}\to\infty$ for all $j\in S_3$}{r2j tends to infinity for all j in S3}}\label{sec:7.7}
\label{subsec:7.7}

Throughout the remaining subsections we are in Case~\ref{case:6.C}. Relabel the exterior bodies so
\begin{equation}
\epsilon\prec \zeta_{S_3,1}\prec \Theta_{\zeta}=|\zeta_{S_3,2}|.
\label{eq:7.10}
\end{equation}
Body 2 has the larger relative $z$ scale and at most one attachment. First suppose that it has none, so all three of its distances to the triangle diverge. Its $z$ equation gives
\begin{equation}
\delta_{12}\longrightarrow\frac{e_1(T_4)}{m_1M}\ne0.
\label{eq:7.23}
\end{equation}
The finite nonzero limit of $r_{12}$ implies $\omega_{S_3,1}-\omega_{S_3,2}=O(1/\Theta_{\zeta})$. In \eqref{eq:6.3}, the coefficient of a common $\omega$ level is $(e_1(S_3)/M)(m_1\zeta_{S_3,1}+m_2\zeta_{S_3,2})$, of modulus comparable to $\Theta_{\zeta}$. The remaining difference term is bounded, so $\omega_{S_3,1},\omega_{S_3,2}=O(1/\Theta_{\zeta})$. Since all three coefficients $\eta_j$, $j\in S_3$, are nonzero, each $w$-separation from body 1 to $S_3$ has modulus comparable to $\epsilon$; therefore $r_{1j}^2\approx\zeta_{S_3,1}\epsilon$ for every $j\in S_3$.

Retaining the binary term in the exact body-1 equation and using \eqref{eq:6.2} and \eqref{eq:7.23} also gives
\begin{equation*}
z_1-m_2Z_{12}
=\sum\limits_{j\in S_3}m_jZ_{1j}
\sim\frac{m_2e_1(S_3)}{m_1M}\zeta_{S_3,2},
\end{equation*}
We estimate the cluster potential directly by applying Lemma~\ref{lem:5.2} to $T_4$ with any fixed reference body $j\in S_3$. Total centering gives $w_1=O(\Theta_{\zeta}^{-1})$, and the finite nonzero limit of $r_{12}$ gives $W_{12}=O(\Theta_{\zeta}^{-1})$. For $i\in S_3\setminus\{j\}$, the divergent distance $r_{i2}$ and \eqref{eq:6.5} give $W_{i2}=o(\epsilon)$. Therefore
\[
\begin{aligned}
f_{1,w}^{T_4}&=w_1-m_2W_{12}=O(\Theta_{\zeta}^{-1}),\\
f_{i,w}^{T_4}&=w_i-m_2W_{i2}=O(\epsilon),
\qquad i\in S_3\setminus\{j\}.
\end{aligned}
\]
Since $z_{1j}=O(\zeta_{S_3,1})$ and $z_{ij}=O(\epsilon)$ for $i\in S_3\setminus\{j\}$, retaining the stronger bound for $f_{1,w}^{T_4}$ gives
\[
\begin{aligned}
U_{T_4}
&=m_1z_{1j}f_{1,w}^{T_4}
+\sum\limits_{i\in S_3\setminus\{j\}}m_iz_{ij}f_{i,w}^{T_4}\\
&=O\!\left(\frac{\zeta_{S_3,1}}{\Theta_{\zeta}}\right)
 +O(\epsilon^2)\longrightarrow0,
\end{aligned}
\]
using the already established $\zeta_{S_3,1}=o(\Theta_{\zeta})$.
The term with $i=j$ vanishes because $z_{jj}=0$. All potential terms between body 2 and $S_3$ vanish, leaving only $m_1m_2/r_{12}$. By \eqref{eq:7.23},
\[
M^{1/3}U=m_1m_2\left(\frac{e_1(T_4)}{m_1}\right)^{1/3}=m_2m_1^{2/3}e_1(T_4)^{1/3}<m_2e_1(T_4)<e_2(S_5).
\]
This contradicts the physical lower bound \eqref{eq:3.1}.

\subsubsection{The case when \texorpdfstring{$\epsilon\prec\zeta_{S_3,1}\prec\zeta_{S_3,2}$}{epsilon smaller than zeta(S3,1), smaller than zeta(S3,2)} and \texorpdfstring{$r_{2j}=O(1)$ for some $j\in S_3$}{r2j is bounded for some j in S3}}\label{sec:7.4}
\label{subsec:7.4}

We now treat the complementary branch, in which body 2 has a bounded attachment. The scale assumptions remain \eqref{eq:7.10}.

For body 1, $z_{1j}\sim\zeta_{S_3,1}$ for every $j\in S_3$. Suppose body 2 has a bounded attachment to vertex $j\in S_3$. Since $z_{2j}\sim\zeta_{S_3,2}$, we have $w_{2j}=r_{2j}^2/z_{2j}=O(1/\Theta_{\zeta})=o(\epsilon)$. Lemma~\ref{lem:7.2} then gives $\omega_{S_3,2}=\eta_jw_{43}+o(\epsilon)$, and \eqref{eq:6.3} gives
\[
\omega_{S_3,1}=\frac{e_1(T_4)}{m_1}\eta_jw_{43}+o(\epsilon).
\]
Since $\eta_j\ne0$, the binary $w$-separation has nonzero leading coefficient $(e_1(S_3)/m_1)\eta_j$. Thus $r_{12}$ and the distances $r_{2k}$, $k\in S_3\setminus\{j\}$, diverge. The $z$ equation for body 2 gives
\begin{equation}
m_j\delta_{2j}\longrightarrow e_1(T_4)/M.
\label{eq:7.11}
\end{equation}
In particular this attachment is finite and nonzero.

If $\frac{e_1(T_4)}{m_1}\eta_j$ differs from all levels $\eta_k$ with $k\in S_3$, then $r_{1k}^2\approx\zeta_{S_3,1}\epsilon$ for every $k\in S_3$. If $\zeta_{S_3,1}\epsilon$ were bounded away from zero on a subsequence, all inverse cubes $\delta_{1k}$, $k\in S_3$, would be bounded, giving $Z_{1k}=O(\zeta_{S_3,1})$ for $k\in S_3$. Since the binary distance diverges, $Z_{12}=o(\Theta_{\zeta})$, and the body-1 equation would give
\[
z_1=m_2Z_{12}+\sum\limits_{k\in S_3}m_kZ_{1k}=o(\Theta_{\zeta}),
\]
contradicting $z_1\sim-(m_2/M)\zeta_{S_3,2}$. Hence $\zeta_{S_3,1}\epsilon\to0$.

Apply the first identity of Lemma~\ref{lem:5.2} to $T_4$ with reference body $j$. For each $i\in T_4\setminus\{j\}$, the only external body is 2, and $r_{i2}$ diverges. Equation~\eqref{eq:6.5} gives $w_i,w_{i2}=O(\epsilon)$, hence $W_{i2}=o(\epsilon)$ and
\[
f_{i,w}^{T_4}=w_i-m_2W_{i2}=O(\epsilon),
\qquad i\in T_4\setminus\{j\}.
\]
The $z$-separations satisfy $z_{1j}=O(\zeta_{S_3,1})$ and $z_{ij}=O(\epsilon)$ for $i\in S_3\setminus\{j\}$. The term with $i=j$ is zero because $z_{jj}=0$, so
\[
U_{T_4}=\sum\limits_{i\in T_4\setminus\{j\}}m_iz_{ij}f_{i,w}^{T_4}
=O(\zeta_{S_3,1}\epsilon)+O(\epsilon^2)\longrightarrow0.
\]
The only remaining potential is $m_2m_j/r_{2j}$; by \eqref{eq:7.11},
\[
\begin{aligned}
M^{1/3}U&=m_2m_j^{2/3}e_1(T_4)^{1/3}\\
&<m_2e_1(T_4)<e_2(S_5),
\end{aligned}
\]
since $m_j<e_1(S_3)<e_1(T_4)$. This excludes every case with no transverse alignment, including all possible signed cancellations among forces on body 1 from $S_3$.

An alignment requires $\frac{e_1(T_4)}{m_1}\eta_j=\eta_k$ for distinct $j,k\in S_3$. Since $e_1(T_4)/m_1>1$, multiplying either endpoint level $\eta_3<0<\eta_4$ by this factor moves it outside the interval $[\eta_3,\eta_4]$. Thus $j=5$. If $m_3>m_4$, the middle level is positive and its only possible target is $k=4$; the case $m_3<m_4$ is its reflection. For $m_3>m_4$ the alignment condition is exactly
\begin{equation}
m_1=\frac{e_1(S_3)\eta_5}{\eta_4-\eta_5}=\frac{(1-\sqrt{m_3m_4})(\sqrt{m_3}-\sqrt{m_4})}{\sqrt{m_4}}.
\label{eq:7.12}
\end{equation}
The following estimates lead to the exclusion of this alignment in Section~\ref{subsec:7.6}.

\subsubsection{Estimates for the isolated \texorpdfstring{$zw$}{zw}-triangle}\label{sec:7.5}
\label{subsec:7.5}

Put
\[
\psi_z=\frac{\mu_z}{z_{43}},\qquad
\psi_w=\frac{\mu_w}{w_{43}}.
\]
Then $\psi_z,\psi_w\to0$ by Lemma~\ref{lem:7.2}. The exact signed distances are
\[
\begin{aligned}
r_{35}&=-\sqrt{m_3}\,r_{34}
\left[\left(1+\frac{\psi_z}{\sqrt{m_3}}\right)
\left(1+\frac{\psi_w}{\sqrt{m_3}}\right)\right]^{1/2},\\
r_{45}&=-\sqrt{m_4}\,r_{34}
\left[\left(1-\frac{\psi_z}{\sqrt{m_4}}\right)
\left(1-\frac{\psi_w}{\sqrt{m_4}}\right)\right]^{1/2},
\end{aligned}
\]
where both square roots use the analytic branch with value one at $(\psi_z,\psi_w)=(0,0)$. On these signed branches, the sums of the two dominant terms in the equations for vertex 5 are
\[
\begin{aligned}
f^{S_3}_{5,z}
&=m_3Z_{53}+m_4Z_{54}=Z_{43}\,\hat f_{5,z},\\
f^{S_3}_{5,w}
&=m_3W_{53}+m_4W_{54}=W_{43}\,\hat f_{5,w}.
\end{aligned}
\]
The two summands in each component have order $\epsilon^{-2}$, while their sum is $o(\epsilon^{-2})$ because the leading terms cancel. Here $\hat f_{5,z},\hat f_{5,w}$ are the normalized sums of these dominant terms, using $Z_{43},W_{43}$ and the masses in \eqref{eq:7.1}. They differ from those in Section~\ref{sec:diagram5:equal-identities}, where the normalizing factors are $4Z_{34},4W_{34}$.
By \eqref{eq:7.2}, $z_{53}=(\sqrt{m_3}+\psi_z)z_{43}$ and $z_{54}=(-\sqrt{m_4}+\psi_z)z_{43}$, with the analogous identities for the $w$ components. Substitution of these identities and the signed-distance formulas gives the following exact expressions for $\hat f_{5,z},\hat f_{5,w}$:
\begin{equation}
\hat f_{5,z}=-\left(1+\frac{\psi_z}{\sqrt{m_3}}\right)^{-1/2}
\left(1+\frac{\psi_w}{\sqrt{m_3}}\right)^{-3/2}
+\left(1-\frac{\psi_z}{\sqrt{m_4}}\right)^{-1/2}
\left(1-\frac{\psi_w}{\sqrt{m_4}}\right)^{-3/2},
\label{eq:7.13}
\end{equation}
\begin{equation}
\hat f_{5,w}=-\left(1+\frac{\psi_z}{\sqrt{m_3}}\right)^{-3/2}
\left(1+\frac{\psi_w}{\sqrt{m_3}}\right)^{-1/2}
+\left(1-\frac{\psi_z}{\sqrt{m_4}}\right)^{-3/2}
\left(1-\frac{\psi_w}{\sqrt{m_4}}\right)^{-1/2}.
\label{eq:7.14}
\end{equation}
All fractional powers use their analytic branches near one. Since $\sqrt{m_3}+\sqrt{m_4}=1$, expansion at $(\psi_z,\psi_w)=(0,0)$ gives
\[
\begin{aligned}
\hat f_{5,z}&=\frac{\psi_z+3\psi_w}{2\sqrt{m_3m_4}}
+O((|\psi_z|+|\psi_w|)^2),\\
\hat f_{5,w}&=\frac{3\psi_z+\psi_w}{2\sqrt{m_3m_4}}
+O((|\psi_z|+|\psi_w|)^2).
\end{aligned}
\]
The derivative with respect to $(\psi_z,\psi_w)$ at the origin is invertible. Unlike the equal-endpoint map, this map need not be odd; its analytic inverse has quadratic remainders:
\begin{equation}
\begin{aligned}
\psi_z&=\frac{\sqrt{m_3m_4}}{4}(3\hat f_{5,w}-\hat f_{5,z})+O((|\hat f_{5,z}|+|\hat f_{5,w}|)^2),\\
\psi_w&=\frac{\sqrt{m_3m_4}}{4}(3\hat f_{5,z}-\hat f_{5,w})+O((|\hat f_{5,z}|+|\hat f_{5,w}|)^2).
\end{aligned}
\label{eq:7.15}
\end{equation}
Since $Z_{43},W_{43}\approx\epsilon^{-2}$, the assumptions $f^{S_3}_{5,z}\approx\Theta_{\zeta}$ and $f^{S_3}_{5,w}=o(\Theta_{\zeta})$ give
\[
\hat f_{5,z}=\frac{f^{S_3}_{5,z}}{Z_{43}}\approx\epsilon^2\Theta_{\zeta}\longrightarrow0,\qquad
\frac{\hat f_{5,w}}{\hat f_{5,z}}
=\frac{f^{S_3}_{5,w}}{f^{S_3}_{5,z}}\frac{Z_{43}}{W_{43}}\longrightarrow0.
\]
The remainders in \eqref{eq:7.15} are therefore $o(\hat f_{5,z})$, so
\begin{equation}
\begin{aligned}
\psi_z&\sim-\frac{\sqrt{m_3m_4}}{4}\hat f_{5,z},\qquad
\frac{\psi_w}{\psi_z}\to-3,\\
\psi_z&\approx\epsilon^2\Theta_{\zeta}.
\end{aligned}
\label{eq:7.16}
\end{equation}

The scalar coefficient $\mathcal{K}$ in \eqref{eq:diagram5:schur-scalar}, evaluated at the normalized masses \eqref{eq:7.1}, is
\begin{equation}
\begin{aligned}
\mathcal{K}
&=-\mathcal{L}_{34}-\frac{\mathcal{L}_{35}\mathcal{L}_{45}}{\mathcal{L}_{35}+\mathcal{L}_{45}}\\
&=m_3m_4\delta_{34}+\frac{(m_3m_4)^2\delta_{35}\delta_{45}}{m_3\delta_{35}+m_4\delta_{45}}.
\end{aligned}
\label{eq:7.17}
\end{equation}
This is the specialization $m_5=m_3m_4$ of the common formula. Expanding the denominator with the signed-distance formulas above, its normalized form has no linear term:
\begin{equation}
\begin{aligned}
-\frac{\delta_{34}}{m_3\delta_{35}}-\frac{\delta_{34}}{m_4\delta_{45}}
&=1+\frac{3}{8\sqrt{m_3m_4}}
(\psi_z^2+\psi_w^2+6\psi_z\psi_w)\\
&\quad+O((|\psi_z|+|\psi_w|)^3).
\end{aligned}
\label{eq:7.18}
\end{equation}
Combining \eqref{eq:7.16}--\eqref{eq:7.18} yields
\begin{equation}
\mathcal{K}\sim-3\sqrt{m_3m_4}\,\delta_{34}\psi_z^2,\qquad \mathcal{K}\approx\epsilon \Theta_{\zeta}^2.
\label{eq:7.19}
\end{equation}
There is no use of real positivity in this expansion: the nonzero coefficient follows from the complex ratio $\psi_w/\psi_z\to-3$.

The pivot is $\mathcal{L}_{55}=-\mathcal{L}_{35}-\mathcal{L}_{45}\sim-\sqrt{m_3m_4}\,\delta_{34}$, and the ratios $\mathcal{L}_{35}/\mathcal{L}_{55}$ and $\mathcal{L}_{45}/\mathcal{L}_{55}$ tend to $-\sqrt{m_4}$ and $-\sqrt{m_3}$, respectively.

These limits verify the hypotheses of Lemma~\ref{lem:5.4} on $\mathcal{L}_{55}$ and the ratios $\mathcal{L}_{35}/\mathcal{L}_{55}$ and $\mathcal{L}_{45}/\mathcal{L}_{55}$. With the scale $\Theta_{\delta}$ defined in Section~\ref{sec:5.2}, its $O(\Theta_{\delta})$ remainder estimate for the Schur complement therefore applies whenever $\epsilon^3\Theta_{\delta}\to0$; the next subsection verifies this last condition for the aligned branch.

\subsubsection{The case when \texorpdfstring{$r_{25}=O(1)$ and $w_{14}=o(\epsilon)$}{r25 is bounded and w14 is of smaller order than epsilon}}\label{sec:7.6}
\label{subsec:7.6}

Assume $m_3>m_4$, body 2 is attached to vertex 5, and body 1 aligns with vertex 4 in the $w$ coordinate. Then $z_{25}\approx \Theta_{\zeta}$, $z_{14}\approx \zeta_{S_3,1}$, and $w_{14}=o(\epsilon)$, while $w_{13},w_{15}\approx\epsilon$. The distances $r_{12},r_{23},r_{24}$ diverge, and \eqref{eq:7.11} gives
\begin{equation}
\delta_{25}\longrightarrow\frac{e_1(T_4)}{Mm_5}.
\label{eq:7.20}
\end{equation}
For $j\in S_3\setminus\{4\}$, we have $z_{1j}\approx z_{14}\approx\zeta_{S_3,1}$ and $w_{1j}\approx\epsilon$. Thus $m_4Z_{14}$ dominates the terms $m_3Z_{13}$ and $m_5Z_{15}$, as the squared-ratio identity gives
\[
\left(\frac{Z_{1j}}{Z_{14}}\right)^2
=\frac{z_{14}}{z_{1j}}
 \left(\frac{w_{14}}{w_{1j}}\right)^3
=O\!\left(\frac{w_{14}^3}{\epsilon^3}\right)\longrightarrow0,
\qquad j\in S_3\setminus\{4\}.
\]
The divergent binary distance gives $Z_{12}=o(\Theta_{\zeta})$, while $z_1/z_{25}\to-m_2/M$. The exact equation $z_1=m_2Z_{12}+\sum\limits_{j\in S_3}m_jZ_{1j}$ therefore yields
\begin{equation}
\begin{gathered}
\frac{m_4Z_{14}}{z_{25}}\longrightarrow-\frac{m_2}{M},\qquad
Z_{14}\approx\Theta_{\zeta},\\
\delta_{14}\approx\frac{\Theta_{\zeta}}{\zeta_{S_3,1}},\qquad
w_{14}^3\zeta_{S_3,1}\Theta_{\zeta}^2\approx1.
\end{gathered}
\label{eq:7.21}
\end{equation}
The last relation follows from the exact identity $Z_{14}^2z_{14}w_{14}^3=1$. For the later Schur estimate, the same ratios also give
\[
\frac{\delta_{1j}}{\delta_{14}}
=\frac{Z_{1j}}{Z_{14}}\frac{z_{14}}{z_{1j}}\longrightarrow0,
\qquad j\in S_3\setminus\{4\}.
\]

The exact identity for $f_{5,z}^{S_3}$ is
\[
f_{5,z}^{S_3}=z_5+m_1Z_{15}+m_2Z_{25}.
\]
Here $Z_{15}=o(Z_{14})=o(\Theta_{\zeta})$. All positions $q_{S_3,j}$ relative to $C_{S_3}$, $j\in S_3$, are $O(\epsilon)$, so $z_5/z_{25}\to-m_2/M$. Together with the finite nonzero limit of $r_{25}$ given by \eqref{eq:7.20}, this gives
\begin{equation}
\frac{f^{S_3}_{5,z}}{z_{25}}\longrightarrow\frac{m_2(e_1(T_4)-m_5)}{Mm_5}\ne0,
\qquad f^{S_3}_{5,w}=o(\Theta_{\zeta}).
\label{eq:7.22}
\end{equation}
For the transverse estimate, the pair components satisfy
\[
\frac{W_{15}}{\Theta_{\zeta}}
=\frac{Z_{15}}{\Theta_{\zeta}}\frac{w_{15}}{z_{15}}
=O\!\left(\frac{\epsilon}{\zeta_{S_3,1}}\right)\longrightarrow0,\qquad
W_{25}=O(\Theta_{\zeta}^{-1}).
\]
Since $w_5=O(\epsilon)$, the identity $f_{5,w}^{S_3}=w_5+m_1W_{15}+m_2W_{25}$ proves the second estimate in \eqref{eq:7.22}. Thus \eqref{eq:7.19} applies.

Every exterior inverse cube is $O(\delta_{14})$, while $\delta_{14}\epsilon^3\approx(\Theta_{\zeta}\epsilon^2)(\epsilon/\zeta_{S_3,1})\to0$. The bounded attachment $2$--$5$ and an edge of the isolated $zw$-triangle have incompatible $z/w$ ratios, so the matrix $\Omega$ has rank at most two. Eliminate the fifth coordinate. Its Schur complement $\mathcal{S}$, indexed by $S_4$, has rank at most one and, by Lemma~\ref{lem:5.4},
\[
\mathcal{S}_{11}=m_1m_4\delta_{14}(1+o(1)),\qquad \mathcal{S}_{13}=O(\delta_{14}),\qquad \mathcal{S}_{33}=\mathcal{K}+O(\delta_{14}).
\]
The leading term of $\mathcal{S}_{11}$ is $m_1m_4\delta_{14}$: the other coefficients satisfy $\delta_{13},\delta_{15}=o(\delta_{14})$ and $\delta_{12}\to0$, while the pivot correction is $O(\delta_{14}^2\epsilon^3)=o(\delta_{14})$. The minor $\mathcal{S}_{11}\mathcal{S}_{33}=\mathcal{S}_{13}^2$ implies $\mathcal{K}=O(\delta_{14})=O(\Theta_{\zeta}/\zeta_{S_3,1})$. But $w_{14}=o(\epsilon)$ and the cubic relation in \eqref{eq:7.21} imply $\epsilon^3\zeta_{S_3,1}\Theta_{\zeta}^2\to\infty$, and hence
\[
\frac{\epsilon \Theta_{\zeta}^2}{\Theta_{\zeta}/\zeta_{S_3,1}}=\epsilon \Theta_{\zeta}\zeta_{S_3,1}=\frac{\epsilon^3\zeta_{S_3,1}\Theta_{\zeta}^2}{\epsilon^2\Theta_{\zeta}}\longrightarrow\infty.
\]
This contradicts \eqref{eq:7.19}. The reflection interchanging endpoints proves the same result for $m_3<m_4$.

Sections~\ref{subsec:7.3}--\ref{subsec:7.6} have excluded all comparable-scale cases and every unequal-scale attachment pattern, including the unique possible transverse alignment. They therefore prove Theorem~\ref{thm:7.1} when $m_3\ne m_4$. Together with Section~\ref{sec:diagram5:equal}, this completes the proof for all positive masses.

A supplementary exact symbolic check of the algebraic identities in Section~\ref{sec:diagram5:unequal} is described in Appendix~\ref{app:unequal-static}.

\section{The exclusion of Diagram 8}\label{sec:10}

We use the index sets defined in \eqref{eq:2.index_sets}.

In diagram 8, $T_3$ is the upper collapsing wing and $S_3$ the lower finite maximal wing.
They share vertex 5, with outer pairs $S_2$ and $T_2$, respectively. Its two scales can be retained
as independent analytic coordinates. After eliminating five bounded variables, the first crossed
interaction gives a nonsingular reduced Hessian unless both outer pairs have equal masses. The
remaining case is excluded by combining the upper bound in Proposition~\ref{prop:8.2} with the physical lower bound.

\begin{theorem}\label{thm:10.1}
For arbitrary positive masses, no complex central-configuration component containing
a physical point can admit a singular sequence approaching AK diagram 8 along which some 4-product tends to infinity.
\end{theorem}

The analytic part gives the following stronger restriction without using physical positivity.

\begin{proposition}\label{prop:10.2}
In any complex multiplier-one central-configuration sequence having AK diagram 8,
the two outer masses of the collapsing wing are equal, and the two outer masses of the finite
wing are equal.
\end{proposition}

The masses $m_i>0$, $i\in S_5$, are fixed throughout the degeneration.
For positive masses, the necessary upper-wing mass condition of Albouy--Kaloshin
\cite[Section~8.8, equation~(40)]{AK} can be written as
\begin{equation}
m_5=\frac{(m_1m_2)^{3/2}}{(m_1+m_2)^2}.
\label{eq:10.2}
\end{equation}

We first establish an upper bound for the moment of inertia, valid for arbitrary positive masses compatible with this diagram. Section~\ref{sec:10.5} will apply it after the analytic reduction has restricted the masses to the paired case.

\begin{proposition}\label{prop:8.2}
In AK diagram 8, let the upper outer masses be $m_1,m_2>0$, the lower outer masses $m_3,m_4>0$, and the shared mass $m_5>0$. Define
\begin{equation}
\mathcal{M}(m_3,m_4;m_5)=\frac{\left((m_3m_4)^{3/2}+m_5(m_3+m_4)^2\right)^{2/3}}{(m_3+m_4)^{1/3}}.
\label{eq:8.2}
\end{equation}
The limiting moment of inertia satisfies $|I|\le \mathcal{M}(m_3,m_4;m_5)$. Consequently this diagram cannot occur on a component containing a physical point whenever $\mathcal{M}(m_3,m_4;m_5)<e_2(S_5)/M^{1/3}$.
\end{proposition}

\begin{proof}
AK Section 8.8 gives $z_3,z_4,z_5=O(\epsilon^2)$ and $m_1z_1+m_2z_2=O(\epsilon^2)$. Define the upper-wing $z$ separation scale by
\[
\Theta_{T_3,z}:=\max\limits_{\substack{i<j\\i,j\in T_3}}|z_{ij}|.
\]
The AK estimates give $\epsilon\prec\Theta_{T_3,z}\prec\epsilon^{-2}$, with $z_1,z_2\approx\Theta_{T_3,z}$. Both upper attachments collapse, so $w_i-w_5=o(1/\Theta_{T_3,z})$ for $i\in S_2$. Moreover $w_5=o(\epsilon^{-2})$. It follows that
\[
m_1z_1w_1+m_2z_2w_2+m_5z_5w_5
=\sum\limits_{i\in S_2}m_iz_i(w_i-w_5)+(m_1z_1+m_2z_2+m_5z_5)w_5\longrightarrow0.
\]
The lower vertices are $w$-circled. Total centering and the preceding upper bounds imply $m_3w_3+m_4w_4=o(\epsilon^{-2})$, so $w_3/w_4\to-m_4/m_3$. Their weighted $z$ center is $O(\epsilon^2)$. The exact two-body decomposition of the moment of inertia now gives
\begin{equation}
I-\frac{m_3m_4}{m_3+m_4}r_{34}^2\longrightarrow0.
\label{eq:8.3}
\end{equation}

After passing to a subsequence as $\epsilon\to0$, define the finite nonzero signed-distance limits
\[
r_{34}^0=\lim\limits_{\epsilon\to0}r_{34},\qquad
r_{35}^0=\lim\limits_{\epsilon\to0}r_{35},\qquad
r_{45}^0=\lim\limits_{\epsilon\to0}r_{45}.
\]
The four crossed distances diverge, since their squared moduli are comparable to $\Theta_{T_3,z}\epsilon^{-2}$. The exact lower outer $w$ equations are
\[
\begin{aligned}
w_3&=m_4W_{34}+m_5W_{35}+\sum\limits_{i\in S_2}m_iW_{3i},\\
w_4&=m_3W_{43}+m_5W_{45}+\sum\limits_{i\in S_2}m_iW_{4i}.
\end{aligned}
\]
For $k\in T_2$ and $i\in S_2$, we have $z_{ki}\approx\Theta_{T_3,z}$ and $w_k\approx\epsilon^{-2}$. Thus
\[
\frac{W_{ki}}{w_k}
=\frac{1}{r_{ki}z_{ki}w_k}
=O\!\left(\frac{\epsilon^2}{r_{ki}\Theta_{T_3,z}}\right)
\longrightarrow0,
\qquad k\in T_2,\quad i\in S_2,
\]
because $r_{ki}\to\infty$ and $\epsilon\prec\Theta_{T_3,z}$. Dividing the two equations by $w_3,w_4$, respectively, and using $w_3/w_4\to-m_4/m_3$ and $w_5/w_k\to0$ for $k\in T_2$, therefore gives
\[
1=\frac{m_3+m_4}{(r_{34}^0)^3}+\frac{m_5}{(r_{35}^0)^3}=\frac{m_3+m_4}{(r_{34}^0)^3}+\frac{m_5}{(r_{45}^0)^3},\qquad (r_{35}^0)^3=(r_{45}^0)^3.
\]
The weighted squared-distance identity gives, in the limit,
\[
(r_{34}^0)^2=\frac{m_3+m_4}{m_4}(r_{35}^0)^2+\frac{m_3+m_4}{m_3}(r_{45}^0)^2.
\]
Subtracting the right side from the left before taking limits yields $-(z_3/m_4+z_4/m_3-(1/m_3+1/m_4)z_5)(m_3w_3+m_4w_4-(m_3+m_4)w_5)$, which tends to zero by the coordinate bounds. Since $(r_{35}^0)^3=(r_{45}^0)^3$, their moduli are equal. The squared-distance identity above therefore gives
\[
|r_{34}^0|^2\le\frac{(m_3+m_4)^2}{m_3m_4}|r_{35}^0|^2.
\]
Multiplying the first limiting force equation by $(r_{34}^0)^3$, taking moduli, and using this bound yields
\[
|r_{34}^0|^3=\left|m_3+m_4+m_5\frac{(r_{34}^0)^3}{(r_{35}^0)^3}\right|\le m_3+m_4+\frac{m_5(m_3+m_4)^3}{(m_3m_4)^{3/2}}.
\]
Combining this bound with \eqref{eq:8.3} gives $|I|\le\mathcal{M}(m_3,m_4;m_5)$, with $\mathcal{M}$ defined in \eqref{eq:8.2}. Since $U=I>0$ on a component containing a physical point, comparison with the physical lower bound in Lemma~\ref{lem:3.1} gives the stated exclusion.
\end{proof}

We use the same function $G=U+I/2$ on centered configurations with multiplier one and signed local distance roots.

\subsection{Coordinates and rescaling with \texorpdfstring{$\mathfrak{b}_1,\mathfrak{b}_2$}{b1 and b2}}\label{sec:10.1}

Choose a local square root $\mathfrak{b}_1$ of the upper signed distance $r_{12}\to0$, so that $r_{12}=\mathfrak{b}_1^2$.
The stroke estimates in \cite[Estimate~2]{AK}, applied to the butterfly clustering in \cite[Sections~8.6--8.8, pp.~576--578]{AK}, give
$Z_{12}\approx W_{34}\approx\epsilon^{-2}$.
Since $r_{34}$ has a finite nonzero limit, the identities
$Z_{ij}=z_{ij}/r_{ij}^3$ and $r_{ij}^2=z_{ij}w_{ij}$ give
\[
z_{12}=r_{12}^3Z_{12},\qquad w_{12}=\frac{1}{r_{12}\,Z_{12}},\qquad
w_{34}\approx W_{34},\qquad z_{34}\approx W_{34}^{-1}.
\]
The upper distance satisfies $r_{12}\succ\epsilon$; equivalently,
$(r_{12}^2Z_{12})^{-1}\to0$.
Introduce a second small parameter $\mathfrak{b}_2$ with
$\mathfrak{b}_2\approx (r_{12}^2Z_{12})^{-1}\to0$.
The centers of mass of the two outer pairs are $C_{S_2}$ and $C_{T_2}$, with the notation defined in Section~\ref{sec:2}.
The unscaled positions $q_{S_2,5},q_{T_2,5}$ of body 5 relative to the centers of mass of the pairs retain their global meanings from Section~\ref{sec:2}.
Their hatted components in this section use the scaling specified in \eqref{eq:10.3}, which differs from that in Section~\ref{sec:9}.
Exact coordinates, after a complex rotation, are
\begin{equation}
\begin{gathered}
\begin{aligned}
z_{12}&=(m_1+m_2)r_{12}/\mathfrak{b}_2,& w_{12}&=\mathfrak{b}_2 r_{12}/(m_1+m_2),\\
z_{34}&=\frac{\mathfrak{b}_2 r_{12}^2r_{34}^2}{m_3+m_4},& w_{34}&=(m_3+m_4)/(\mathfrak{b}_2 r_{12}^2),\\
\widehat{\zeta}_{S_2,5}&=\frac{\zeta_{S_2,5}}{\mathfrak{b}_2 r_{12}^2},&
\widehat{\omega}_{S_2,5}&=\frac{\omega_{S_2,5}}{\mathfrak{b}_2 r_{12}},\\
\widehat{\zeta}_{T_2,5}&=\frac{\zeta_{T_2,5}}{\mathfrak{b}_2 r_{12}^2},&
\widehat{\omega}_{T_2,5}&=\frac{r_{12}^2\omega_{T_2,5}}{\mathfrak{b}_2},
\end{aligned}\\
\begin{aligned}
z_5&=\frac{\mathfrak{b}_2 r_{12}^2}{M}
\bigl((m_1+m_2)\widehat{\zeta}_{S_2,5}+(m_3+m_4)\widehat{\zeta}_{T_2,5}\bigr),\\
w_5&=\frac{\mathfrak{b}_2}{M}
\bigl((m_1+m_2)r_{12}\widehat{\omega}_{S_2,5}+(m_3+m_4)r_{12}^{-2}\widehat{\omega}_{T_2,5}\bigr).
\end{aligned}
\end{gathered}
\label{eq:10.3}
\end{equation}
The last two lines determine $q_5=(z_5,w_5)$ from $\widehat{\zeta}_{S_2,5},\widehat{\omega}_{S_2,5},\widehat{\zeta}_{T_2,5},\widehat{\omega}_{T_2,5}$
by imposing that the center of mass of the whole system is at the origin.
In fact the invariant formula
\begin{equation*}
\mathfrak{b}_2^2=\frac{(m_1+m_2)(m_3+m_4)}{r_{12} z_{12}w_{34}}
\end{equation*}
chooses $\mathfrak{b}_2$, and the rotation then makes the first and fourth coordinates of \eqref{eq:10.3} hold.
Its factor and inverse are bounded. The seven parameters $(\mathfrak{b}_1,\mathfrak{b}_2,\widehat{\omega}_{S_2,5},r_{34},\widehat{\zeta}_{T_2,5},\widehat{\zeta}_{S_2,5},\widehat{\omega}_{T_2,5})$ give local
coordinates on the centered configuration space modulo rotations whenever $\mathfrak{b}_1\mathfrak{b}_2\ne0$.
The finite covering choices of the local roots cause no loss of configurations.

We record the bounds needed for compactness of the other five parameters. The lower-cluster estimate in \cite[Section~8.8, p.~578]{AK} gives
$z_j\preceq W_{34}^{-1}$ for $j\in S_3$, and the total center-of-mass equation gives
$C_{S_2,z}\preceq W_{34}^{-1}$. Since
$\mathfrak{b}_2 r_{12}^2\approx W_{34}^{-1}$ by \eqref{eq:10.3}, both $\widehat{\zeta}_{S_2,5}$ and $\widehat{\zeta}_{T_2,5}$ are bounded.
The upper distances are comparable to $r_{12}$, and the upper coordinate estimates give
$w_{i5}\preceq w_{12}$ for $i\in S_2$; this bounds $\widehat{\omega}_{S_2,5}$.
Using the exact identity $W_{12}=\mathfrak{b}_2/((m_1+m_2)r_{12}^2)$, the upper and crossed force estimates become
\[
\begin{aligned}
W_{kl}&\preceq W_{12}, &&k,l\in T_3,\quad k\ne l,\\
W_{ij}&\prec W_{12}, &&i\in S_2,\quad j\in T_2.
\end{aligned}
\]
For the crossed terms, the ratio bound is
$W_{ij}/W_{12}=O(\mathfrak{b}_2 \mathfrak{b}_1^3)\to0$.
The individual upper $w$ force equations therefore give
$w_1,w_2\preceq W_{12}$.
The small upper transverse differences give the same bound for $w_5$,
and the total center-of-mass equation gives $C_{T_2,w}\preceq W_{12}$.
Consequently, the exact identity
\[
\omega_{T_2,5}=(m_1+m_2)W_{12}\widehat{\omega}_{T_2,5}
\]
shows that $\widehat{\omega}_{T_2,5}$ is bounded.
Finally, $r_{34}$ has a finite nonzero limit and is bounded.
The map $r_{34}\mapsto r_{34}^2/(m_3+m_4)^2$ is locally invertible near that limit,
so the signed distance $r_{34}$ itself may be used as a coordinate.
Passing to a subsequence fixes a limiting point and every local signed sheet.
All lower distances and upper distance ratios used below have nonzero finite limits.
The scaled crossed squared distances $\mathfrak b_1^2\mathfrak b_2^2r_{ij}^2$, $i\in S_2$, $j\in T_2$, also have nonzero limits, as displayed in \eqref{eq:10.23}. Choose a fixed polydisc about this limiting point on which these scaled factors, the lower squared distances, and the squared upper distance ratios stay away from zero. Their local signed roots are holomorphic there; the expansions and divisibility argument in Section~\ref{sec:10.4} use this same neighborhood.
Write $r_{34}^0=\lim\limits_{\epsilon\to0}r_{34}\ne0$. A superscript $0$ on any of the four rescaled relative coordinate components likewise denotes its limit along this subsequence as $\epsilon\to0$, along which $(\mathfrak{b}_1,\mathfrak{b}_2)\to(0,0)$.

Rotational invariance of $G$ supplies the omitted rotational equation. Thus its seven critical
equations in \eqref{eq:10.3} characterize the full centered multiplier-one configurations for $\mathfrak{b}_1\mathfrak{b}_2\ne0$.

\subsection{Distance formulas and the expansion of \texorpdfstring{$U_{T_3}$}{U(T3)}}\label{sec:10.2}

The other five wing squared distances are
{\footnotesize\begin{equation}
\begin{aligned}
r_{12}^2&=\mathfrak{b}_1^4,\\
r_{15}^2&=r_{12}^2(m_2-\mathfrak{b}_2^2r_{12}\widehat{\zeta}_{S_2,5})(m_2/(m_1+m_2)^2-\widehat{\omega}_{S_2,5}),\\
r_{25}^2&=r_{12}^2(m_1+\mathfrak{b}_2^2r_{12}\widehat{\zeta}_{S_2,5})(m_1/(m_1+m_2)^2+\widehat{\omega}_{S_2,5}),\\
r_{35}^2&=\left(\frac{m_4r_{34}^2}{(m_3+m_4)^2}-\widehat{\zeta}_{T_2,5}\right)(m_4-\mathfrak{b}_2^2\widehat{\omega}_{T_2,5}),\\
r_{45}^2&=\left(\frac{m_3r_{34}^2}{(m_3+m_4)^2}+\widehat{\zeta}_{T_2,5}\right)(m_3+\mathfrak{b}_2^2\widehat{\omega}_{T_2,5}).
\end{aligned}
\label{eq:10.5}
\end{equation}}
The four crossed squared distances between $S_2$ and $T_2$ are
{\footnotesize\begin{equation}
\begin{aligned}
r_{13}^2 &={}\frac{1}{\mathfrak{b}_2^2r_{12}} \left[m_2+\mathfrak{b}_2^2r_{12}\left(-\widehat{\zeta}_{S_2,5}+\widehat{\zeta}_{T_2,5}-\frac{m_4r_{34}^2}{(m_3+m_4)^2}\right)\right]\\
&\qquad{}\cdot\left[-m_4+\mathfrak{b}_2^2\widehat{\omega}_{T_2,5}+\mathfrak{b}_2^2r_{12}^3\left(-\widehat{\omega}_{S_2,5}+\frac{m_2}{(m_1+m_2)^2}\right)\right],\\[2pt]
r_{14}^2 &={}\frac{1}{\mathfrak{b}_2^2r_{12}} \left[m_2+\mathfrak{b}_2^2r_{12}\left(-\widehat{\zeta}_{S_2,5}+\widehat{\zeta}_{T_2,5}+\frac{m_3r_{34}^2}{(m_3+m_4)^2}\right)\right]\\
&\qquad{}\cdot\left[m_3+\mathfrak{b}_2^2\widehat{\omega}_{T_2,5}+\mathfrak{b}_2^2r_{12}^3\left(-\widehat{\omega}_{S_2,5}+\frac{m_2}{(m_1+m_2)^2}\right)\right],\\[2pt]
r_{23}^2 &={}\frac{1}{\mathfrak{b}_2^2r_{12}} \left[-m_1+\mathfrak{b}_2^2r_{12}\left(-\widehat{\zeta}_{S_2,5}+\widehat{\zeta}_{T_2,5}-\frac{m_4r_{34}^2}{(m_3+m_4)^2}\right)\right]\\
&\qquad{}\cdot\left[-m_4+\mathfrak{b}_2^2\widehat{\omega}_{T_2,5}+\mathfrak{b}_2^2r_{12}^3\left(-\widehat{\omega}_{S_2,5}-\frac{m_1}{(m_1+m_2)^2}\right)\right],\\[2pt]
r_{24}^2 &={}\frac{1}{\mathfrak{b}_2^2r_{12}} \left[-m_1+\mathfrak{b}_2^2r_{12}\left(-\widehat{\zeta}_{S_2,5}+\widehat{\zeta}_{T_2,5}+\frac{m_3r_{34}^2}{(m_3+m_4)^2}\right)\right]\\
&\qquad{}\cdot\left[m_3+\mathfrak{b}_2^2\widehat{\omega}_{T_2,5}+\mathfrak{b}_2^2r_{12}^3\left(-\widehat{\omega}_{S_2,5}-\frac{m_1}{(m_1+m_2)^2}\right)\right].
\end{aligned}
\label{eq:10.6}
\end{equation}}
The moment of inertia is exactly
\begin{equation}
\begin{aligned}
I={}&\frac{m_1m_2}{m_1+m_2}r_{12}^2+\frac{m_3m_4}{m_3+m_4}r_{34}^2\\
&+\mathfrak{b}_2^2\left\{(m_1+m_2)r_{12}^3\widehat{\zeta}_{S_2,5}\widehat{\omega}_{S_2,5}+(m_3+m_4)\widehat{\zeta}_{T_2,5}\widehat{\omega}_{T_2,5}\right.\\
&\left.\hspace{3em}-\frac{\begin{gathered}((m_1+m_2)\widehat{\zeta}_{S_2,5}+(m_3+m_4)\widehat{\zeta}_{T_2,5})\\
{}\cdot((m_1+m_2)r_{12}^3\widehat{\omega}_{S_2,5}+(m_3+m_4)\widehat{\omega}_{T_2,5})\end{gathered}}{M}\right\}.
\end{aligned}
\label{eq:10.7}
\end{equation}
The reciprocal of each crossed distance is $\mathfrak{b}_1\mathfrak{b}_2$ times a holomorphic function that is nonzero near the limiting point.
Every dependence of that factor on the five bounded variables is divisible by $\mathfrak{b}_2^2$.
Thus the derivatives of these reciprocals with respect to those variables are divisible by $\mathfrak{b}_2^3\mathfrak{b}_1$.

Let $\gamma_1(\widehat{\omega}_{S_2,5})$ and $\gamma_2(\widehat{\omega}_{S_2,5})$ be the signed holomorphic roots inherited
from $r_{15}/r_{12}$ and $r_{25}/r_{12}$ at $\mathfrak{b}_2=0$.
The inverse square roots below use the analytic branch with value one when their
perturbation is zero. Equation~\eqref{eq:10.5} gives the exact upper potential
\[
\begin{aligned}
U_{T_3}
&=\frac{m_1m_2}{r_{12}}+\frac{m_1m_5}{r_{15}}+\frac{m_2m_5}{r_{25}}\\
&=\mathfrak{b}_1^{-2}\Biggl[
m_1m_2+
\frac{m_1m_5}{\gamma_1}
\left(1-\frac{\mathfrak{b}_2^2\mathfrak{b}_1^2\widehat{\zeta}_{S_2,5}}{m_2}\right)^{-1/2}
\\
&\qquad{}+\frac{m_2m_5}{\gamma_2}
\left(1+\frac{\mathfrak{b}_2^2\mathfrak{b}_1^2\widehat{\zeta}_{S_2,5}}{m_1}\right)^{-1/2}
\Biggr].
\end{aligned}
\]
Separating the constant and linear terms of these two inverse square roots yields
\begin{equation}
U_{T_3}=\mathfrak{b}_1^{-2}v_1(\widehat{\omega}_{S_2,5})+\mathfrak{b}_2^2\widehat{\zeta}_{S_2,5}v_2(\widehat{\omega}_{S_2,5})
+\mathfrak{b}_2^4\mathfrak{b}_1^2v_3(\mathfrak{b}_1,\mathfrak{b}_2,\widehat{\omega}_{S_2,5},\widehat{\zeta}_{S_2,5}),
\label{eq:10.8}
\end{equation}
where the functions $v_1,v_2,v_3$ are given as follows:
\begin{equation*}
\begin{aligned}
v_1(\widehat{\omega}_{S_2,5})&=m_1m_2+\frac{m_1m_5}{\gamma_1}+\frac{m_2m_5}{\gamma_2},\\[4pt]
v_2(\widehat{\omega}_{S_2,5})&=\frac{m_5}{2}
\left(\frac{m_1}{m_2\gamma_1}-\frac{m_2}{m_1\gamma_2}\right),\\[4pt]
v_3(\mathfrak{b}_1,\mathfrak{b}_2,\widehat{\omega}_{S_2,5},\widehat{\zeta}_{S_2,5})
&=\frac{m_1m_5\widehat{\zeta}_{S_2,5}^2}{m_2^2\gamma_1}
\frac{2+\sqrt{1-\frac{\mathfrak{b}_2^2\mathfrak{b}_1^2\widehat{\zeta}_{S_2,5}}{m_2}}}
{2\sqrt{1-\frac{\mathfrak{b}_2^2\mathfrak{b}_1^2\widehat{\zeta}_{S_2,5}}{m_2}}
\left(1+\sqrt{1-\frac{\mathfrak{b}_2^2\mathfrak{b}_1^2\widehat{\zeta}_{S_2,5}}{m_2}}\right)^2}\\[4pt]
&\quad+\frac{m_2m_5\widehat{\zeta}_{S_2,5}^2}{m_1^2\gamma_2}
\frac{2+\sqrt{1+\frac{\mathfrak{b}_2^2\mathfrak{b}_1^2\widehat{\zeta}_{S_2,5}}{m_1}}}
{2\sqrt{1+\frac{\mathfrak{b}_2^2\mathfrak{b}_1^2\widehat{\zeta}_{S_2,5}}{m_1}}
\left(1+\sqrt{1+\frac{\mathfrak{b}_2^2\mathfrak{b}_1^2\widehat{\zeta}_{S_2,5}}{m_1}}\right)^2}.
\end{aligned}
\end{equation*}
where
\[
\gamma_1^2=m_2(m_2/(m_1+m_2)^2-\widehat{\omega}_{S_2,5}),\qquad
\gamma_2^2=m_1(m_1/(m_1+m_2)^2+\widehat{\omega}_{S_2,5}).
\]
The square roots in the exact remainder $v_3$ take the value one at
$\mathfrak{b}_1=\mathfrak{b}_2=0$.
Its denominators are nonzero near the normalized limiting point, so $v_3$ is
holomorphic in $(\mathfrak{b}_1,\mathfrak{b}_2,\widehat{\omega}_{S_2,5},\widehat{\zeta}_{S_2,5})$ there. In particular,
\[
v_3(0,0,\widehat{\omega}_{S_2,5},\widehat{\zeta}_{S_2,5})
=\frac{3m_5\widehat{\zeta}_{S_2,5}^2}{8}
\left(\frac{m_1}{m_2^2\gamma_1}+\frac{m_2}{m_1^2\gamma_2}\right).
\]
The leading upper $z$ force cancellation gives $\gamma_1^3=\gamma_2^3=-m_5/(m_1+m_2)$.
Together with \eqref{eq:10.2} and the exact identity $m_1\gamma_1^2+m_2\gamma_2^2=m_1m_2/(m_1+m_2)$, this forces
\begin{equation*}
\gamma_1=\gamma_2=-\frac{\sqrt{m_1m_2}}{m_1+m_2},\qquad \widehat{\omega}_{S_2,5}^0=\frac{m_2-m_1}{(m_1+m_2)^2}.
\end{equation*}
Indeed both squared roots have modulus $m_1m_2/(m_1+m_2)^2$. Equality in the triangle inequality for their
positive weighted sum makes both squared roots positive real; their negative cubes then fix their
signed roots. In particular,
\begin{equation*}
\begin{aligned}
v_1(\widehat{\omega}_{S_2,5}^0)&=v_1'(\widehat{\omega}_{S_2,5}^0)=0,\\
v_1''(\widehat{\omega}_{S_2,5}^0)&=\frac{3m_1m_2m_5(m_1+m_2)}{4\gamma_1^5}=-\frac{3(m_1+m_2)^4}{4}\ne0,\\
v_2(\widehat{\omega}_{S_2,5}^0)&=\frac{m_2-m_1}{2}.
\end{aligned}
\end{equation*}

\subsection{Critical points and the Hessian of \texorpdfstring{$G_{S_3}$}{G(S3)}}\label{sec:10.3}

In this subsection, $r_{35}$ and $r_{45}$ denote the signed distance functions in
\eqref{eq:10.5} evaluated at $\mathfrak{b}_2=0$, on their inherited local holomorphic
branches. We denote the lower contribution to $G=U+I/2$ on this chart by
$G_{S_3}(r_{34},\widehat{\zeta}_{T_2,5})$. Its first three terms are the potential terms for the pairs in
$S_3$; its fourth term is one half of $m_3m_4r_{34}^2/(m_3+m_4)$, the second
term in the first line of \eqref{eq:10.7}. Thus
\begin{equation}
\begin{gathered}
G_{S_3}(r_{34},\widehat{\zeta}_{T_2,5})=\frac{m_3m_4}{r_{34}}+\frac{m_3m_5}{r_{35}}+\frac{m_4m_5}{r_{45}}+\frac{m_3m_4}{2(m_3+m_4)}r_{34}^2,\\
r_{35}^2=m_4\left(\frac{m_4r_{34}^2}{(m_3+m_4)^2}-\widehat{\zeta}_{T_2,5}\right),\qquad
r_{45}^2=m_3\left(\frac{m_3r_{34}^2}{(m_3+m_4)^2}+\widehat{\zeta}_{T_2,5}\right).
\end{gathered}
\label{eq:10.12}
\end{equation}
Because $r_{34}\ne0$ near the limiting point, its critical equations give
\begin{equation}
r_{35}^{-3}=r_{45}^{-3},\qquad (m_3+m_4)r_{34}^{-3}+m_5r_{35}^{-3}=1.
\label{eq:10.13}
\end{equation}
Every admissible leading lower point is nondegenerate:
\begin{equation}
\det D_{r_{34},\widehat{\zeta}_{T_2,5}}^2 G_{S_3}=\frac{9m_3^3m_4^3m_5r_{34}^2}{4(m_3+m_4)^2r_{35}^5r_{45}^2}\ne0.
\label{eq:10.14}
\end{equation}
At a critical point, changing from the normalized squared-distance coordinate to $r_{34}$
multiplies the Hessian determinant by $(2r_{34}/(m_3+m_4)^2)^2$.
For an explicit verification in $(r_{34},\widehat{\zeta}_{T_2,5})$, \eqref{eq:10.13} cancels the terms
from the second derivatives of the squared distances and from the moment of inertia.
The Hessian is therefore the sum of the three rank-one matrices formed from
the squared-distance gradients, with coefficients $(3/4)m_3m_4/r_{34}^5$, $(3/4)m_3m_5/r_{35}^5$, and $(3/4)m_4m_5/r_{45}^5$.
Expansion of the determinant gives
\[
\frac{9m_3^2m_4^2m_5r_{34}^2}{4(m_3+m_4)^2}\left\{
\frac{(m_3+m_4)^2}{r_{34}^5}\left(\frac{m_4}{r_{35}^5}+\frac{m_3}{r_{45}^5}\right)+\frac{m_3m_4m_5}{r_{35}^5r_{45}^5}
\right\}.
\]
The identity $m_3r_{35}^2+m_4r_{45}^2=(m_3m_4/(m_3+m_4))r_{34}^2$ and \eqref{eq:10.13} simplify it to \eqref{eq:10.14}.
Since all three lower distances are finite and nonzero,
\eqref{eq:10.14} proves nondegeneracy on every signed or nonreal lower distance
branch admitting a critical point.

\subsection{Reduction of \texorpdfstring{$G$}{G} to \texorpdfstring{$\tilde G(\mathfrak{b}_1,\mathfrak{b}_2)$}{Gtilde(b1,b2)}}\label{sec:10.4}

We treat $\mathfrak b_1,\mathfrak b_2$ as independent complex parameters in the fixed signed neighborhood of Section~\ref{sec:10.1}.

In the coordinates $(\mathfrak{b}_1,\mathfrak{b}_2,\widehat{\omega}_{S_2,5},r_{34},\widehat{\zeta}_{T_2,5},\widehat{\zeta}_{S_2,5},\widehat{\omega}_{T_2,5})$,
combining \eqref{eq:10.7}, \eqref{eq:10.8}, and \eqref{eq:10.12} gives the exact
expression for the full function $G=U+I/2$:
\[
\begin{aligned}
G={}&\mathfrak{b}_1^{-2}v_1(\widehat{\omega}_{S_2,5})
+\mathfrak{b}_2^2\widehat{\zeta}_{S_2,5}v_2(\widehat{\omega}_{S_2,5})
+\mathfrak{b}_2^4\mathfrak{b}_1^2v_3(\mathfrak{b}_1,\mathfrak{b}_2,\widehat{\omega}_{S_2,5},\widehat{\zeta}_{S_2,5})\\
&+G_{S_3}(r_{34},\widehat{\zeta}_{T_2,5})
+\frac{m_1m_2}{2(m_1+m_2)}\mathfrak{b}_1^4\\
&+\frac{m_3m_5}{r_{35}}
\left[\left(1-\frac{\mathfrak{b}_2^2\widehat{\omega}_{T_2,5}}{m_4}\right)^{-1/2}-1\right]\\
&+\frac{m_4m_5}{r_{45}}
\left[\left(1+\frac{\mathfrak{b}_2^2\widehat{\omega}_{T_2,5}}{m_3}\right)^{-1/2}-1\right]
+\sum\limits_{\substack{i\in S_2\\j\in T_2}}\frac{m_im_j}{r_{ij}}\\
&+\frac{\mathfrak{b}_2^2}{2M}\Big[
(m_1+m_2)(m_3+m_4+m_5)\mathfrak{b}_1^6\widehat{\zeta}_{S_2,5}\widehat{\omega}_{S_2,5}\\
&\qquad +(m_3+m_4)(m_1+m_2+m_5)\widehat{\zeta}_{T_2,5}\widehat{\omega}_{T_2,5}\\
&\qquad -(m_1+m_2)(m_3+m_4)
(\widehat{\zeta}_{S_2,5}\widehat{\omega}_{T_2,5}+\mathfrak{b}_1^6\widehat{\zeta}_{T_2,5}\widehat{\omega}_{S_2,5})\Big].
\end{aligned}
\]
Here $r_{35}$ and $r_{45}$ are the signed distance functions of \eqref{eq:10.12}
at $\mathfrak{b}_2=0$, still depending on $(r_{34},\widehat{\zeta}_{T_2,5})$.
The crossed distances $r_{ij}$, $i\in S_2$, $j\in T_2$, use the full expressions
in \eqref{eq:10.6} with $r_{12}=\mathfrak{b}_1^2$.
The two inverse square roots take the value one at $\mathfrak{b}_2=0$, preserving
the inherited signed branches.

Order the five bounded variables as $(\widehat{\omega}_{S_2,5},r_{34},\widehat{\zeta}_{T_2,5},\widehat{\zeta}_{S_2,5},\widehat{\omega}_{T_2,5})$. The five equations
\begin{equation}
\mathfrak{b}_1^2\partial_{\widehat{\omega}_{S_2,5}}G=0,\qquad \partial_{r_{34}}G=\partial_{\widehat{\zeta}_{T_2,5}}G=0,\qquad
\mathfrak{b}_2^{-2}\partial_{\widehat{\zeta}_{S_2,5}}G=\mathfrak{b}_2^{-2}\partial_{\widehat{\omega}_{T_2,5}}G=0
\label{eq:10.15}
\end{equation}
extend holomorphically across $(\mathfrak{b}_1,\mathfrak{b}_2)=(0,0)$. The apparent divisions in the last two equations
are removable: the upper $\widehat{\zeta}_{S_2,5}$ dependence, the lower $\widehat{\omega}_{T_2,5}$ dependence, and the contribution of the centers to the moment of inertia
contain $\mathfrak{b}_2^2$; crossed derivatives contain $\mathfrak{b}_2^3\mathfrak{b}_1$. The Laurent term $\mathfrak{b}_1^{-2}v_1(\widehat{\omega}_{S_2,5})$ is
independent of the other four bounded variables.

For $i,j\ge0$, write $G_{\mathfrak{b}}^{i,j}$ for the coefficient
$[\mathfrak{b}_1^i\mathfrak{b}_2^j](G-\mathfrak{b}_1^{-2}v_1(\widehat{\omega}_{S_2,5}))$, with all five bounded coordinates held fixed.
The subscript $\mathfrak{b}$ identifies the ordered parameter pair $(\mathfrak{b}_1,\mathfrak{b}_2)$; the first and second superscript indices give the powers of $\mathfrak{b}_1$ and $\mathfrak{b}_2$, respectively.
In the following formula, $G_{\mathfrak{b}}^{0,2}(\widehat{\zeta}_{S_2,5},\widehat{\omega}_{T_2,5})$ denotes this coefficient evaluated at
$(\widehat{\omega}_{S_2,5},r_{34},\widehat{\zeta}_{T_2,5})=(\widehat{\omega}_{S_2,5}^0,r_{34}^0,\widehat{\zeta}_{T_2,5}^0)$.
Thus
\begin{equation*}
\begin{aligned}
G_{\mathfrak{b}}^{0,2}(\widehat{\zeta}_{S_2,5},\widehat{\omega}_{T_2,5})={}&\frac{m_2-m_1}{2}\widehat{\zeta}_{S_2,5}
-\frac{(m_1+m_2)(m_3+m_4)}{2M}\widehat{\zeta}_{S_2,5}\widehat{\omega}_{T_2,5}\\
&+\frac{m_3+m_4}{2}
\left(\frac{m_1+m_2+m_5}{M}-m_5r_{35}^{-3}\right)\widehat{\zeta}_{T_2,5}^0\widehat{\omega}_{T_2,5}.
\end{aligned}
\end{equation*}
Here $r_{35}$ is evaluated at $(r_{34}^0,\widehat{\zeta}_{T_2,5}^0)$ on the chart $\mathfrak{b}_2=0$.
As in Section~\ref{sec:9}, this formula gives the values and the derivatives with respect to the last two bounded variables on that locus. Derivatives with respect to the first three bounded variables are taken from the exact distance expressions before specialization. The last two equations in \eqref{eq:10.15} give
\begin{equation}
\begin{aligned}
\widehat{\omega}_{T_2,5}^0&=\frac{M(m_2-m_1)}{(m_1+m_2)(m_3+m_4)},\\
\widehat{\zeta}_{S_2,5}^0&=\frac{m_1+m_2+m_5-Mm_5r_{35}^{-3}}{m_1+m_2}\widehat{\zeta}_{T_2,5}^0.
\end{aligned}
\label{eq:10.17}
\end{equation}
The Jacobian of these five equations in the bounded variables is block lower triangular. The block corresponding to the first three variables consists of $v_1''(\widehat{\omega}_{S_2,5}^0)$ and the invertible lower Hessian \eqref{eq:10.14}; the block corresponding to the last two variables is
\[
-\frac{(m_1+m_2)(m_3+m_4)}{2M}
\begin{pmatrix}0&1\\1&0\end{pmatrix}.
\]
This block is invertible because all masses are positive. Hence the analytic implicit-function
theorem solves all five variables uniquely as
holomorphic functions of $(\mathfrak{b}_1,\mathfrak{b}_2)$ near the limiting point. Every approaching full
central-configuration sequence must lie on this solution.

Differentiating the exact expression for $G$ gives the following explicit form
of the first equation in \eqref{eq:10.15}:
\begin{equation}
\begin{aligned}
0={}&v_1'(\widehat{\omega}_{S_2,5})
-\frac{m_1m_2m_5}{2\gamma_2^3}
\left[\left(1+\frac{\mathfrak{b}_2^2\mathfrak{b}_1^2\widehat{\zeta}_{S_2,5}}{m_1}\right)^{-1/2}-1\right]\\
&+\frac{m_1m_2m_5}{2\gamma_1^3}
\left[\left(1-\frac{\mathfrak{b}_2^2\mathfrak{b}_1^2\widehat{\zeta}_{S_2,5}}{m_2}\right)^{-1/2}-1\right]\\
&+\frac{(m_1+m_2)\mathfrak{b}_2^2\mathfrak{b}_1^8}{2M}
\big[(m_3+m_4+m_5)\widehat{\zeta}_{S_2,5}-(m_3+m_4)\widehat{\zeta}_{T_2,5}\big]\\
&+\frac{\mathfrak{b}_2\mathfrak{b}_1^4}{2}
\sum\limits_{\substack{i\in S_2\\j\in T_2}}m_im_jZ_{ij}.
\end{aligned}
\label{eq:10.18}
\end{equation}
Here $v_1'(\widehat{\omega}_{S_2,5})=\frac{m_1m_2m_5}{2}(\gamma_1^{-3}-\gamma_2^{-3})$.
For the crossed pairs, \eqref{eq:10.3} gives
$\partial_{\widehat{\omega}_{S_2,5}}r_{ij}^2=-\mathfrak{b}_2\mathfrak{b}_1^2z_{ij}$, which yields
the last term. The two inverse-square-root differences are divisible by
$\mathfrak{b}_2^2\mathfrak{b}_1^2$. By \eqref{eq:10.3} and \eqref{eq:10.6},
each crossed $Z_{ij}=z_{ij}/r_{ij}^3$ is $\mathfrak{b}_2^2\mathfrak{b}_1^5$
times a holomorphic function. Hence every correction term after $v_1'(\widehat{\omega}_{S_2,5})$
in \eqref{eq:10.18} is holomorphic and divisible by
$\mathfrak{b}_2^2\mathfrak{b}_1^2$.
Since $v_1''(\widehat{\omega}_{S_2,5}^0)\ne0$, the quotient
\begin{equation}
\frac{\widehat{\omega}_{S_2,5}(\mathfrak{b}_1,\mathfrak{b}_2)-\widehat{\omega}_{S_2,5}^0}{\mathfrak{b}_2^2\mathfrak{b}_1^2}
\qquad (\mathfrak{b}_1\mathfrak{b}_2\ne0)
\label{eq:10.19}
\end{equation}
extends holomorphically to a neighborhood of $(0,0)$.
Together with $v_1(\widehat{\omega}_{S_2,5}^0)=v_1'(\widehat{\omega}_{S_2,5}^0)=0$, this proves that the quotient
\[
\mathfrak{b}_1^{-2}v_1(\widehat{\omega}_{S_2,5}(\mathfrak{b}_1,\mathfrak{b}_2))
\]
is holomorphic and divisible by $\mathfrak{b}_2^4\mathfrak{b}_1^2$. Thus the reduced function
\begin{equation*}
\tilde{G}(\mathfrak{b}_1,\mathfrak{b}_2)=G\bigl(\mathfrak{b}_1,\mathfrak{b}_2,\text{the five solved variables}\bigr)
\end{equation*}
extends holomorphically to a neighborhood of the origin.

At $\mathfrak{b}_2=0$ the first three solved variables are exactly $(\widehat{\omega}_{S_2,5}^0,r_{34}^0,\widehat{\zeta}_{T_2,5}^0)$, for every small $\mathfrak{b}_1$. Therefore
\begin{equation*}
\tilde{G}(\mathfrak{b}_1,0)=G_{S_3}(r_{34}^0,\widehat{\zeta}_{T_2,5}^0)+\frac{m_1m_2}{2(m_1+m_2)}\mathfrak{b}_1^4.
\end{equation*}
In particular $\partial_{\mathfrak{b}_1}^2\tilde{G}(0,0)=0$. The variables $(r_{34},\widehat{\zeta}_{T_2,5})$ differ from their limiting values by $O(\mathfrak{b}_2^2)$, because the crossed leading coefficient is independent of all five bounded variables. These shifts preserve the quadratic terms of $G$ by stationarity of the limiting lower function. Shifts in the last two bounded variables also preserve those terms, since their limiting values in \eqref{eq:10.17} are critical for $G_{\mathfrak{b}}^{0,2}$. After substitution, the upper singular term starts at total degree six by \eqref{eq:10.19}.
On the inherited signed distance branches, define the scaled limits
\[
\chi_{ij}:=\lim\limits_{(\mathfrak{b}_1,\mathfrak{b}_2)\to(0,0)}
\mathfrak{b}_1\mathfrak{b}_2r_{ij},\qquad i\in S_2,\ j\in T_2.
\]
For $i,j\ge0$, define $\tilde{G}_{\mathfrak{b}}^{i,j}:=[\mathfrak{b}_1^i\mathfrak{b}_2^j]\tilde{G}$ in its Taylor expansion at $(0,0)$.
In particular, $\tilde{G}_{\mathfrak{b}}^{1,1}$ is the coefficient of $\mathfrak{b}_1\mathfrak{b}_2$.
The preceding stationarity and divisibility show that it equals the corresponding
coefficient in the holomorphic part $G-\mathfrak{b}_1^{-2}v_1(\widehat{\omega}_{S_2,5})$,
which is independent of the five bounded variables.
Consequently the quadratic expansion is
\begin{equation}
\tilde{G}(\mathfrak{b}_1,\mathfrak{b}_2)=G_{S_3}(r_{34}^0,\widehat{\zeta}_{T_2,5}^0)+\tilde{G}_{\mathfrak{b}}^{0,2}\mathfrak{b}_2^2+\tilde{G}_{\mathfrak{b}}^{1,1}\mathfrak{b}_1\mathfrak{b}_2
+O\bigl((|\mathfrak{b}_1|+|\mathfrak{b}_2|)^3\bigr),
\label{eq:10.22}
\end{equation}
where $\tilde{G}_{\mathfrak{b}}^{0,2}:=G_{\mathfrak{b}}^{0,2}(\widehat{\zeta}_{S_2,5}^0,\widehat{\omega}_{T_2,5}^0)$, and
\begin{equation}
\begin{gathered}
\tilde{G}_{\mathfrak{b}}^{1,1}=\frac{m_1m_3}{\chi_{13}}+\frac{m_1m_4}{\chi_{14}}+\frac{m_2m_3}{\chi_{23}}+\frac{m_2m_4}{\chi_{24}},\\
\chi_{13}^2=-m_2m_4,\qquad \chi_{14}^2=m_2m_3,\qquad
\chi_{23}^2=m_1m_4,\qquad \chi_{24}^2=-m_1m_3.
\end{gathered}
\label{eq:10.23}
\end{equation}
All independent signed choices of these four roots are retained. For $\mathfrak{b}_1\mathfrak{b}_2\ne0$ the two
remaining critical equations are exactly $\partial_{\mathfrak{b}_1}\tilde{G}=\partial_{\mathfrak{b}_2}\tilde{G}=0$, because all five
eliminated derivatives vanish there.

\subsection{The cases when \texorpdfstring{$\tilde G_{\mathfrak{b}}^{1,1}\ne0$}{Gtilde(b;1,1) is nonzero} and \texorpdfstring{$\tilde G_{\mathfrak{b}}^{1,1}=0$}{Gtilde(b;1,1) is zero}}\label{sec:10.5}

If $\tilde{G}_{\mathfrak{b}}^{1,1}\ne0$, \eqref{eq:10.22} gives
\begin{equation*}
\det D^2\tilde{G}(0,0)=-\bigl(\tilde{G}_{\mathfrak{b}}^{1,1}\bigr)^2\ne0.
\end{equation*}
The analytic inverse-function theorem applied to $\nabla_{(\mathfrak{b}_1,\mathfrak{b}_2)}\tilde{G}$ implies that its only zero sufficiently
near the origin is the origin itself. No approaching diagram-8 sequence can occur on that sheet.

For positive masses, the two diagonal terms in \eqref{eq:10.23} are purely imaginary and the two
off-diagonal terms are real. If $\tilde{G}_{\mathfrak{b}}^{1,1}=0$, cancellation in the imaginary part forces
\[
\frac{m_1m_3}{\sqrt{m_2m_4}}=\frac{m_2m_4}{\sqrt{m_1m_3}},\qquad\text{hence }m_1m_3=m_2m_4.
\]
Cancellation in the real part similarly forces $m_1m_4=m_2m_3$. Positivity then gives $m_1=m_2$ and $m_3=m_4$.
This proves Proposition~\ref{prop:10.2}; it makes no assertion excluding complex ends in that residual paired case.

For a component containing a physical point, the residual masses satisfy $m_2=m_1$, $m_4=m_3$, and $m_5=m_1/4$ by \eqref{eq:10.2}.
Since $U=I>0$ on this component, Proposition~\ref{prop:8.2}, specialized to these masses, gives
\[
U\le\mathcal{M}(m_3,m_3;m_1/4)=2m_3\left(\frac{m_1+m_3}{4}\right)^{2/3},\qquad
U^3\le\frac{m_3^3(m_1+m_3)^2}{2}.
\]
Here the mass sum and pair sum are
\[
M=\frac{9m_1+8m_3}{4},\qquad e_2(S_5)=\frac{3m_1^2+9m_1m_3+2m_3^2}{2}.
\]
The strict comparison with the physical lower bound follows from the exact identity
\begin{equation*}
\begin{aligned}
8\left(e_2(S_5)^3-\frac{M m_3^3(m_1+m_3)^2}{2}\right)={}&27m_1^6+243m_1^5m_3+783m_1^4m_3^2\\
&+1044m_1^3m_3^3+496m_1^2m_3^4+83m_1m_3^5>0.
\end{aligned}
\end{equation*}
Thus $U<e_2(S_5)/M^{1/3}$, contradicting the physical lower bound. This proves Theorem~\ref{thm:10.1}.

\section{Completion of the five-body argument}
\label{sec:12}

\begin{theorem}\label{thm:12.1}
The conclusion of Theorem~\ref{thm:1.1} holds for every positive mass vector for the planar five-body problem.
\end{theorem}
\begin{proof}
If a physical point were not isolated on the normalized rotation slice, a positive-dimensional
irreducible component through it would, by Lemma~\ref{lem:2.3}, contain a singular sequence along which some 4-product tends to infinity. Proposition~\ref{prop:3.3}
excludes diagram 3, Proposition~\ref{prop:4.1} excludes diagrams 1 and 4, Proposition~\ref{prop:8.1} excludes
diagram 7, Theorem~\ref{thm:9.1} excludes diagram 6, Theorem~\ref{thm:11.1} excludes diagram 2,
Theorem~\ref{thm:7.1} excludes diagram 5, and Theorem~\ref{thm:10.1} excludes diagram 8. These exhaust the possible diagrams for this singular sequence, a
contradiction. An affine algebraic variety has finitely many irreducible components, hence finitely
many isolated points. Every physical similarity class has a representative in the fixed rotation gauge,
with only finite remaining sign and reflection choices. There are therefore finitely many physical similarity classes.
\end{proof}

At masses $(1,1,1,1,1/4)$, the signed Roberts family has $U^3=32$ under multiplier-one normalization, while the cube of the physical lower bound is $1372/17>32$. By Lemmas~\ref{lem:2.1} and~\ref{lem:3.1}, a component containing this family has no physical point.

\begin{corollary}\label{cor:12.2}
The number of physical similarity classes is at most $4^{17}$ for every positive mass vector for the planar
five-body problem.
\end{corollary}
\begin{proof}
Eliminate the distance variables in favor of inverses $s_{ij}$. The ten relations
$s_{ij}^2((x_i-x_j)^2+(y_i-y_j)^2)=1$ and the force equations have degree at most four. Centering and the linear rotation gauge
remove three position coordinates, leaving seventeen variables. The affine B\'ezout bound for isolated
zeros gives at most $4^{17}$ isolated points, even if other components have positive dimension \cite[Example~8.4.6]{Fulton}. One
justification is to take seventeen generic linear combinations of the equations, preserving isolation
at each isolated zero, and perturb these combinations generically. Each isolated zero contributes
at least one nearby zero, while ordinary B\'ezout bounds the perturbed total by $4^{17}$. Every physical
class has a representative in the gauge.
\end{proof}

Appendix~\ref{app:energy} establishes the polynomial factorization and coefficient positivity used in
the three-point potential obstruction. Appendix~\ref{app:unequal-static} records supplementary exact
checks of the unequal-endpoint identities and expansion coefficients. The asymptotic case
classification, analytic chart coverage, and global algebraic selection are proved in the main text.

\appendix
\section{Algebraic verification for Section~\ref*{subsec:6.3}}\label{app:energy}

This appendix proves the elimination and positivity statements used in the collinear part of
Lemma~\ref{lem:6.2}. Section~\ref{subsec:6.3} derives the common-root condition from the three-point
equations and applies Proposition~\ref{prop:app-energy} below, including the continuity argument
at the exceptional mass values.

Fix $(\sigma_{10},\sigma_{20})\in\{1,-1\}\times\{1,-1\}$ and take as inputs the polynomials
$g_5,g_2,g_3,g_6$ defined in \eqref{eq:6.10}, \eqref{eq:6.11}, \eqref{eq:6.12}, and
\eqref{eq:6.14}, respectively.
Throughout the symbolic calculation, $m_1,m_2,\Delta_U$ are independent variables, and
\[
M=m_1+m_2+\frac94,\qquad
e_2(S_5)=m_1m_2+\frac94(m_1+m_2)+\frac32.
\]
The mass polynomial in the elimination statement is
\[
\mathcal{G}_{\sigma_{10},\sigma_{20}}(m_1,m_2)
=(-\sigma_{10}m_1+\sigma_{20}m_2-4m_1m_2)^2+4\sigma_{10}\sigma_{20}m_1m_2.
\]
The interpretation $\Delta_U=(U^0)^3-e_2(S_5)^3/M$ is imposed when applying the result to a configuration.

\begin{proposition}\label{prop:app-energy}
For each sign combination there is a polynomial
$\mathcal{D}_{\sigma_{10},\sigma_{20}}\in\mathbb Z[m_1,m_2,\Delta_U]$, specified by the construction below, with the following properties.
\begin{enumerate}[label=(\roman*)]
\item For $m_1,m_2>0$ and any $\Delta_U\in\mathbb C$, if $g_5$ and $g_6$ have a common root in $\mathbb C$, then
\[
\mathcal{G}_{\sigma_{10},\sigma_{20}}(m_1,m_2)
\mathcal{D}_{\sigma_{10},\sigma_{20}}(m_1,m_2,\Delta_U)=0.
\]
\item Every nonzero coefficient of $\mathcal{D}_{\sigma_{10},\sigma_{20}}$ is a positive integer, and its constant term with respect to $\Delta_U$ is a nonzero polynomial in $m_1,m_2$. Consequently,
\[
\mathcal{D}_{\sigma_{10},\sigma_{20}}(m_1,m_2,\Delta_U)>0
\quad\text{for }m_1,m_2>0,\ \Delta_U\ge0.
\]
\end{enumerate}
\end{proposition}

\subsection{Polynomial elimination for Lemma~\ref*{lem:6.2}}\label{app:energy:construction}

Introduce the coefficient expansions
\begin{equation*}
g_5(\xi)=\sum\limits_{i=0}^5g_{5i}\xi^i,\qquad
g_6(\xi)=\sum\limits_{i=0}^6g_{6i}\xi^i,\qquad
g_{55}=-36(m_1+m_2).
\end{equation*}
All these coefficients belong to $\mathbb Z[m_1,m_2,\Delta_U]$; the factors in \eqref{eq:6.14} clear the denominators in $M$ and $e_2(S_5)$.
The pseudo-remainder can be written directly as
\begin{equation}
\begin{aligned}
g_4(\xi)&=g_{55}^2g_6(\xi)-\bigl(g_{55}g_{66}\xi+g_{55}g_{65}-g_{66}g_{54}\bigr)g_5(\xi),\\
\deg_{\xi} g_4&\le4.
\end{aligned}
\label{eq:A.2}
\end{equation}
Indeed, the two highest-degree terms cancel. Thus $g_4$ is the unique polynomial of degree at most four satisfying
\[
g_4(\xi)\equiv g_{55}^2g_6(\xi)\pmod{g_5(\xi)}.
\]
Formula~\eqref{eq:A.2} also shows directly that $g_4$ has integer-polynomial coefficients.
For positive masses, $g_{55}\ne0$ and $g_4(\xi)=g_{55}^2g_6(\xi)$ at every root of $g_5$; hence $g_5,g_6$ have a common root if and only if $g_5,g_4$ do.

Write $g_4(\xi)=\sum\limits_{i=0}^4g_{4i}\xi^i$ and set $g_{45}=0$.
For independent formal variables $\xi_1,\xi_2$, define a $5\times5$ matrix $\mathcal{B}$ by
\begin{equation}
\sum\limits_{i,j=0}^4\mathcal{B}_{ij}\xi_1^i\xi_2^j
=\frac{g_5(\xi_1)g_4(\xi_2)-g_4(\xi_1)g_5(\xi_2)}{\xi_1-\xi_2}.
\label{eq:A.3}
\end{equation}
The numerator vanishes when $\xi_1=\xi_2$, so the division is exact.
The quotient has degree at most four in each variable, and its coefficient of $\xi_1^i\xi_2^j$ is the matrix entry $\mathcal{B}_{ij}$.
Interchanging $\xi_1$ and $\xi_2$ leaves the quotient unchanged, proving that $\mathcal{B}$ is symmetric.
Its entries belong to $\mathbb Z[m_1,m_2,\Delta_U]$.
Equivalently, starting from the zero matrix, for each $0\le j<i\le5$ add $g_{5i}g_{4j}-g_{4i}g_{5j}$ to the entries
$(i-1-k,j+k)$ for $k=0,\ldots,i-j-1$.

If $\check{\xi}$ is a common root of $g_5$ and $g_4$, substituting $\xi_2=\check{\xi}$ into the polynomial identity \eqref{eq:A.3} gives
\[
\mathcal{B}(1,\check{\xi},\check{\xi}^2,\check{\xi}^3,\check{\xi}^4)^T=0.
\]
The displayed vector is nonzero, so $\det\mathcal{B}=0$.
This is the elimination step that removes the shape variable from the common-root condition.

\subsection{Factorization and positivity for Lemma~\ref*{lem:6.2}}\label{app:energy:positivity}

For each of the four sign combinations, exact polynomial division gives the factorization
\begin{equation}
\begin{aligned}
\det\mathcal{B}
={}&-\sigma_{10}\sigma_{20}\,2^{35}3^{18}(m_1+m_2)^9M^4\\
&{}\times\mathcal{G}_{\sigma_{10},\sigma_{20}}(m_1,m_2)
\mathcal{D}_{\sigma_{10},\sigma_{20}}(m_1,m_2,\Delta_U).
\end{aligned}
\label{eq:A.4}
\end{equation}
The polynomial $\mathcal{D}_{\sigma_{10},\sigma_{20}}$ is defined as the exact quotient in this identity, with the displayed normalization.
The checks described in Section~\ref{app:energy:verification} verify that this quotient belongs to $\mathbb Z[m_1,m_2,\Delta_U]$ and that multiplication by the entire displayed factor reproduces $\det\mathcal{B}$ coefficient by coefficient.
Because the quotient is constructed as a polynomial, it is also defined at parameter values where $\mathcal{G}_{\sigma_{10},\sigma_{20}}=0$.

Every nonzero coefficient of the integer polynomial $\mathcal{D}_{\sigma_{10},\sigma_{20}}$ is strictly positive.
For each of the four sign combinations $(\sigma_{10},\sigma_{20})\in\{1,-1\}\times\{1,-1\}$, the corresponding polynomial $\mathcal{D}_{\sigma_{10},\sigma_{20}}(m_1,m_2,\Delta_U)$ has exactly 952 nonzero monomials.
In every case, the numbers of monomials at successive degrees $0,1,\ldots,5$ in $\Delta_U$ are respectively
\begin{equation}
302,\qquad234,\qquad174,\qquad122,\qquad78,\qquad42.
\label{eq:A.5}
\end{equation}
These assertions are verified by inspecting the complete exact coefficient lists.
In particular, the constant term with respect to $\Delta_U$ is a nonzero polynomial with positive coefficients in $m_1,m_2$.
It is strictly positive for positive masses, and all the terms containing positive powers of $\Delta_U$ are nonnegative when $\Delta_U\ge0$.
This proves the strict positivity required in Proposition~\ref{prop:app-energy}(ii), including at $\Delta_U=0$.

\begin{proof}[Proof of Proposition~\ref{prop:app-energy}]
For positive masses, Section~\ref{app:energy:construction} shows that a common root of $g_5,g_6$ forces $\det\mathcal{B}=0$.
The factor $-\sigma_{10}\sigma_{20}2^{35}3^{18}(m_1+m_2)^9M^4$ in \eqref{eq:A.4} is nonzero, giving part~(i) for every complex $\Delta_U$.
The exact coefficient verification and the preceding argument for the constant term with respect to $\Delta_U$ give part~(ii).
The treatment of the exceptional mass values with $\mathcal{G}_{\sigma_{10},\sigma_{20}}=0$ remains the continuity argument in the proof of Lemma~\ref{lem:6.2}.
\end{proof}

\subsection{Exact coefficient verification for Lemma~\ref*{lem:6.2}}\label{app:energy:verification}

The factorization \eqref{eq:A.4} and the coefficient signs are verified by computer-assisted exact
polynomial arithmetic in $\mathbb Z[m_1,m_2,\Delta_U]$, with $m_1,m_2,\Delta_U$ treated as symbolic indeterminates, for all four sign combinations. The determinant is
evaluated in two ways: a division-free Laplace expansion using the coefficient construction of
Section~\ref{app:energy:construction}, and fraction-free elimination after reconstructing the
matrix from \eqref{eq:A.3}. In each calculation, exact polynomial division yields the quotient
$\mathcal{D}_{\sigma_{10},\sigma_{20}}$, and multiplication by the full prefactor verifies
\eqref{eq:A.4} coefficient by coefficient.

The complete coefficient lists and the source code for both calculations are supplied as
supplementary material, together with the software requirements and instructions for reproducing
the verification. Each list specifies the coefficients in
\[
\mathcal{D}_{\sigma_{10},\sigma_{20}}
=\sum_{i,j,k\ge0}c_{ijk}^{\sigma_{10},\sigma_{20}}m_1^im_2^j\Delta_U^k.
\]
The exact verification compares every coefficient with the reconstructed quotient, checks that
each nonzero coefficient is a positive integer, and confirms the counts in \eqref{eq:A.5}.
Thus the positivity assertion uses all coefficients, including those of degree zero in $\Delta_U$;
the monomial counts serve as consistency checks.

\section{Supplementary symbolic verification for Section~\ref*{sec:diagram5:unequal}}\label{app:unequal-static}

For the unequal-endpoint case of Section~\ref{sec:diagram5:unequal}, supplementary exact
calculations verify the centered-coordinate formulas in \eqref{eq:7.3}, the alignment identity
\eqref{eq:7.12}, the Taylor coefficients of $\hat f_{5,z},\hat f_{5,w}$ and the inverse derivative
in \eqref{eq:7.13}--\eqref{eq:7.15}, and the Schur identity and quadratic term in
\eqref{eq:7.17}--\eqref{eq:7.19}.

Under the normalization \eqref{eq:7.1}, we have $\sqrt{m_4}=1-\sqrt{m_3}$, so the coefficients
are rational functions of $\sqrt{m_3}$. The identities and retained Taylor coefficients are
checked exactly over $\mathbb Q(\sqrt{m_3})$, on the analytic branches fixed in
Section~\ref{sec:7.5}. Here $\sqrt{m_3}$ is treated as a symbolic parameter.
For the Taylor calculations, all coefficients through total degree
two in $(\psi_z,\psi_w)$ are retained exactly. These comparisons use exact rational arithmetic
throughout, without numerical evaluation of the masses or radicals.

These calculations verify the algebraic identities and Taylor coefficients. Section~\ref{sec:diagram5:unequal} establishes the analytic inversion, remainder estimates, and classification of limiting cases.
The source code and instructions for reproducing these checks are supplied as supplementary material.

\section*{Acknowledgements}
\begin{sloppypar}
The first named author is partially supported by the National Key Research and Development
Program of China (NKPs) [Grant Number 2024YFA1013201], the Natural Science Foundation of
Chongqing (NSFCQ) [Grant Number \mbox{CSTB2024NSCQ-LZX0040}], and the Special Project of Chongqing
Municipal Science and Technology Bureau [Grant Number 2025CCZ015].

\par\bigskip
\noindent\textbf{Declaration of AI use}\par\nobreak\smallskip

The initial draft of this manuscript was prepared around September 20, 2026, with the assistance of ChatGPT 6.0 (OpenAI), which also helped develop the potential estimates and diagram-exclusion arguments in Sections 3–9 and the symbolic verification programs for Appendices A and B, and assisted with revision and exposition. The authors take full responsibility for the content, including the correctness of the mathematical arguments and the accuracy and completeness of the citations. 

\end{sloppypar}

\end{document}